\documentclass[a4paper,twoside,10pt]{amsart}

       \usepackage{color}

 \usepackage[active]{srcltx}

\usepackage{amsthm,amsmath,amssymb,amscd}
\usepackage{graphicx}
\usepackage{epic}

\newtheorem{theorem}{Theorem}[section]
\newtheorem{proposition}[theorem]{Proposition}
\newtheorem{lemma}[theorem]{Lemma}
\newtheorem{corollary}[theorem]{Corollary}

\theoremstyle{definition}
\newtheorem{definition}[theorem]{Definition}
\newtheorem{notation}[theorem]{Notation}

\newtheorem{example}[theorem]{Example}
\theoremstyle{remark}
\newtheorem{remark}[theorem]{Remark}

\def\1{\mathbf{1}}

\def\a{\alpha}
\def\A{{\mathbf A}}
\def\ac{\mathrm{ac}}

\def\b{\beta}

\def\C{{\mathbf C}}

\def\D{\Delta}
\def\d{\delta}

\def\f{\phi}

\def\g{\gamma}

\def\var{{\mathrm{Var}_\C}}
\def\kvar{{K_0 (\mathrm{Var}_\C)}}
\def\kmon{{K_0^{\hat{\mu}}(\mathrm{Var}_\C)}}
\def\khs{K_0 (\mathrm{HS}^{\mathrm{mon}})}

\def\l{\lambda}
\def\L{{\mathbf L}}
\def\Lcal{{\mathcal L}}

\def\N{{\mathbf N}}
\def\Newton{{\mathcal N}}

\def\orb{\mathrm{orb}}
\def\ord{{{\mathrm{ord}}_t}}

\def\orden{\mathrm {ord}}

\def\p{\pi}
\def\P{{\mathbf P}}

\def\Q{{\mathbf Q}}

\def\r{\rho}

\def\R{{\mathbf R}}

\def\z{\zeta}

\def\Z{{\mathbf Z}}

\def\s{\sigma}
\def\t{\tau}

\def\elem(#1,#2){  \{ \frac{#1}{\overline{\ #2\ }}\}  }

\title{Motivic Milnor fiber of a quasi-ordinary hypersurface}
\author{Pedro D. Gonz\'alez P\'erez}

\address{ICMAT.  Fac. de CC. Matem\'aticas.
Univ. Complutense de Madrid.
Plaza de las Ciencias 3. 28040. Madrid. Spain.}

\email{pgonzalez@mat.ucm.es}
\thanks{Both authors are supported by MCI-Spain grant MTM2010-21740-C02.}

\author{Manuel Gonz\'alez Villa}
\address{ Fac. de CC. Matem\'aticas.
Univ. Complutense de Madrid.
Plaza de las Ciencias 3. 28040. Madrid. Spain
 and 
Match, Univ. Heidelberg, Im Neuenheimer Feld 288, 69120 Heidelberg, Germany
}

\email{mgv@mat.ucm.es, villa@mathi.uni-heidelberg.de}



\begin{document}

\begin{abstract}
Let $f: (\C^{d+1}, 0) \to (\C,0)$ be a germ of complex analytic function such that
its zero level defines an irreducible germ of
quasi-ordinary hypersurface $(S, 0) \subset (\C^{d+1}, 0)$.
 We describe the motivic Igusa zeta function, the motivic Milnor fibre
and the Hodge-Steenbrink spectrum of $f$ at $0$ in terms of
topological invariants of $(S, 0) \subset (\C^{d+1}, 0)$.
\end{abstract}
\maketitle

\section* {Introduction}
An equidimensional germ  $(S,0)$ of complex analytic variety  of
dimension $d$ is \textit {quasi-ordinary} (q.o.) if there exists a
finite map $\pi: (S,0)\rightarrow(\C^d,0)$ which is unramified
outside a normal crossing divisor of $(\C^d,0)$.
Q.o.~surface singularities originated in the classical Jung's method to parametrize and resolve
surface singularities (see \cite{Jung}) and are related to the
 classification of singularities by Zariski's
dimensionality type (see \cite{Lipman-Eq}).
In addition to the applications to
equisingularity problems, the
 class of q.o.~singularities is also of interest
to test and study various open questions and conjectures for
singularities in general, particularly in the hypersurface case. In many cases the results 
use the fractional power series parametrizations of q.o.
singularities. The parametrizations allow  explicit computations combining
analytic, topological and combinatorial arguments (see for
instance \cite{Lipman, Gau, PPP-Duke, McEwan-Nemethi, ACNLM-AMS}). 

It is natural to investigate new invariants of singularities arising in the development of motivic integration, 
such as
the motivic Igusa zeta functions for the 
class  of q.o.~singularities with the hope to extend the methods or
results to wider classes (for instance using Jung's
approach).

Denef and Loeser have introduced several (topological, naive motivic Igusa and  motivic Igusa) 
zeta functions   associated 
to a germ  $f: (\C^{d+1}, x) \to (\C, 0)$ of complex analytic
function. All these zeta functions are inspired by the Igusa zeta function in the
$p$-adic setting and admit an expression in terms of  an embedded resolution of singularities (or log-resolution) of 
$f$.  The  \textit{monodromy conjecture} predicts a relation between    the poles of these zeta functions 
 and the eigenvalues of the local monodromy
of $f$ at points of $f(x_1,\dots,x_{d+1})=0$ (see \cite{DL-JAMS, Denef-LoeserIgusa, Denef-LoeserBarca}). 

Let us consider, for instance, the motivic Igusa zeta function  $Z(f, T)_x$ of $f$. 
It is known to be a rational function with respect to the indeterminate $T$  with 
coefficients in certain  Grothendieck ring of varieties. Furthermore, 
it determines some classical invariants of the Milnor fiber of $f$ at $x$ 
(see \cite{Denef-LoeserLef,Denef-LoeserIgusa, Denef-LoeserBarca, GLM-Duke}).
Its  definition can be given in terms of certain jet spaces.

The set $\mathcal{L}_n(\C^{d+1})_x$ (resp. $\mathcal{L} (\C^{d+1})_x$)  of $n$-jets of arcs (resp.  of arcs)
 in $\C^{d+1}$ consists of
$\varphi \in (\C[t] / (t^{n+1})^{d+1}$ (resp. $\varphi \in
\C[[t]]$) such that $\varphi (0) =x$. 
If $n \geq 1$ the set
 \[ \mathcal{X}_{n,1} \colon =
\{ \varphi \in \mathcal{L}_n(\C^{d+1})_x \mid f \circ \varphi =
t^n \mod t^{n+1} \}, \] enjoys a natural action of the group
$\mu_n$ of $n$-th roots of unity. This set defines a class in the
\textit{monodromic Grothendieck ring} $\kmon$ (see its definition
in Section \ref{motmeasure}). 
The \textit{motivic Igusa zeta function}
$Z(f, T)_x$ is the generating series
\[
Z(f, T)_x :=  \sum_{n \geq 0} [ \mathcal{X}_{n, 1} ] \L^{- n (d+1) }  T^n,
\]
where $\L$ denotes the class of the affine line $\A_\C^1$ with the
trivial $\mu_n$-action. The \textit{motivic Milnor
fiber}  of $f$ at $x$ is defined as a certain limit
\[
 \mathcal{S}_{f, x} := - \lim_{T \to \infty} Z(f, T)_x    \in \kmon[\L^{-1}].
\]
Denef and Loeser have proven that  
 $\mathcal{S}_{f,x}$ and ${F}_x$ have
the same \textit{Hodge characteristic}. Hence     $\mathcal{S}_{f, x}$
determines the \textit{Hodge-Steenbrink spectrum} of $f$ at $x$ (see
\cite{Denef-LoeserIgusa} and Section \ref{motmeasure}).

In this paper we give formulas for these invariants when $f =0$ defines an
irreducible germ $(S,0)$
of quasi-ordinary hypersurface singularity at $0 \in \C^{d+1}$.
The germ $(S, 0) $
 is parametrized by a \textit{q.o.~branch},
a fractional power series $\z(x_1^{1/n}, \dots,
x_d^{1/n})$, for some integer $n$. The q.o.~branch has 
a finite set of \textit{characteristic
monomials}, generalizing the notion of characteristic pairs of
plane branches (see \cite{Abhyankar, Lipman}).  These monomials
classify the embedded topology of the germ $(S, 0) \subset
(\C^{d+1}, 0)$ (see \cite{Lipman, Gau}). 
The semiroots of $f$  are auxiliary functions $f_i :
(\C^{d+1}, 0) \to (\C, 0)$, $i=0, \dots, g$ such that  the
hypersurface germ with equation
 $f_i = 0$, is a q.o.~hypersurface parametrized by  a suitable truncation of $\z$ (see
\cite{GPSemi,PPP-Duke}).
Embedded resolutions of $(S,0)$ are built in \cite{GP-Fourier}, as a composition of local 
toric modifications in terms of the characteristic monomials. 
The functions 
$x_1, \dots, x_d, f_0, \dots, f_{g-1}, f_g =f$ 
define suitable local coordinates at intermediate points in these resolutions. 
See Sections \ref{qoh} and \ref{sec-toric}.

We begin by  studying  the order of contact of
an arc $\varphi$ in $\C^{d+1}$ with $\varphi (0) =0$, with the $(d+g+1)$-tuple of functions 
$F:= (x_1, \dots, x_d,   f_0,
\dots, f_g)$,  i.e., we study the values of
\[
 \orden_t (F \circ \varphi ) : = ( \orden_t (x_1 \circ \varphi), \dots, \orden_t (x_d \circ \varphi) ,
\orden_t (f_0 \circ \varphi), \dots, \orden_t (f_g \circ \varphi)).
\]
If the arc $\varphi$ does not factor through the hypersurface $V$
of equation $x_1 \dots x_d f_0 \dots f_g =0$ then $\orden_t (F \circ
\varphi)$ belongs to $\Z^{d+g+1}_{\geq 0}$. We characterize these
orders as certain integral points in the support of a  fan
$\Theta$, consisting of cones in $\R^{d+g+1}$ of dimension $\leq
d+1$ (see Theorem \ref{trop}). This fan is determined explicitely
in terms of the characteristic monomials. It  coincides
essentially with the one given in \cite{GP-Fourier} Section 5.4 in connection
with the combinatorial description of the  toric embedded resolution. The fan $\Theta$ can be seen as the
\textit{tropical set} associated to the embedding of $(\C^{d+1},
0)$ in $(\C^{d+g+1}, 0)$  defined by this sequence of functions.

The second key step is  to determine the motivic measure $\mu( H_{k,1} )$ of 
the set
\[
 H_{k,1}:= \{ \varphi \in \mathcal{L} (\C^{d+1})_0
 \mid \ord (F \circ \varphi) =k, \quad f \circ \varphi = t^{k_{d+g+1}} \mod k_{d+g+1} \}.
\]
We reach this goal by applying the change of variables formula for motivic integrals to
a toric embedded resolution of $(S, 0)$.
Since $H_{k,1}$ has a natural action of the group of $k_{d+g+1}$-th roots of unity,
this measure is an element of the monodromic Grothendieck ring $\kmon[\L^{-1}]$.

Taking into account that the set $\mathcal{X}_{n,1}$ is a disjoint
union of sets of the form $H_{k,1}$ 
we get that: 
\[
 Z(f, T)_0 = \sum
\mu( H_{k,1} ) T^{k_{d+g+1}}\]
(see Proposition \ref{zf1}).
The main result provides formulas for the motivic Igusa zeta
function, the motivic naive Igusa zeta function, the topological
zeta function and the motivic Milnor fiber of a q.o.~germ $f :
(\C^{d+1}, 0) \to (\C, 0)$ in terms of the characteristic
monomials (see Theorem \ref{main-zeta}). As an application we
obtain a formula for the Hodge-Steenbrink spectrum of $f$ at $0$
(see Corollary \ref{main-hs}) and also formulas for a finite set
of candidate poles of the topological zeta function $Z(f, s)_0$
(see Corollary \ref{poles}). We deduce that these invariants are
determined by the embedded topological type of the
q.o.~hypersurface $(S, 0)$. Our results generalize those obtained when $d=1$
by Guibert \cite{Guibert}, Schrauwen, Steenbrink and Stevens \cite{SSS} 
Saito \cite{Saito00}, and Artal et al.~ \cite{ACNLM-AMS}.

The naive motivic zeta function of a germ of complex analytic
function defining a q.o.~ hypersurface, not necessarily
irreducible, was studied by Artal Bartolo, Cassou-Nogu\'es, Luengo and
Melle-Hern\'andez in \cite{ACNLM-AMS} to prove the monodromy
conjecture for this class of singularities. In Section
\ref{comparison} we compare the method of \textit{Newton maps}
used in \cite{ACNLM-AMS} and  the toric resolutions of
\cite{GP-Fourier}, with the help of quotients of $\C^{d+1}$ by finite
abelian groups. We point out in Example \ref{false} the phenomenon
of \textit{false reducibility}, discovered in \cite{MGV},  which
may appear when using Newton maps.  This phenomenon, in the case $d \geq 2$,  leads   to
inaccurate formulas for the zeta functions in \cite{ACNLM-AMS}. 
However,  we check that the set of candidate poles 
which arise by our method coincides with the one given in 
\cite{ACNLM-AMS}. 
 See Remarks \ref{comp},  \ref{comp-1} and \ref{comp-4}.

An application to the description of the \textit{log-canonical 
threshold} for this class of singularities will be given elsewhere. 
In the case $d = 1$, see for instance  \cite{ACNLM-CM-08}.

\smallskip

\noindent { {\bf Acknowledgments:}
We are grateful to Nero Budur, Alejandro Melle and to Andr\'as N\'emethi 
for several discussions and their useful suggestions concerning this paper. 
The second author is also grateful to Enrique Artal-Bartolo, 
Mirel C\u aibar, Pierrette Cassou-Nogu\'es, Gary Kennedy, Anatoly Libgober, Ignacio Luengo, 
Laurentiu Maxim, Lee McEwan and to the Alfred Renyi Institute of Mathematics (Hungary) 
and the Department of Mathematics, Statistics and Computer Sciences of the University of Illinois at Chicago (USA).}

\section{Invariants of singularities and motivic integration}  \label{motmeasure}

In this section we recall the definitions of the invariants we study in this paper.

\subsection*{The arc space of an algebraic variety}
Let us denote by $\var$ the category of complex algebraic varieties, i.e.,
reduced and separated schemes of finite type over the complex field $\C$.
If $0 \ne \varphi \in \C [[t]]$ we denote by $\ord(\varphi)$ (resp. $\ac (\varphi)$  the \textit{order}
 (resp. \textit{angular coefficient})
of the series  $0 \ne \varphi \in \C [[t]]$, that is, if
$\varphi = t^{n} u$,  where $u$ is a unit in $\C[[t]]$, then $\ord(\varphi) = n$ and $\ac(\varphi) = u(0) \in \C^*$.

If $n \geq 0$ the space $\mathcal{L}_n(X)$ of arcs modulo
$t^{n+1}$ on a complex algebraic variety  $X$ has also the
structure of algebraic variety over $\C$. The canonical morphisms
$\pi^{n'}_n : \mathcal{L}_{n'}(X) \rightarrow \mathcal{L}_n(X)$,
$n' >n$ induced by truncation define a projective system. The arc
space $\mathcal{L}(X)$ of $X$ is the projective limit of
$\mathcal{L}_n(X)$. We consider the spaces $\mathcal{L}_n(X)$ and
$\mathcal{L}(X)$ with their reduced structures. The arc space
$\mathcal{L}(X)$ is a separated scheme over $\C$ but it is not of
finite type in general. If $X$ is an affine variety
given by the equations $f_i (x_1, \dots, x_d)=0$, $i \in I$ then
the $\C$-rational points of $\mathcal{L}_n (X)$ (resp. of
$\mathcal{L}(X)$) are $d$-tuplas $\varphi=(\varphi_1
,\varphi_2,\dots,\varphi_d) \in (\C[t]/ (t^{n+1}))^d$ (resp.
$\varphi \in (\C[[t]])^d)$) such that  $f_i \circ \varphi =0 \mod
t^{n+1}$, (resp. $f_i \circ \varphi =0$),  for all $i \in I$.

Truncating arcs $\mod t^{n+1}$ defines a morphism of schemes
$\pi_n:\mathcal{L}(X) \rightarrow \mathcal{L}_n(X)$.
Note that $\pi_0(\mathcal{L}(X))=X$. For any ${x} \in X$ we denote $(\pi_0^n)^{-1}({x})$
by  $\mathcal{L}_{n}(X)_{x}$ and $(\pi_0)^{-1}({x})$ by  $\mathcal{L}(X)_{x}$.
If the variety $X$ is smooth  the maps $\pi_n$ are surjective and the maps $\pi^{n'}_n$ are locally
trivial fibrations with fibre ${\A}^{(n'-n) \dim X }_\C$.

\subsection*{Grothendieck rings of algebraic varieties over $\C$}

 The \textit{Grothendieck group} of algebraic varieties over $\C$,
is the abelian group $\kvar$ generated by symbols $[X]$, where $X$
is a complex algebraic variety, with the relations $[X]=[Y]$ if
$X$ and $Y$ are isomorphic, and $[X] = [Y] + [X \backslash Y]$
where $Y$ is a Zariski closed subset of $X$. Setting $[X][X'] :=
[X \times_\C X']$ defines a ring structure on $\kvar$ in such a
way that $1$ is equal to the class of a closed point. We denote by
${\L}$ the class of ${\A}^1_\C$ in $\kvar$ and by $\mathcal{M}_\C$
the ring $\kvar [\L^{-1}]$.

We denote by $\mu_n$ the algebraic group of $n$-roots of unity,
i.e., $\mu_n = {\rm Spec}(\C[x]/(x^n-1))$ for $n \geq 1$. The
groups $\mu_n$, $n \geq 1$, together with the morphisms $\mu_{nd}
\rightarrow \mu_n$ given by $x \mapsto x^d$ form a projective
system. We denote by   $\hat{\mu}$ its projective limit. A
\textit{good $\mu_n$-action} on a complex algebraic variety $X$ is
a group action $\mu_n \times X \rightarrow X$ such that it is a
morphism of algebraic varieties over $\C$ and
 each orbit is contained in an open affine subvariety of $X$.
A \textit{good $\hat{\mu}$-action} is a group action $\hat{\mu} \times X \rightarrow X$
which factors through a good $\mu_n$-action for some $n$.

The \textit{monodromic Grothendieck ring} $\kmon$ of algebraic varieties over $\C$ is
the abelian group on the symbols $[X,\hat{\mu}]$
where $X$ an algebraic variety with a good $\hat{\mu}$-action, with the relations $[X,\hat{\mu}]
= [Y,\hat{\mu}]$ if $X$ and $Y$ are isomorphic as varieties with good $\hat{\mu}$-actions and
 $[X,\hat{\mu}] = [Y,\hat{\mu}] + [X \backslash Y,\hat{\mu}]$ , where $Y$ is a Zariski
 closed subset of $X$, considered with the $\hat{\mu}$-action induced by the one on $X$.
We denote by ${\L}$
the affine line ${\A}^1_\C$ with the trivial $\hat{\mu}$-action.
The product is defined as in $\kvar$, by  adding  the relation
  $[X\times V,\hat{\mu}]=[X ]  {\L}^d$, for $V$ any
$d$-dimensional affine space over $\C$ equipped by any linear
$\hat{\mu}$-action. We denote by  $[\mu_n] \in \kmon $ the class
of $\mu_n$ equipped with the $\mu_n$ action; notice that
$[\mu_n]=n \in \kvar$ but $[\mu_n] \not = n \in \kmon$ if $n> 1$.
Another example is the following:
\begin{example}\label{P4basicexampleGT}
Let $f \in \C[x_1^{\pm 1},\dots,x_{m}^{\pm 1}]$
be a quasi-homogeneous Laurent polynomial, i.e.,
there exists  integers $D$ and  $w_1, \dots, w_{m}$ such that
$f(t^{w_1}x_1,\dots,t^{w_{m}} x_{m})=t^Df(x_1,\dots,x_{m})$.
Then $f = 1 $ defines a \textit{quasi-homogeneous hypersurface} $V$ in the torus $(\C^*)^m$,
which is endowed with the    $\mu_D$-action
given by $ \a \cdot (x_1, \dots, x_m) \mapsto (\a^{w_1} x_1, \dots, \a^{w_d} x_d)$ for $\a \in  \mu_D$.
\end{example}

\subsection*{Motivic integration}

Let $X$ be a complex algebraic variety of pure dimension. A subset
$A \subset \mathcal{L}(X)$ is \textit{cylinder} if
$A=\pi^{-1}_n(C)$ for some $n$ and some  constructible subset $C \subset
\mathcal{L}_n(X)$.
A subset  $A \subset \mathcal{L}(X)$ is \textit{stable} if for some $m \in \mathbf{N}$ the set
 $\pi_m (A)$ is constructible,
 $A = \pi^{-1}_m \pi_m (A)$ and the morphism
$\pi_{n+1}(A) \rightarrow \pi_n(A)$ is a piecewise locally trivial
fibration with fibre $\mathbf{A}^{\dim X}_\mathbf{C}$, for $n >
m$. If $X$ is smooth any cylinder is stable (see
\cite{Denef-LoeserInvent}).

We denote by    $\hat{\mathcal{M}}_\C$ the
completion of ${\mathcal{M}}_\C$ with respect to the filtration
$F^m\mathcal{M}_\C$, with $m \in {\Z}$, where $F^m\mathcal{M}_\C$ is the subgroup of $\mathcal{M}_\C$ generated
by elements $[S]{\L}^{-i}$ with $S \in \var$ and $i-\dim S \geq m$.

Denef and Loeser have introduced a class of \textit{measurable}
subsets of $\mathcal{L}(X)$, containing the stable subsets,
together with a measure $\mu_X(A) \in \hat{\mathcal{M}}_\C$, for
any  measurable set $A \subset \mathcal{L}(X)$, see
\cite[Appendix]{Denef-LoeserMcKay}. This motivic measure have the
following properties:
 \begin{proposition} \label{motmes}     \
\begin{enumerate}
\item[(i)] If $A \subset \mathcal{L}(X)$ is stable at level $n$,
then $\mu_X(A)=[\pi_n(A)]{\L}^{-n \dim X}.$

\item[(ii)] If $A \subset \mathcal{L}(X)$ is contained in
$\mathcal{L}(S)$ with $S$ a reduced closed subscheme of $X$ with
$\dim S < \dim X$,  then $\mu_X(A)=0.$

\item[(iii)] Let $A_i$ with $i \in \N$ be a family of mutually
disjoint measurable subsets of $\mathcal{L}(X)$. Assume that $A :=
\cup_{i \in \N}A_i$ is a measurable set.  Then $\sum_{i \in \N}
\mu(A_i)$ converges to $\mu(A)$ in $\hat{\mathcal{M}}_\C$.
\end{enumerate}
\end{proposition}

If $A \subset \mathcal{L}(X)$ is a measurable set,
and
$\alpha:A\to {\Z} \cup \{\infty\}$ is a function, then
${\L}^{-\alpha}$ is \textit{integrable} on $A$ if the fibres
of $\alpha$ are measurable,
$\alpha^{-1}(\infty)$ has measure zero and the
\textit{motivic integral}
$
\int_{A}{\L}^{-\alpha} d\mu_X:=\sum_{n\in {\Z}}
\mu_X(A\cap \alpha^{-1}(n)){\L}^{-n}
$
converges in $\hat{\mathcal{M}}_\C.$

\begin{theorem}{\rm (Change variables formula \cite[Lemma 3.3]{Denef-LoeserInvent})} \label{P4changevariables}
Let $\pi: Y \to X$ be a proper modification
between two smooth complex algebraic varieties of pure dimension.
Let $A$  be a measurable subset of $\mathcal{L}(X)$.
Assume that $\pi$ induces a bijection between $B \subset\mathcal{L}(Y)$ and $A$. Then, for any
function
$\alpha:A\to {\Z}\cup\{\infty\}$ such that
${\L}^{-\alpha}$ is integrable on $A$, we have
$$
\int_{A} {\L}^{-\alpha} d\mu_X = \int_{B} {\L}^{-\alpha \circ \pi - {\rm ord}_t J_{\pi}(y)} d\mu_Y,
$$
where $\ord J_{\pi}(y)$, for any $y\in \mathcal{L}(Y),$ denotes the order of
the Jacobian of $\pi$ at $y$.
\end{theorem}

\begin{remark}
If the sets of arcs considered in Proposition \ref{motmes} and Theorem  \ref{P4changevariables}
are equipped with a good $\hat{\mu}$-action then these results hold replacing $\mathcal{M}_\C$ (resp.
$\hat{\mathcal{M}}_\C$) by
$\mathcal{M}_\C^{\hat{\mu}}$ (resp. by $\hat{\mathcal{M}}_\C^{\hat{\mu}}$), where
$\mathcal{M}_\C^{\hat{\mu}} := \kmon[\L^{-1}]$ and
$\hat{\mathcal{M}}_\C^{\hat{\mu}}$ is the completion of $\mathcal{M}_\C^{\hat{\mu}} $
with respect to the filtration introduced above.
\end{remark}

\subsection*{Hodge characteristic and spectrum}

A \textit{Hodge structure} is a finite dimensional vector space
over $\Q$, equipped with a bigrading $H \otimes \C = \oplus_{p, q
\in \Z} H^{p,q}$ such that $H^{p,q}$ is the complex conjugate of
$H^{q, p}$ and each weight $m$ summand $\oplus_{p+q = m} H^{p,q}$
is defined over $\Q$. Let $\mathrm{HS}$ (resp.
$\mathrm{HS}^{\mathrm{mon}}$)  be the abelian category of Hodge
structures (resp. equipped in addition with a quasi-unipotent
endomorphism). The elements of the Grothendieck ring
$K_0(\mathrm{HS})$ (resp. $\khs$) are formal differences of
isomorphic classes of Hodge structures (resp. equipped with a
quasi-unipotent endomorphism). The product is induced by the
tensor product of representatives.

A mixed Hodge structucture is a finite dimensional vector space
$V$ with a finite increasing filtration ${W}_{\bullet} V$, such
that the associated graded vector space $\mathrm{Gr}_\bullet^W V$
underlies a Hodge structure with $\mathrm{Gr}_m ^W V$ as weight
$m$ summand. 

The canonical class $[V] := \sum_l [\mathrm{Gr}_l^W]
\in  K_0(\mathrm{HS}) $ is associated to the mixed Hodge structure
$V$. If $X$ is a complex algebraic variety the simplicial
cohomology groups $H_c^i( X, \Q)$ of $X$, with compact support,
have a mixed Hodge structure.

The \textit{Hodge characteristic} 
of a complex algebraic variety $X$ is 
       $\chi_h (X) = \sum_i (-1)^i [ H^i _c (X,
\Q)] \in K_0 (\mathrm{HS})$. 
Since $\chi_h(\L) $ is invertible in $K_0 (\mathrm{HS})$ we get by additivity that
$ \chi_h $ extends to a ring homomorphism $ \chi_h : \mathcal{M}_\C \to K_0
(\mathrm{HS})$. 
If in addition $X$ is equipped with a good
${\mu}_n$-action then $\chi_h (X)$ can be seen as an element of
$\khs$ and we get in this way a ring homomorphism
\[\chi_h^{\mathrm mon} : \mathcal{M}_\C^{\hat{\mu}} \to \khs,\]
called the \textit{monodromic Hodge characteristic}.
We refer to \cite{Chris-Steenbrink, Saito-MHM} for more on Hodge structures. 

We denote by ${\Z}[t^{1/{\Z}}]$ the ring $\bigcup_{n\geq 1}{\Z}[t^{1/n},t^{-1/n}]$. 
The  {\it Hodge
spectrum} is a linear map
$$hsp: \khs \rightarrow {\Z}[t^{1/{\Z}}]$$
given by 
\[
hsp([H]) = \sum_{\alpha \in \mathbb{Q} \cap [0,1[}t^{\alpha}(
\sum_{p,q \in {\Z}} h^{p,q}_{\alpha}(H)t^p),\] for $H$  a  Hodge
structure with a quasi-unipotent endomorphism, where
$h^{p,q}_\alpha$ denotes the  dimension of the eigenspace
$H^{p,q}_\alpha$ corresponding to the eigenvector $e^{2\pi i
\alpha}$. 
The composition
\[
 Sp:= hsp \circ \chi_h^{\mathrm{mon}} : {\mathcal{M}_\C^{\hat{\mu}}}  \longrightarrow  {\Z}[t^{1/{\Z}}]
\]
 is a group
homomorphism.

In the following Lemmas we recall some elementary properties of
the map $Sp$.   See \cite{Steenbrink76, Steenbrink-Asterisque} for more properties and applications of this map.

Let $Q_{n,r}(t)$ denote the polynomial
$\frac{t^{1/n}-t}{1-t^{1/n}} \frac{t^{1/r}-t}{1-t^{1/r}} \in
{\Z}[t^{\frac{r}{n}}].$ We expand the polynomial  $Q_{n,r}(t) =
\sum_{\alpha} c_\alpha t^\alpha$ and define  the polynomials
$Q^{<}_{n,r} (t): = \sum_{\alpha <1} c_\alpha t^\alpha$ and
$Q^{>}_{n,r} (t) := \sum_{\alpha >1} c_\alpha t^\alpha$.
\begin{lemma}\label{P4ejemsp} We have the following properties:
\begin{enumerate}
\item[(i)] If $X$ is a complex algebraic variety, endowed with a
good $\hat{\mu}$-action, then the equality   $Sp([X]  \L^d) =
Sp([X]) t^d$ holds for any $d \in \Z_{\geq 0}$. 

\item[(ii)] $Sp([\mu_r]) = \frac{1-t}{1-t^{1/r}}$.

 \item[(iii)] If  $\gcd(n,r)=1$ then we have
$Sp([\{(x,y) \in {\A}^2_\C \mid y^n - x^r= 1\}]) = t -
\frac{t^{1/n}-t}{1-t^{1/n}} \frac{t^{1/r}-t}{1-t^{1/r}}$.
\item[(iv)]  If  $\gcd(n,r)=1$ and $e > 1$  the following
equality holds
\[Sp([\{(x,{y}) \in \A^2_\C \mid ({y}^n - x^r)^e=1\}])=
\frac{1-t}{1-t^{1/e}} \big ( t - Q^{<}_{n,r}(t^{1/e}) - t^{\frac{e
-1}{e}}Q^{>}_{n,r}(t^{1/e})\big ).\]
\end{enumerate}
\end{lemma}
\textit{Proof.} 
Claim (i) can be seen as a consequence of (2.1.2) in \cite{Saito-MA91}.
For (ii),  (iii)
and (iv)  we refer to Lemma 3.4.2 and Lemma 3.4.3 of
\cite{Guibert}.  \hfill $\square$

\begin{notation}\label{notxhsp}
We denote the polynomial $Sp([\mu_r])$ by $P_r$ (see Lemma
\ref{P4ejemsp}). We set also
\[
 P_{n,r}:=Sp([\{(x,y) \in (\C^*)^2 \mid  y^n - x^r= 1\}]) \mbox{
 and }
P_{n,r}^{(e)}:=Sp([\{(x,y) \in (\C^*)^2 \mid  (y^n - x^r)^e =
1\}]).
\]
\end{notation}

\begin{corollary}\label{cornotxhsp} We have the following
identities:
\[
\begin{array} {lcl}
P_{n,r} & = &  t - \frac{t^{1/n}-t}{1-t^{1/n}}
\frac{t^{1/r}-t}{1-t^{1/r}} -  P_n  - P_r
\\
P_{n,r}^{(e)}& =  &
P_e  \big ( t - Q^{<}_{n,r}(t^{1/e}) - t^{\frac{e
-1}{e}}Q^{>}_{n,r}(t^{1/e})\big ) - P_{ne} - P_{re}
\end{array}
\]
 \end{corollary}
\textit{Proof.} Let $n, r, e \in {\Z_{> 0}}^d$. The polynomial
$f=y^n-{x}^{r}$ is quasi-homogeneuos with  weights $w_x=n, w_y=r$ and
degree $D=nr$.   The  identity
\[  [\{({x},y) \in \C^2 \mid (y^n-{x}^{r})^e=1\}]
 = [\mu_{ne}] + [\mu_{re}] + [\{({x},y) \in (\C^*)^2 \mid (y^n-{x}^{r})^e=1\}]
\]
holds in $K_0^{\hat{\mu}}({\rm Var} _\C)$ and  follows from the
partition of   $\{({x},y) \in \C^{2} \mid (y^n - {x}^{r})^e=1\}$
into three disjoint subsets:
$\{({x},y) \in \C^{2} \mid y^{ne} =0, {x}^{re} = 1\}$, $\{({x},y)
\in \C^{2} \mid y^{ne} =1, {x}^{re} =0\}$ and $\{(x,y) \in
(\C^*)^2 \mid (y^n - {x}^{r})^e=1\}$. The $\mu_{D}$-action
introduced in Example \ref{P4basicexampleGT} is compatible with
this partition.
Then the identities hold by Lemma \ref{P4ejemsp}. \hfill $\square$

\subsection*{The Milnor fiber and the  Hodge-Steenbrink spectrum}

 Let $f : \C^{d+1} \rightarrow \C$ be a non constant morphism. Suppose that  $f(x)=0$.
 Take $0 < \delta << \epsilon < 1$ and denote by ${B}(x,\epsilon)$
 the open ball  of radius $\epsilon$ centered at  $0 \in \C^{d+1}$,
  by ${D}_\delta$ the open disk of radius $\delta$ centered at $x \in \C$
  and ${D}^*_\delta:= {D}_\delta \setminus \{0\}$. If  $\epsilon$ and $\delta$ are small enough
 $$f : f^{-1}({D}^*) \cap  {B}(x,\epsilon) \rightarrow {D}^*_\delta.$$
 is a locally trivial fibration. The fibre
 ${F}_x := f^{-1}(z) \cap {B}(x,\epsilon)$ for $z \in {D}^*_\delta$,
which is defined up to isotopy,
 is called the  \textit{Milnor fibre} of $f$ at the point  $x$.
 See \cite{Milnor}.

The Milnor fibre is equipped  with the monodromy operator $M_{f,x}
:  H^\cdot ({F}_x, \Q) \rightarrow  H^\cdot ({F}_x, \Q)$ on the
cohomology ring.
By the Monodromy Theorem the operator $M_{f,x}$ is
quasi-unipotent, that is, there exist  $A, B \in {\Z_{\geq 0}}$
such that $(M^A_{f,x} - I)^B =0$
(see \cite{MR0344248}). The cohomology groups $H^i ({F}_x, \Q)$
carry a natural mixed Hodge structure which is compatible with the
monodromy operator $M_{f,x}$ (see \cite{Steenbrink76} for isolated singularities and \cite{Navarro}
 for the general case). The \textit{Hodge characteristic} of
the Milnor fiber is defined as
\[
\chi_h^{\mathrm{mon}} (F_{x}) = \sum_i (-1)^i [H^i ({F}_x,\Q) ]   \in 
\khs.
\]

      Following \cite{Saito-MA91, Denef-LoeserIgusa, Denef-LoeserBarca} we
introduce the     Hodge-Steenbrink spectrum   $Sp ' (f,x)$
of $f$ at $x$ as:
\begin{definition}
We set  
  \[ 
Sp ' (f,x) := (-1)^d hsp (  \chi_h^{\mathrm{mon}} (F_{x}) -1)
  \quad \mbox{ and } \quad 
 Sp (f, x) := t^{d+1} \iota (Sp '(f, x)),  
\]
where $\iota : \Z [t^{1/\Z}] \to \Z [ t^{1/ \Z} ]$ is the isomorphism such that $t^\a \mapsto t^{-\a}$.  
\end{definition}

\begin{remark}     \label{sp-normal}
There are several normalizations for the Hodge-Steenbrink spectrum in the literature. 
We have followed those of Saito \cite{Saito-MA91}. In     \cite{Denef-LoeserIgusa} 
$Sp' (f,x) $ (resp. $Sp (f,x)$)  is denoted by $HSp' (f, x)$ (resp. $HSp(f,x)$), while 
in \cite{Denef-LoeserBarca} and \cite{Guibert} the polynomial $Sp'(f, x)$ is denoted by $hsp(f,x)$. 
In the paper \cite{Steenbrink76} the Hodge-Steenbrink spectrum of $f$ at $x$ is 
defined to be the polynomial $Sp(f, x)$, while in \cite{Steenbrink-Asterisque,Chris-Steenbrink,Varchenko} 
the Hodge-Steenbrink spectrum of
$f$ at $x$ differs of the polynomial $Sp(f, x)$ by multiplication by $t$ (see \cite{Kulikov}). 
If $f$ has an isolated singular point at $x$ we have that 
$ Sp' (f, x) =  Sp (f, x)$  (see \cite{Steenbrink-Asterisque}). 
\end{remark}

\subsection*{Motivic Igusa zeta function of a function at a point}
Let $f: {\C}^{d+1} \rightarrow \C$ be a non constant morphism such
that $f(x)=0$.
For each natural $n$, we set
$$     \mathcal{X}_{n}:= \{ \varphi \in \mathcal{L}_{n}(\C^{d+1})_{x} \, | \, {\rm ord}_t f \circ \varphi =n \}.$$
The set $\mathcal{X}_{n}$ is a locally closed subvariety of
$\mathcal{L}_{n}(\C^{d+1})_{x}$. We also consider the morphism
$${\ac}(f) : \mathcal{X}_{n} \rightarrow \C^*, \varphi \mapsto \ac (f \circ \varphi).$$
There is a natural action of $\C^*$ on $\mathcal{X}_{n}$ given by
$a\cdot\varphi(t)= \varphi(at)$. Since ${\ac}(f)(a\cdot \varphi)=
a^n{\ac}(f)(\varphi)$, the morphism ${\ac}(f)$ is a locally
trivial fibration. We denote by $\mathcal{X}_{n,1}$ the fibre
${\rm ac}(f)^{-1}(1)$. The previous $\C^*$-action on
$\mathcal{X}_n$ restricts to a $\mu_n$-action on
$\mathcal{X}_{n,1}$.

\begin{definition}  \cite{Denef-LoeserIgusa} \label{igusa}
\begin{enumerate} \item[(i)] The \textit{naive motivic Igusa zeta function}
of $f$ at $x$ is
$$Z^{\mathrm{naive}}(f, T)_x:=\sum_{n \geq 1}[\mathcal{X}_{n}]{\L}^{-n(d+1)}T^n \in \hat{\mathcal{M}}_\C[[T]].$$
\item[(ii)]  The \textit{motivic Igusa zeta function} of  $f$ at $x$ is
$$Z(f,T)_x:=\sum_{n \geq 1}[\mathcal{X}_{n,1}]{\L}^{-n(d+1)}T^n \in \hat{\mathcal{M}}_\C^{\hat{\mu}}[[T]].$$
\end{enumerate}
\end{definition}

If $A$ is $\mathcal{M}_\C$ or $\mathcal{M}^{\hat{\mu}}_\C$ we
denote by $A[[T ]]_{{\mathrm{sr}}}$ the $A$-submodule of $A[[T ]]$
generated by 1 and by finite products of terms $\L^e T^i (1- \L^e
T^i)^{-1}$ with $e \in \Z$ and $i \in \Z_{>0}$.

Denef and Loeser gave formulas for $Z^{\mathrm{naive}}(f,T)_x$ and
$Z(f,T)_x$
 in terms data associated to an embedded resolution of
singularities of $f$ (see \cite{Denef-LoeserLef} Theorem 2.4,
\cite{Denef-LoeserIgusa} Theorem 2.2.1). In particular they proved
that 
\[
Z^{\mathrm{naive}}(f,T)_x \in \mathcal{M}_\C[[T]]_{{\mathrm{sr}}}
\quad \mbox{ and } \quad  Z(f,T)\in
\mathcal{M}_\C^{\hat{\mu}}[[T]]_{{\rm sr}}.\]

The \textit{topological zeta function} of  $f$ at $x$ is a
rational function $Z_{\mathrm{top}}(f,s) \in \Q (s)$ defined by
Denef and Loeser in \cite{DL-JAMS} in terms of an embedded
resolution of $f$. The motivic zeta function
$Z^{\mathrm{naive}}(f,T)_x$ determines $Z_{\mathrm{top}}(f,s) \in
\Q (s)$.
\begin{proposition}  \cite{Denef-LoeserBarca,Denef-LoeserIgusa}
\label{ztop-rel} We have that  $Z_{\mathrm{top}}(f,s) =
\chi_{\mathrm{top}}(Z^{\mathrm{naive}}(f, \L^{-s}))$, where
$\chi_{\mathrm{top}}$ denotes the topological Euler
characteristic.
\end{proposition}

\begin{lemma}\cite[Lemma 4.1.1]{Denef-LoeserIgusa} \label{P4limiteremark}
There exists a unique ring  homomorphism
$\mathcal{M}_\C^{\hat{\mu}}[[T]]_{{\rm sr}} \rightarrow
\mathcal{M}_\C^{\hat{\mu}}$, $\psi \mapsto {\lim}_{T \rightarrow
\infty}  \psi$, such that $ {\lim}_{T \rightarrow \infty}
 {\L}^{e}T^{i}(1 - {\L}^{e}T^{i})^{-1} = -1$.
\end{lemma}

\begin{lemma}\label{P4limiteremarkb}
If $\frac{p(T)}{q(T)} \in \mathcal{M}_\C^{\hat{\mu}} [[T]]_{{\rm
sr}}$ and
 ${\deg}_T p(T) < {\deg}_T q(T)$
 then we get
${\lim}_{T \rightarrow \infty} \frac{p(T)}{q(T)} =0$.
\end{lemma}
\textit{Proof.} It follows easily from Lemma \ref{P4limiteremark}.
\hfill $\square$

\begin{definition} \cite{Denef-LoeserLef, Denef-LoeserIgusa}
The \textit{motivic Milnor fiber} of $f$ at the  point $x$ is
defined as
\[
\mathcal{S}_{f,x} := - {\rm lim}_{T \rightarrow \infty} Z(f,T)_x
\in \mathcal{M}^{\hat{\mu}}_\C.
\]
\end{definition}

The  motivic Milnor fiber   $\mathcal{S}_{f,x}$  can be seen as a
motivic incarnation of the classic Milnor fibre ${F}_x$ of $f$ at
$x$. For instance $\mathcal{S}_{f,x}$ and ${F}_x$ have the
same Hodge characteristic.
\begin{theorem}[Denef and Loeser \cite{Denef-LoeserLef, Denef-LoeserIgusa}]
\label{P4thmincarnation}
The following equality holds in $\khs$:
$$\chi^{\rm mon}_h({F}_x)= \chi^{\rm mon}_h(\mathcal{S}_{f,x}).$$
\end{theorem}

\begin{corollary}\label{incar}
We have the equality $Sp'(f,x) = (-1)^d hsp( \chi_h (\mathcal{S}_{f,x}) -1)$.
\end{corollary}

\begin{remark}
The definitions and results of this section extend to the case
when   $f \colon (\C^{d+1}, x) \to (\C, 0)$ is a germ of complex
analytic function.
\end{remark}

\section{Quasi-ordinary hypersurface singularities}\label{qoh}

A germ $(S, 0)$ of complex analytic variety equidimensional of dimension $d$
 is \textit{quasi-ordinary} (q.o.) if there exists a
finite projection $\p: (S, 0) \rightarrow (\C^d,0)$ which is a
local isomorphism outside a normal crossing divisor. If $(S, 0)$
is a hypersurface there is an embedding $(S, 0) \subset (\C^{d+1},
0) $, defined by an equation $f= 0$,  where $f \in \C \{ x_1,
\dots, x_d  \} [{y}]$ is a \textit{q.o.~polynomial}, that is,  a
Weierstrass polynomial in ${y}$  with discriminant $\D_{{y}} f$ of the
form $\D_{{y}} f = x^\d  u$ for a unit $u$ in the
ring $ \C \{ x \}$ of convergent power series in the variables $x=
(x_1, \dots, x_d)$ and $\d \in \Z^d_{\geq 0}$.     In this coordinates
the q.o. projection $\p$ is the restriction of the projection
\begin{equation}
\label{projection} \C^{d+1}_\C \to \C^d_\C, \quad (x_1, \dots, x_d,
y) \mapsto (x_1,\dots, x_d).
\end{equation}
We abuse of notation denoting the projection (\ref{projection}) also by $\p$. 
The Jung-Abhyankar theorem guarantees that the roots of a q.o.~
polynomial $f$, called {\it q.o.~ branches}, are fractional power
series in the ring $\C \{ x^{1/m}\}$, for some integer $m \geq 1$ (see
\cite{Abhyankar}).

In this paper we suppose that the germ $(S,0)$ is analytically irreducible, that is,
the polynomial  $f $ is irreducible in  $\C \{ x_1, \dots, x_d  \} [{y}]$.

 \begin{notation}
If $\a, \b \in \Q^d$ we consider the preorder relation given by
$\a \leq \b$ if $\b \in \a + \Q^d_{\geq 0}$. We set also $\a < \b$
if       $\a \leq \b$ and $ \a \ne  \b$. The notation $\a \nleq
\b$ means that the relation $\a \leq \b$ does not hold. In $\Q^d
\sqcup \{ \infty \}$ we set that $\a < \infty$.
 \end{notation}

\begin{lemma} \label{expo}   \cite[Proposition 1.3]{Gau}
Let  $f \in \C \{ x_1, \dots, x_d \} [{y}] $ be an irreducible
q.o.~polynomial. Let $\z$ be a root of $f$ with expansion:
\begin{equation} \label{expan}
\z = \sum c_{\l} x^\l.
\end{equation}
There exist unique vectors $\l_1, \dots,
\l_g \in \Q^d_{\geq 0}$ such that $\l_1 \leq \cdots \leq \l_g$ and
the three conditions below hold.
 We set the notations  $\l_0 = 0$, $\l_{g+1} = \infty$,
  and  introduce the lattices
$M_0 := \Z^d$, $M_j :=
M_{j-1} + \Z \l_j$, for $j=1, \dots, g$.
\begin{enumerate}
\item [(i)]  We have that      $c_{\l_j} \ne 0$ for $j =1, \dots, g$.
\item [(ii)] If $c_{\l} \ne 0$ then the vector $\l$ belongs to the lattice $M_j$, where $j$ is the unique integer
such that  $\l_j \leq \l$ and
    $\l_{j+1} \nleq \l$.

\item [(iii)] For $j=1, \dots, g$, the vector $\l_j $ does not
belong to $ M_{j-1}$.
\end{enumerate}
   If $\z \in \C \{ x^{1/m}\}$ is a fractional power series satisfying the three conditions above then
$\z$ is a q.o.~branch.
\end{lemma}

\begin{remark} Let us denote by $K$ the field of fractions of $\C \{ x_1, \dots, x_d \}$.
If $\t \in \C \{ x_1^{1/m}, \dots, x_d^{1/m} \} $ is a q.o.~branch then  the minimal polynomial $F \in K[y]$ of $\t$
over $K$ has coefficients in    the ring
$\C \{ x_1, \dots, x_d \}$.
Then the equation $ F = 0 $  defines the q.o.~hypersurface parametrized by $\t$.
\end{remark}

\begin{definition}
The exponents $\l_1 , \dots, \l_g $ in Lemma \ref{expo}
 (resp. the monomials $x^{\l_1}, \dots, x^{\l_g}$) are called {\em characteristic}
of the q.o.~ branch $\z$. The integers  $n_j := [M_{j-1} : M_j]$,  $j=1, \dots, g$,
are strictly greater than $0$ by Lemma \ref{expo}. We set also
$n_0 =1$, $e_0 = n_1 \cdots n_g$ and $e_{j} := e_0 / (n_1 \dots
n_j)$ for $j=1, \dots, g$.  By convenience we denote by  $M$  the
lattice $M_g$ and by $N$ its dual lattice.
\end{definition}

\begin{notation} \label{lamb-coord}
We denote by $( \l_{j,1}, \dots, \l_{j,d})$ the coordinates of the
characteristic exponent $\l_j$ with respect to the canonical basis
of $\Q^d$ and by  $\geq_{\mbox{\rm lex}}$ the lexicographic order.
\end{notation}

\begin{definition}
  The q.o.~
branch $\z$
 is normalized if
\[
(\l_{1,1}, \dots, \l_{g,1}) \geq_{\mbox{\rm lex}} \cdots
\geq_{\mbox{\rm lex}} (\l_{1,d}, \dots, \l_{g,d}),
\]
and if $\l_1$ is not of the form $(\l_{1,1}, 0, \dots, 0)$ with
$\l_{1,1} < 1$.
\end{definition}
 Lipman proved that if the  q.o. branch $\z$ is not normalized then
 there is a normalized q.o. branch $\z'$ parametrizing the same
 germ and  gave \textit{inversion formulas} relating the
 characteristic monomials of $\z$ and $\z'$ (see \cite{Gau}, Appendix).
Lipman and Gau studied q.o.~singularities from a topological point
of view. They proved that the characteristic exponents of a
normalized q.o.~branch $\z$ parametrizing $(S,0)$, classify the
\textit{embedded topological type} of the hypersurface germ $(S,
0) \subset (\C^{d+1}, 0)$ (see \cite{Gau, Lipman}).

\begin{lemma}
If $\z$ is a q.o.~branch of the form (\ref{expan}) then
the series $\z_{j-1} := \sum_{ \l \ngeq \l_j} c_\l x^\l$ is a
q.o.~branch with characteristic exponents $\l_1, \dots, \l_{j-1}$
for $j=0, \dots, g$.
\end{lemma}
\textit{Proof.}
This is consequence of  Lemma \ref{expo}.
\hfill $\square$

 \begin{definition}  \label{semi}
The q.o.~ branch  $\z_j$ parametrizes the q.o.~hypersurface $(S_j,
0)$. We denote  by $f_j \in \C \{x_1, \dots, x_d\} [{y}]$ the
q.o.~polynomial defining $(S_j,0)$, for $j=0, \dots, g-1$. By
convenience we also denote $\z$ also by $\z_g$, hence $f= f_g$ and $S=
S_g$.
\end{definition}

\begin{definition}
The elements of $M$ defined by:
\begin{equation}\label{rel-semi}
{\g}_1 =  \l_1 \quad \mbox{ and } \quad {\g}_{j+1} - n_j {\g}_{j}
= \l_{j+1} -  \l_{j}, \quad \mbox{ for } \quad j= 1, \dots, g-1,
\end{equation}
together with the semigroup $\Gamma := \Z^d_{\geq 0} +  \g_1
\Z_{\geq 0} + \cdots + \g_g \Z_{\geq 0} \subset \Q^d_{\geq 0}$
are associated with the q.o. polynomial $f$.  We set $ \g_0 = 0 \in M$. 
\end{definition}

\begin{lemma}  \label{gen-semi}
  \cite[Proposition 3.1]{GPSemi}
 We have that    $f(\z_{j-1}) = x^{e_0 \g_{j}} \cdot u_j$,
where $u_j$ is a unit in $\C \{x^{1/e_0} \}$, for $j=1, \dots, g$.
\end{lemma}

\begin{definition}  \cite[Definition 4.1]{GPSemi} 
A q.o.~hypersurface $(S', 0) \subset (\C^{d+1}, 0)$ is a $(j-1)$-th \textit{semi-root} of $(S,0)$ if
it is parametrized by a q.o. branch $\t$ with $j-1$
characteristic exponents and such that
 $f_{|y =\t} = x^{e_0 \g_{j}} \cdot u_j$, where $u_j$ is a unit in $\C \{x^{1/e_0} \}$, for $j=1, \dots, g$.
\end{definition}

By Lemma \ref{gen-semi} the germs $(S_j,0)$, for $j=0, \dots, g$
are semi-roots of $(S,0)$. The semi-roots play an important role
in the definition of the toric embedded resolution of the germ
$(S, 0)$ below.

\begin{lemma} \label{semi-relation}
The order of $\g_{j} + M_{j-1}$ in the finite abelian
group $M_j/M_{j-1}$ is equal to $n_j$, for $j = 1, \dots,g$.  We
have the inequality
\begin{equation} \label{preforma}
n_j\gamma_{j} < \g_{j+1},  \textrm{ for }   j=1, \dots, g.
\end{equation}
\end{lemma}

\begin{definition} \label{integral}
We denote by $r_j$ the integral length of the vector $n_j \g_j$ in
the lattice $M_{j-1}$, for $j=1, \dots, g$.
\end{definition}
The integer $r_j$ is the  greatest common divisor of the
coordinates of the vector
   $n_j \g_j$ in terms of a basis of the lattice $M_{j-1}$.

\begin{lemma} \label{gcd}
We have that  $\gcd(r_j, n_j)=1$ for $j=1, \dots, g$.
\end{lemma}
\textit{Proof}. This follows from the first statement of Lemma \ref{semi-relation}.
\hfill $\square$

\begin{lemma}\label{P4eqexp} Set $\l_0 := 0 \in M$. We have the equality
$
\sum_{i=0}^j e_{i}(\lambda_{i+1} - \lambda_i) = e_j\gamma_{j+1}$,
for $0 \leq j \leq g-1$.
 \end{lemma}
\textit{Proof.}  We proceed by induction in $j$. For $j=0$ we have
that  $e_0 \lambda_1 = e_0 \gamma_1$ because $\lambda_1=\gamma_1$.
Suppose that it is true for  $j$. We deduce the result from the
equalities:
\[
\begin{array}{c}
\sum_{i=1}^{j+1} e_{i}(\lambda_{i+1} - \lambda_i) =
\sum_{i=1}^j e_{i}(\lambda_{i+1} - \lambda_i) + e_{j+1}(\lambda_{j+2} - \lambda_{j+1}) =
\\
e_{j}\gamma_{j+1}  + e_{j+1}(\lambda_{j+2} - \lambda_{j+1}) = e_{j+1} ( n_{j+1}\gamma_{j+1} +
(\lambda_{j+2} - \lambda_{j+1}) )= e_{j+1}\gamma_{j+2}.
\end{array}
\]
\hfill $\square$

\begin{notation}  \label{basis}
    The canonical basis $\varepsilon_1, \dots, \varepsilon_d$ of $\Q^d$ is also
a basis of the lattice $M_0$ and a minimal set of generators
of the semigroup  $\Z^d_{\geq 0} \subset M_0$.  The dual basis $\epsilon_1, \dots, \epsilon_d$
 of the dual lattice $N_0$ spans a regular cone $\r$ in
$N_{0,\R}$. It follows that $\R^{d}_{\geq 0} = \r^\vee$ and
$\Z^d_{\geq 0} = \r^\vee \cap M_0$, where $\r^\vee $ denotes the
dual cone of $\r$ (see Notations of Section \ref{sec-toric}). We
denote by $\langle, \rangle : N_\R\times M_\R \to \R$ the duality
pairing associated to the dual vector spaces $N_\R$ and $M_\R$.
\end{notation}

\begin{notation}\label{P4coordnulas}
If $ 1\leq j \leq g$ we denote by   $\ell_j$ the number of nonzero
coordinates of the characteristic exponent $\lambda_{j}$, with
respect to the canonical basis $\varepsilon_1, \dots,
\varepsilon_d$ of $\Q^d$. We set $\ell_0:
=0$.
\end{notation}

\begin{remark}
              We have that $0 = \ell_0 <  \ell_1 \leq \ell_2 \leq \cdots \leq \ell_g$.
\end{remark}

\begin{notation}       \label{cones}
We denote by $\epsilon_1,  \dots, \epsilon_{d+g+1}$ the canonical
basis of the lattice $\Z^{d+g+1}$. We consider the cone  $\r$
embedded as $\r_0:= \r \times \{ 0 \} \subset \R^d \times
\R^{g+1}$ (see Notation \ref{basis}).
The following cones in  $\R^{d+g+1}$ are $d$-dimensional and simplicial
\[
\begin{array}{c}
 \rho_j := \{ a + \sum_{i=1}^j
 \langle a, \gamma_i
\rangle \epsilon_{d + i}  +  \sum_{i=j+1}^{g+1}    n_j \cdots n_{i-1} \langle a, \gamma_j \rangle \epsilon_{d + i}
\mid
a \in \rho \}, \mbox{ for }    j= 1,\dots, g.
\end{array}
\]
We denote by $\nu_i^{(j)}$ the primitive integral vector in the edge of the cone $\r_j$
defined by the vector 
$\epsilon_i 
     + \sum_{r=1}^j
 \langle \epsilon_i, \gamma_r
\rangle \epsilon_{d + r}  +  \sum_{r=j+1}^{g+1}    n_j \cdots n_{r-1} \langle \epsilon_i , 
\gamma_j \rangle \epsilon_{d + r} \in \Q^{d+g+1}$.
We introduce also the following cones of dimension $d+1$
 \[
\sigma^+_j :=  \rho_j + \R_{\geq 0}  \epsilon_{d+j}, \quad \sigma^-_j := \rho_{j-1} + \rho_j,
\quad \mbox{ for } j =1,\dots, g, \quad \mbox{ and }  \quad  \sigma_{g+1} := \rho_g +
\R_{\geq 0}  \epsilon_{d+g+1}.
\]
\end{notation}

 We introduce a fan $\Theta$ in $\R^{d+g+1}$ associated to the q.o.~ polynomial $f$.
This fan appears in  \cite[Section 5.4 ]{GP-Fourier}, in connexion
with the  definition of a toric embedded resolution of the q.o.
hypersurface $(S,0) \subset (\C^{d+1})$.

\begin{definition}  \label{complex}
We set $\Theta_j := \{ \s_j^+, \s_j^-, \r_j \}$ for $j =1, \dots, g$ and $\Theta_{g+1} := \{ \s_{g+1} \}$.
The set $\Theta$ consisting of the faces of the cones in $\cup_{j=1}^{g+1} \Theta_j$, is a fan
associated to the q.o.~polynomial $f$.
\end{definition}

\begin{definition} \label{xi}  The following functions are linear forms on $\Z^{d+g +1}$
\[
\eta (k) = k_{d+g+1},  \  \xi_1 (k) :=    \sum_{i=1}^{d+1} k_i , \ \mbox{ and }   \
\xi_j ( k )   :=    \sum_{i=1}^d k_i + \sum_{i=1}^{j-1} (1 - n_i) k_{d+i}  + k_{d+j}, \  \mbox{ for } 
j =2, \dots, g+1.
\]

\end{definition}

\begin{lemma}\label{not00}
If $k = \sum_{i=1}^{d+g+1} k_i \epsilon_i \in \Z^{d + g +1} $ is a
primitive integral vector generating an edge of some cone $\theta
\in \Theta_j$ for some $1 \leq j \leq g +1$ then $\xi_j(k) \geq
0$. If in addition $k \in \r_{j-1}$ then we get $\xi_{j}  (k) = \xi_{j-1} (k)$.
If $k_{d+g +1} =0$ then we have that $\xi_j(k) =1$, morever in
this case either  $k = \epsilon_{d+j+1}$ or   $k = \epsilon_i$ for some
$1 \leq i \leq d$ such that  the $i$-th coordinate of the
characteristic exponent $\l_j$ is equal to zero.
\end{lemma}
\textit{Proof}.  Let us check that $\xi_j(k) \geq 0$. The claim is
obvious for $j=1$ because the definition of $\xi_1(k)$. If  $j
>1$ and $k \in \rho_{j}$ then  we get that $k_{d+i+1} \geq  n_i k_{d+i}$
for $1\leq i \leq j-1$ and $k_{d+i+1} = n_i k_{d+i}$ for $j \leq i
\leq g$ because of the definition of $\rho_j$ and Lemma \ref{semi-relation}.
In this case we get
that  $\xi_j(k) = \sum_{i=1}^{d+1} k_i + \sum_{i=1}^j (k_{d+i+1} -
n_i k_{d+i}) \geq 0$. If $k \in \rho_{j-1}$ we have that
$k_{d+j+1}=n_j k_{d+j}$ hence
 $\xi_j(k) = \sum_{i=1}^{d+1} k_i + \sum_{i=1}^{j-1} (k_{d+i+1} - n_ik_{d+i})
 = \xi_{j-1}(k) \geq 0$. Finally, if $k = \epsilon_{d+j+1}$  and $1 \leq j \leq g$
 we get that
 $\xi_j(\epsilon_{d+j+1}) =1$.

If $k \ne  \epsilon_{d+j+1}$,  $k_{d+g+1} =0$ and if
 $\nu := (k_1,\dots,k_d)$ we get that
$\langle \nu, \g_j \rangle =0$ (see  Notation \ref{cones}),  that is,  the vector
$\nu \in N_0$ is orthogonal to $\g_j$. Since $ \nu$ defines an edge of $\r$
it follows that $k = \epsilon_i$ for some $ 1 \leq i \leq d$ such
that the i-$th$ coordinate of $\g_j$ is equal to zero. Then the
conclusion follows by  relations (\ref{rel-semi}). \hfill
$\square$

\section{Statement of main results}    \label{results}

In this section we introduce the main results of the paper.

We recall first some basic facts on generating functions.
The \textit{generating function} of a set
$A \subset \Z^{m}$  is the formal power series $\sum_{u \in A}
x^u$.
\begin{lemma}\cite[Chapter 4, Theorem 4.6.11]{stanley} \label{gf} If  $\theta \in
\R^m_{\geq 0}$ is a rational polyhedral cone for the lattice
$\Z^m$ the generating function $G_{{\theta}}$ of the set
$
\mathrm{int} ({\theta}) \cap \Z^m$ is a rational function of
the form $ G_{{\theta}} = {P_\theta}{\prod_{v} (1 - x^{v})^{-1}}$,
where  $v$ runs through the set of  primitive integral vectors in
the edges of the cone $\theta$ and $P_\theta$ is a polynomial with
integral coefficients and exponents in $ \theta \cap \Z^{m}$. If in addition
the cone $\theta$ is simplicial,  we denote by $v_1, \dots, v_r
\in \Z^m$ the primitive integral vectors in the edges of $\theta$
and by $D_\theta$ the set $D_\theta:= \{ \nu = \sum_{i=1}^r a_i
v_i \in \mathrm{int} ({\theta}) \cap \Z^m \mid 0 < a_i \leq 1
\}$; then we get that $P_\theta =  \sum_{\nu \in D_\theta} x^\nu$.
\end{lemma}

\begin{lemma}\label{simplicial} With the hypothesis and notations
of Lemma \ref{gf} we  assume in addition that the cone $\theta$ is
simplicial of dimension $r$. We consider two $\Z$-linear maps
$\kappa, \iota : \Z^m \rightarrow \Z$ such that for any primitive
integral vector in an edge of the cone $\theta$ we have $\iota (v)
\geq 0$ and if $\iota(v)=0$ then $\kappa(v) =1$. We denote by
$\tilde{G}_\theta$ the image of $G_\theta$ by the monomial map
$x^\nu \mapsto \L^{-\kappa(\nu)}T^{\iota(\nu)}$ and  by $r_0$ the
number of rays $\R_{\geq 0} v$  of $\theta$ such that $\iota(v) =
0$.  We have that
$\lim_{T\rightarrow \infty} \tilde{G}_\theta = (-1)^{{r-r_0}} (\L-1)^{{-r_0}}$.
\end{lemma}
\textit{Proof.} By Lemma \ref{gf} we get that $\tilde{G}_\s =
(\sum_{{\nu} \in D_\theta } \L^{-\kappa({\nu})} T^{\iota({\nu})})
\prod_{i=1}^r (1-\L^{-\kappa(v_i)}T^{\iota(v_i)})^{-1} $. The
vector ${\nu_0} := \sum_i v_i$ belongs to $D_\theta$ and verifies
that  $\iota(\nu) < \iota(\nu_0)$ for all $\nu_0 \ne \nu \in
D_\theta$. Lemma \ref{P4limiteremarkb} implies that $\lim_{T
\rightarrow \infty} \tilde{G}_\theta = \lim_{T \rightarrow \infty}
\prod_{i=1}^r \L^{-\kappa(v_i)}T^{\iota(v_i)}
(1-\L^{-\kappa(v_i)}T^{\iota(v_i)})^{-1}$. The result follows from
Lemma \ref{P4limiteremark}. \hfill $\square$

Let $\theta \subset \R^m$ be a rational simplicial cone of dimension $r$.
We denote by $v_1, \dots, v_r \in  \Z^m$ the primitive vectors generating the edges of the cone
$\theta$ and by $\Z (\theta \cap \Z^m)$ the sublattice of $\Z^m$ spanned by
$\theta \cap \Z^m$.  The lattice $\oplus_{i=1}^r \Z v_i$ is a sublattice of finite index
$\mathrm{mult} (\theta)$  of $\Z (\theta \cap \Z^m)$. Notice that $\mathrm{mult} (\theta)$ is equal to
the cardinality of the set $D_\theta$ introduced in Lemma \ref{simplicial}.

\begin{notation}       \label{ztop-not}
 If $\theta \in \Theta_j$, for $1 \leq j \leq g+1$,  is a simplicial cone
 we denote by 
$J_\theta (f, s)$ the element of $\Q (s)$ given by
$
 J_\theta (f, s) := { \mathrm{mult} (\theta) }{\prod_{v} (\xi_j (v)  + \eta(v)  s)^{-1}}$, 
where $v$ runs through the primitive vectors defining the edges of $\theta$.
If $\theta \in \Theta_j$ is not simplicial then $\theta = \s_j^-$  and we
set
$
 J_{\s_j^-} (f, s) := (1 + n_j \dots n_g s)^{-1} (  J_{\r_{j-1}} (f, s)     - J_{\s_j^+} (f,s) )
$, 
where $J_{\r_{0}} (f, s) :=1$.
\end{notation}

\begin{notation} \label{gf2}
We use the notations introduced in Lemma \ref{gf}. Let $\theta$ be
a cone in $\Theta_j$, for some $j=1, \dots, g+1$. We denote by
$\tilde{P}_\theta$ the image of $P_\theta$ by the monomial mapping
\begin{equation} \label{polymap}
 \Z [x_1, \dots, x_{d+g+1} ] \to \Z[\L, T], \quad x^\nu \mapsto \L^{-\xi_j(\nu)} T^{\eta(\nu)}.
\end{equation}
We define  the rational function
\begin{equation} \label{stheta}
 S_\theta :=    \tilde{P}_\theta \prod_{v}   (1 - \L ^{-\xi_j (v) } T^{\eta(\nu)}) ^{-1},
\end{equation}
where  $v$ run through the primitive
integral vectors in the edges of $\theta$.
 \end{notation}

\begin{remark}
The rational function $S_\theta$ is well defined. By Lemma    \ref{not00} for any vector $v$ as above 
we have that $(\xi_j (v), \eta(v)) \ne
(0,0)$.
\end{remark}

\begin{definition}          \label{csigma}
 If $\theta \in \Theta_j$  for $ 1\leq j \leq g$, we
define the elements $c_1(\theta) \in \kmon$ (resp. $c(\theta) \in
\kvar$) as
\[ \left \{
\begin{array}{lcl}
[\mu_{n_{j}e_{j}}]  \quad \mbox{  (resp. } \L -1  \mbox{ ) }  &
{\rm if } &\theta = \s^+_{j},
\\
\protect [\{(x,{y}) \in (\C^*)^2 | ({y}^{n_{j}} -
x^{r_{j}})^{e_{j}}=1\}]     \quad \mbox{  (resp. } (\L -1)(\L -2)
\mbox{ ) } & {\rm if} &\theta = \rho_{j},
\\
\protect [\mu_{r_{j}e_{j}}]  \quad  \mbox{  (resp. } \L -1 \mbox{ ) }  &{\rm if} &\theta = \s^-_{j}.
 \end{array}
\right. \]
 If $\theta = \s_{g+1} \in \Theta_{g+1}$ we set
$c_1(\theta) := 1$ (resp. $c(\theta) := \L -1$).
\end{definition}

\begin{theorem}\label{main-zeta} Let $f : (\C^{d+1},0) \to
(\C,0)$ be the complex analytic function defined by an irreducible
q.o.~polynomial. The  zeta functions $Z(f,T)_0$, $ Z^{{\rm naive}}
(f , T )_0 $, $Z_{\mathrm{top}} (f, s)_0$
 and the motivic Milnor fibre
$\mathcal{S}_{f,0}$ of the $f$ at $0$ are determined by the
embedded topological type of the q.o.~hypersurface germ $(S,0)
\subset (\C^{d+1}, 0)$ defined by $\{ f =0\}$. We have the
following formulas in terms of the characteristic exponents of the
q.o. branch $\z$ parametrizing $(S,0)$:
\[
\begin{array}{c}
Z_{\mathrm{top}}(f, s)_0 =    J_{\s_{g+1}} (f, s) +

\displaystyle{\sum_{j=1}^{g}  }
(J_{\s_j^+}  (f, s)  - J_{\r_j}  (f, s)  +  J_{\s_j^-}  (f, s))
\\
 Z^{{\rm
naive}} (f , T )_0 = \displaystyle{\sum_{j =1}^{g+1} \sum_{\theta \in \Theta_j}}
c(\theta) ({\L}-1)^{\dim \theta -1}  S_\theta,
\\
       Z (f , T )_0  =   \displaystyle{\sum_{j =1}^{g+1} \sum_{\theta \in \Theta_j}}
c_1(\theta) ({\L}-1)^{\dim \theta -1} S_\theta,
\\

\mathcal{S}_{f,0}  = (1 - \L)^{\ell_g} + \displaystyle{\sum_{{j}=1}^g \big (  
(1 - \L)^{\ell_j -1} ( c_1(\rho_j) + c_1(\s^+_j) ) + (1 - \L)^{\ell_{j-1}} c_1(\s^-_j)  \big )}
\end{array}
\]
 \end{theorem}

\begin{remark}\label{P4remarkmainthm}
For $d=1$ we recover a reformulation of results of Guibert (see
\cite{Guibert}, Proposition 3.3.1) and Artal et al.~(see \cite{ACNLM-AMS}, Theorem 4.3, 4.9 and 5.3).  
Notice  that Guibert asssumes that the line $\{ x=0 \}$ is non-tangent to the curve.
\end{remark}

\begin{remark}      
If $g=1$ then $f$ is Newton non-degenerate and the formula above
for $Z_{\mathrm{top}}(f, s)_0$ reformulates in this particular
case the one given by Denef and Loeser (see  Th\'eor\`eme 5.3
\cite{DL-JAMS}).
\end{remark}

\begin{remark}        \label{comp}
If $d \geq 2$ the phenomenon of false reducibility, discovered in \cite{MGV}, may lead 
to inaccurate formulas for 
the zeta functions of q.o.~singularities in \cite{ACNLM-AMS} (see Example \ref{false}). 
In Section \ref{comparison} we compare the 
method of Newton maps used in \cite{ACNLM-AMS} with the toric methods used in this paper. We check that 
the sets of candidate poles we get below (see Corollary \ref{poles})  is the same as the one given in \cite{ACNLM-AMS}.
\end{remark}

\begin{corollary}
\label{main-hs} With the hypothesis of Theorem \ref{main-zeta} the
Hodge-Steenbrink spectrum of $f$ at $0$ is determined by the
characteristic exponents. We have the formula:
\[
Sp' (f, 0) = (-1)^d ( \sum_{j=1}^g  \big ( (1-t)^{\ell_j - 1} (
P_{n_j,r_j}^{(e_j)} + P_{n_j e_j} ) + (1-t)^{\ell_{j-1}} P_{r_j
e_j} \big) + (1-t)^{\ell_g} -1 ). 
\]
\end{corollary}
\textit{Proof}. It follows  by using the formula for the motivic
Milnor fiber in Theorem \ref{main-zeta}, Notation \ref{notxhsp},
Lemma \ref{P4ejemsp}  and Corollary \ref{cornotxhsp}. \hfill
$\square$

 By Corollary \ref{main-hs} the Hodge-Steenbrink spectrum 
of an irreducible q.o. polynomial $f$ is determined by the embedded topological type of 
the germ of singularity defined by  $f=0$. 
\begin{remark}\label{P4hspcurvasSaito} For $d=1$ we recover
Saito's Formula  (\cite{Saito00} Theorem 1.5 or \cite{Guibert},
Corollary 3.4.1). See Remark \ref{sp-normal}.  
In this case the Hodge-Steenbrink spectrum of
$f$ at $0$ classifies the embedded topological type of the
analytically \textit{irreducible} germ of plane curve $(S, 0) \subset
(\C^2,0)$ defined by $f=0$. 
However, there are examples of functions $f_i :
(\C^2, 0) \to (\C,0)$ defining non-irreducible germs 
at the origin such that the Hodge-Steenbrink spectrum of
$f_i$ at $0$ is the same for $i=1,2$ meanwhile the germs defined
by $f_i=0$, $i=1, 2$ have different embedded topology, see
\cite{Steenbrink-Asterisque, SSS}.  \end{remark}
\begin{remark}
If $d \geq 2$ the  Hodge-Steenbrink
spectrum 
associated to 
an irreducible q.o. polynomial $f \in \C \{ x_1, \dots, x_d\}[y]$
 is not a complete invariant of the embedded topological type of 
the germ $f =0$ at the origin. For instance, 
the q.o.~polynomials
$f_1 := y^3 - x_1 x_2$ and $f_2 := y^3 - x_1^2x_2$ define
irreducible germs with different embedded topological types
meanwhile Corollary \ref{main-hs} implies their Hodge-Steenbrink
spectrums coincide. \end{remark}

\begin{remark}
If $d \geq 3$ the Seifert form on the middle homology of the Milnor fiber of an hypersurface singularity 
$F: (\C^{d +1 }, 0) \to (\C, 0)$ defining an isolated singularity at the origin 
determines the topological type of the singularity (see \cite{Durfee}, Theorem 3.1).
However, there are examples of germs   $f_i  :
(\C^2, 0) \to (\C,0)$  defining non-irreducible germs 
at the origin such that they have isomorphic integral Seifert form 
but the germs    defined
by $f_1 =0$ and $f_2=0$ have different embedded topology and different spectrum 
(see \cite{DuBois-Michel}). 
The germs
$F_i : (\C^4 , 0) \to (\C,0)$ given by  $F_i = f_i(x_1,x_2)  + x_3^2 + x_4^2$, $i=1, 2$, 
have isomorphic Seifert forms    (see \cite{Sakamoto}). 
It follows that the  germs $F_i= 0$ at $0 \in \C^4$ have the same embedded topological type, meanwhile 
the application of Thom Sebastiani Theorem for the spectrum implies that 
 $F_1$ and $F_2$ have different spectrum at $0$ (see \cite{Denef-LoeserBarca} Theorem 5.1 or \cite{Varchenko}). 
In this example 
the topological zeta functions of $F_i$ at $0$, $i=1, 2$, 
 are also different (see \cite{ACNLM-JLMS}). 
\end{remark}

\begin{definition}  \label{special}       (see Notations \ref{cones} and 
\ref{lamb-coord}) 
The vector $\nu_i^{(j)} \in \r_j$ is special if  for any $1 \leq i' \leq d$ the
following equalities 
\[ \l_{j', i'} = 0 \mbox{ for }   1 \leq j' < j, \mbox{ and }  \l_{j, i'} 
= 1/n_j,\]
are verified iff $i' = i$ , 
and in addition we have that $\l_{j+1, i} > \l_{j, i} $ if $j < g$. 
\end{definition}

\begin{example} 
If $\{x_i =0\} \cap S$ is the unique coordinate hyperplane section of 
$(S, 0)$  which is not contained in the singular locus of $(S,0)$ 
but it is contained in the critical locus of the quasi-ordinary projection $\pi$,
 then $\nu_i^{(j)}$ is a special vector. This is a consequence of  
the characterization of the singular locus of a q.o.~hypersurface (see \cite[Therorem 7.3]{Lipman}).
\end{example}

\begin{corollary}     \label{poles}     (compare with  \cite[Definition 5.5 and Proposition 5.7]{ACNLM-AMS}) 
With the hypothesis of Theorem \ref{main-zeta} and the Notation \ref{cones}
there exists a polynomial $Q \in \Z[s]$ such that 
$$Z(f,s)_0 = Q (1 + s)^{-1}  \prod_{j = 1}^g  
\prod_{\nu_i^{(j)} \mathrm{non-special} }  (\xi_j ( \nu_i^{(j)} ) + \eta(\nu_i^{(j)} ) s )^{-1}.$$
There exists $R \in \Z [ \L^{\pm 1} ] [T]  $ such that 
\[ Z^{{\rm naive}} (f , T )_0 = R (1 - \L^{-1} T) ^{-1} 
\prod_{j = 1} ^{g} 
\prod ^{i=1, \dots, d}_{\nu_i^{(j)} \mathrm{non-special} } 
( 1  - \L^{ -\xi_j (\nu_i^{(j)})} T^{\eta( \nu_i^{(j)}) } )^{-1}  .\]
\end{corollary}
Notice that if $\nu_i^{(j)} = \epsilon_i $ then it does not contribute to a pole of these zeta functions. 

\begin{remark}   \label{comp-1}
We prove Corollary \ref{poles} in Section \ref{comparepoles}. 
We can use Corollary \ref{poles}, Remark \ref{comp-4} and the arguments 
from \cite[Chapter 6]{ACNLM-AMS} to complete the proof of the monodromy conjecture in this case. 
The argument consists of showing that the candidate poles 
are related with monodromy eigenvalues associated to the restriction of $f$ to suitable  plane sections of $\C^{d+1}$ 
in a neighbourhood of $0$.   
\end{remark}

In order to prove these results we use an 
embedded resolution $\Pi' \colon  Z' \to \C^{d+1}$ of $(S,0) \subset (\C^{d+1}, 0)$.
This resolution is a composition of toric modifications which we recall in Section \ref{resolution}
following \cite{GP-Fourier}.
We remark a relation between the denominators of the zeta functions in Corollary \ref{poles}
and certain exceptional divisors of the modification $\Pi'$. 

\begin{remark} \label{vanish-2}
For any $1 \leq i \leq d$ and $1 \leq j \leq g$ 
there exists an irreducible divisor $E_i^{(j)} \subset Z'$ such that
the order of vanishing of the jacobian of $\Pi'$ (resp.~of $f \circ \Pi'$) 
along an  irreducible divisor 
$E_i^{(j)} \subset Z'$ is equal to 
$\xi_j (\nu_i^{(j)}) - 1$  (resp. $\eta (\nu_i^{(j)}$).
 Notice that 
$E_i^{(j)}$ is contained in the critical locus of $\Pi'$, that is $E_i^{(j)}$ is exceptional divisor for $\Pi'$, 
  if and only if $\xi_j (\nu_i^{(j)}) > 1$. These facts are 
  consequence of Remarks \ref{vanish} and \ref{vanish-3} below.
\end{remark}

\section{Some notions on toric geometry and Newton polyhedra}     \label{sec-toric}

In this section we introduce the notations and basic results of
toric geometry which we use in this paper (see \cite{Fulton,Ewald,
Kempf,Oda}).

If $N \cong \Z^{d}$ is a lattice we denote by $N_\R$ the vector
space spanned by $N$ over the field $\R$. In what follows a {\it
cone} mean a {\it rational convex polyhedral cone}: the set of non
negative linear combinations of vectors $v_1, \dots, v_r \in N$.
The cone $\t$ is {\it strictly convex} if it contains no lines, in
such case we denote by $0$ the $0$-dimensional face of $\t$.  The
dimension of $\t$ is the dimension the linear span of $\t$ in
$N_\R$. The relation $\theta \leq \t$ (resp. $\theta < \t$)
denotes that $\theta$ is a face of $\t$ (resp. $\theta \ne \t$ is
a face of $\t$). The cone $\t$ is \textit{simplicial} the number
of edges of $\t$ is equal to the dimension of $\t$. The cone $\t$
is \textit{regular} if it is generated by a sequence of vectors $v_1, \dots, v_r$ 
which is contained in a basis of the lattice $N$.

We denote by $M$ the dual lattice, i.e., $M =\mathrm{Hom} (N,\Z)$
and by $\langle, \rangle : N \times M \to \Z$ the duality pairing
between the lattices $N$ and $M$.  
The {\em dual cone}  $\t^\vee$ (resp. {\em orthogonal cone}
$\t^\bot$) of $\t$ is the set $ \{ w  \in M_\R\ |\ \langle u, w
\rangle \geq 0,$  (resp. ${ \langle u, w \rangle } = 0$) $ \; \forall
u \in \t \}$. We denote the \textit{relative interior} of the cone $\t$ by 
 $\stackrel{\circ}{\t}$ (or also by $ \mathrm{int}(\t)$ when 
it is more convenient for typographical reasons).

A {\em fan} $\Sigma$ is a family of strictly convex
  cones  in $N_\R$
such that for any $\s, \s' \in \Sigma$ we have $\s \cap \s' \in
\Sigma$ and if $\t \leq \s$ then  $\t \in \Sigma$.  The {\em
support} of the fan $\Sigma$ is the set $|\Sigma | := \bigcup_{\t
\in \Sigma} \t \subset N_\R$. We say that a fan $\Sigma'$ is a
{\it subdivision} of the fan $\Sigma$ if both fans have the same
support and if every cone of $\Sigma'$ is contained in a cone of
$\Sigma$. The fan $\Sigma'$ is a \textit{regular subdivision} of
$\Sigma$ if all the cones in $\Sigma'$ are regular and  $\theta
\in \Sigma'$ for any regular cone $\theta \in \Sigma$. Every fan
admits a regular subdivision.

Let ${\s}$ be an strictly convex cone rational for the lattice
$N$. By Gordan's Lemma the semigroup ${\s}^\vee \cap M$ is
finitely generated. We denote by $\C[{\s}^\vee \cap M] := \{ \sum_{\mathrm{finite}}  \ c_\d x^\d \mid
c_\d \in \C, \d \in {\s}^\vee \cap M \}$ the
semigroup algebra of $\s^\vee \cap M$ with complex coefficients.
The toric variety $Z_{{\s}} := \mathrm{Spec} \C [ {\s}^\vee \cap
M]$, denoted also by $Z_{{\s}, N}$ or $Z^{{\s}^\vee \cap M}$, is
normal. The torus $ T_N:= Z^{M}$ is an open dense subset of
$Z_{\s}$, which acts on $Z_{\s}$ and the action extends the action
of the torus on itself by multiplication. There is a one to one
correspondence between the faces $\t $ of ${\s}$ and the orbits
$\orb_\t$ of the torus action on $Z_{\s}$,
 which reverses the
inclusions of their closures.

The
 affine varieties $Z_\s$ corresponding to cones in a fan $\Sigma$
 glue up to define a {\em toric variety} $
 Z_\Sigma$.
A subdivision $\Sigma'$ of  the fan $\Sigma$ defines a {\em toric
 modification} $ \phi : Z_{\Sigma '}
 \rightarrow   Z_\Sigma$.
If $\Sigma'$ is a regular subdivision of $\Sigma$ then the
modification $\phi$ is a resolution of singularities of
$Z_{\Sigma}$.

If $\theta$ is $d$-dimensional strictly convex cone then its orbit
$\orb_\theta$ is reduced to a closed point called the origin $0$
of the toric variety $Z_\theta$. We denote by $ \C [[\theta^\vee \cap
M]]$ the ring of formal power series with exponents in $\t^\vee
\cap M$. This ring is the completion with respect to the maximal
ideal of the ring $\C \{ \theta^\vee \cap M \} \subset  \C [[\theta^\vee
\cap M]] $ of germs of holomorphic functions at $0 \in Z_\t$.
The \textit{Newton polyhedron}  $\Newton (\psi)$ of a nonzero series $\psi = \sum
c_\d x^{\d} \in \C [[\theta^\vee \cap M]]$ is the convex hull of the Minkowski sum $\{ \d \mid c_\d \ne 0
\} + \theta^\vee$. The support function of the polyhedron $\Newton
(\psi)$ is the piecewise linear function
\[
 h: \theta \to \R, \quad \nu \mapsto   \inf
\{    \langle \nu, \omega \rangle \mid \omega \in {\mathcal N}
(\psi )  \}.
\]
A vector $\nu \in \theta$  defines the face  $ {\mathcal F}_\nu :=
\{ \omega \in {\mathcal N} (\psi ) \mid \langle \nu, \omega
  \rangle   = h (\nu ) \}$
of  the polyhedron ${\mathcal N} (\psi)$. All faces of ${\mathcal
N} (\psi)$ are of this form, in particular, the compact faces are
defined by vectors $\nu  \in \stackrel{\circ}{\theta}$. The set
consisting of the cones $ \s( {\mathcal F} ) := \{ \nu \in \theta
\; \mid \langle \nu , \omega \rangle   = h(\nu),
 \; \forall \omega \in {\mathcal F}\}$,
for ${\mathcal F}$ running through the faces of ${\mathcal N}
(\psi)$ is a fan subdividing the cone $\theta$.
     If $\nu \in \theta$ we denote by $\psi_{| \nu}$ the symbolic restriction
$
 \psi_{| \nu} := \sum_{\d \in     {\mathcal F}_\nu} c_\d x^\d$
 of the series $\psi$ to the face   $
{\mathcal F}_\nu$ of its Newton polyhedron.   If $\nu \in
\stackrel{\circ}{\theta}$ then
 $\psi_{| \nu} $ is a polynomial, i.e., $\psi _{| \nu} \in \C [\theta^\vee \cap M]$.

\section{Toric embedded resolutions of q.o.~hypersurfaces}   \label{resolution}

In this section we give an outline of the construction of a
toric embedded resolution of a q.o. hypersurface following \cite{GP-Fourier}.

\begin{notation}
Recall that we denote by $\r$ the cone $\R^{d}_{\geq 0} \subset
N_\R$ (see Notation \ref{basis}). We denote by  $N'_j$ the lattice
$N_j \times \Z$,  and by $M'_j$ its dual lattice $M_j \times \Z$,
for $j=0, \dots, g$. The cone  $\r' := \r \times \R_{\geq 0}$ is
rational for the lattice $N'_j$.
\end{notation}

We keep notations of Section \ref{qoh}, that is,
 $\z$ is  a q.o.~branch with characteristic exponents $\l_1, \dots, \l_g$.
The q.o.~branch $\z$  parametrizes a q.o.~hypersurface $(S, 0)$
defined by a q.o. polynomial $f \in \C\{ x_1, \dots, x_d \}[y]$.
We assume, making a change of coordinates, that the coordinate $y$ is equal to $f_0$
(see  Definition \ref{semi}).

With this assumptions the Newton polyhedron $\Newton (f) \subset
M'_\R$ has only one compact edge $\mathcal E_1$ with vertices $(0,
e_0)$ and $(e_0 \l_1, 0)$.  It defines a subdivision $\Sigma_1$ of
the cone $\r'$, with only two $(d+1)$-dimensional cones
$\bar{\s}_1^+$ and $\bar{\s}_1^-$ which intersect along the
$d$-dimensional face $\bar{\r}_1$. The support function of the
Newton polyhedron $\Newton (f)$ is defined for $(\nu, r) \in \r'$ by:
\[
h^{(1)} (\nu, r) = \left \{
\begin{array}{lcl}
e_1 \langle \nu,  n_1 \lambda_1 \rangle & {\rm if } &  (\nu,r) \in
\bar{\s}^+_1,
\\
e_1 n_1  r & {\rm if } &  (\nu,r) \in \bar{\s}^-_1,
\end{array}
\right.
\]
that is,  the cone  $\bar{\s}_1^-$  (resp. $\bar{\s}_1^+$)
corresponds to the vertex $(0, e_0)$ (resp. $(e_0 \l_1, 0)$) of
$\Newton (f)$. The polynomial $f$ is of the form,
\begin{equation} \label{forma-f}
f = (y^{n_1} - c_1 x^{n_1 \lambda_1})^{e_1} + \cdots
\end{equation}
where $c_1 \in \C^*$ and  the terms which are not written lie
above the edge $ \mathcal E_1$.

\begin{remark}
We choose the notation   $\bar{\s}_1^+$,  $\bar{\s}_1^-$ and
$\bar{\r}_1$ for the cones in $\Sigma_1$ since these cones are
related with the cones $\s_1^+$, $\s_1^-$ and $\r_1$  in
$\Theta_1$ (see Definition \ref{complex} and the proof of Theorem
\ref{trop}).
\end{remark}

We describe the transformation of the hypersurface $(S, 0)$ under
the toric modification \[ \Pi_1 : Z_1 \to Z_0:= \C^{d+1} \]
defined by the subdivision $\Sigma_1$. Notice that $Z_1$ is not a
smooth variety. First we need an elementary lemma to describe the
chart $Z_{\bar{\r}_1, N'_0} \subset Z_1$ associated to the cone
$\bar{\r}_1$  and the lattice $N'_0$.

\begin{lemma}\label{iso-semi}
Let us consider the lattice homomorphism
\[
 \phi_1: M'_0 \to M_1, \quad (\l, a) \mapsto   \l + a \l_1.
\]
We have the following properties:
\begin{enumerate}
 \item[(i)]
The homomorphism $\phi_1$ is surjective. \item[(ii)] $\ker
(\phi_1) \cap  M'_0  =  \bar{\r}_1^\perp \cap M'_0 =  (-n_1 \l_1,
n_1) \Z$.

\item[(iii)] $\phi_1^{-1} (\r^\vee \cap M_1) = \bar{\r}_1^\vee
\cap M'_0$.

\item[(iv)] We have an isomorphism of semigroups $\bar{\r}_1^\vee
\cap M'_0 \to    (\r^\vee \cap M_1) \times \Z$.

\item[(iv)] The dual homomorphism $ \phi_1^* : N_1 \to N'_0$ is
injective and verifies that $\phi_1^* (\nu) = (\nu, \langle \nu,
\l_1  \rangle )$.

\item[(v)] $\phi_1^* (N_1)$ is the lattice generated by
$\bar{\r}_1 \cap N_0'$.

\item[(vi)] $\phi_1^* (\r \cap N_1) = \bar{\r}_1 \cap N_0'$.
\end{enumerate}
\end{lemma}
\textit{Proof.} Recall that by definition we have the equality
$N_1 = \{ \nu \in N_0 \mid \langle \nu, \l_1 \rangle \in \Z \}$.
The details of the proof are left to the reader (see \cite{GP-Fourier} Lemma 17).
 \hfill $\square$

By Lemma  \ref{iso-semi} the term
\[
w_1 := y^{n_1} x^{-n_1 \l_1}
\]
 is a holomorphic
function on the chart  $Z_{\bar{\r}_1, N'_0}$.
\begin{lemma}  \label{iso-chart}    We have the following properties:
\begin{enumerate}
 \item[(i)]
The orbit $\orb (\bar{\r}_1)$ is a  one-dimensional torus embedded as a closed subset in the chart
$Z_{\bar{\r}_1, N'_0}$.
\item[(ii)]
The  coordinate ring of
   $\orb (\bar{\r}_1)$ is
equal to $ \C [ w_1^{\pm 1}]$.
          \item[(iii)]
The hypersurface $\{ w_1 = c \}   \subset Z_{\bar{\r}_1, N'_0}$ is isomorphic to  $Z_{\r, N_1}$
for any $c \in \C^*$.
\item[(iv)]
The chart $Z_{\bar{\r}_1, N'_0}$ is  isomorphic to
$Z_{\r, N_1} \times \C^*$.
\item[(v)] If $o \in \orb (\bar{\r}_1)$ the germ
    $(Z_{\bar{\r}_1, N'_0}, o)$ is isomorphic to  $(Z_{\r, N_1} \times \C, (0,0))$.

\end{enumerate}

\end{lemma}
\textit{Proof.} The first assertion is well-known. The statement
(ii) is consequence of (ii) in     Lemma \ref{iso-semi}. The
assertion (iv) is consequence of the isomorphism of semigroups in
Lemma \ref{iso-semi} (iv). The statement (iii) is consequence of
(iv) and (ii). The assertion (v) is deduced from (iv). \hfill
$\square$

The modification $\Pi_1$ is an isomorphism over the torus  $T_{N_0'}= (\C^*)^{d+1} \subset Z_0$.
The \textit{strict transform}  $S^{(1)}$ of the hypersurface $S$ by $\Pi_1$
is the closure of $\Pi_1^{-1} (S \cap (\C^*)^{d+1})$ in $Z_1$.
It is shown in \cite{GP-Fourier} that    $S^{(1)}$ intersects the exceptional fiber $\Pi_1^{-1} (0)$ only at
one point $o_1$ which belongs to the orbit $\orb (\bar{\r}_1)$ of $Z_1$.
Hence, in the chart $Z_{\bar{\r}_1} = Z_{\bar{\r_1}, N_0'}$  we can factorize
\begin{equation}
\label{factor} ( f \circ \Pi_1 )_{|Z_{\bar{\r}_1}} =  x^{e_1 n_1
\l_1} \left( (w_1 - {c_1})^{e_1} + \cdots    \right),
\end{equation}
where the terms which are not written vanish on $\orb (\bar{\r}_1)$.
The strict transform of $(S, 0)$ by $\Pi_1$ is a hypersurface germ
$(S^{(1)}, o_1) \subset (Z_{\bar{\r}_1, N'_0}, o_1)$ defined by
the vanishing of $ ( x^{-e_1 n_1 \l_1} \cdot f \circ \Pi_1)
_{|Z_{\bar{\r}_1}}$. The intersection $S^{(1)} \cap \orb
(\bar{\r}_1)$,  defined by the equation $
      (w_1 - {c_1})^{e_1} =0
$, consists of the point $o_1$ counted with multiplicity $e_1$. We
get that the strict transform $S_1^{(1)}$ of the semi-root $S_1$
intersects the orbit    $\orb (\bar{\r}_1)$  at the point $o_1$
with multiplicity one. It follows that the function
\[
y_1 :=      x^{-n_1 \l_1} \cdot (f_1 \circ
\Pi_1)_{|Z_{\bar{\r}_1}},
\]
can be taken as a local coordinate at $o_1$, replacing $w_1 -
{c_1}$. By Lemma \ref{iso-chart} (iii) the strict transform
$S_1^{(1)}$ is analytically isomorphic to the germ of toric
variety $Z_{\r, N_1}$ at the origin. Taking into account the isomorphism of 
Lemma 5.4 (iv) we get.

\begin{proposition}\cite[Proposition 19]{GP-Fourier}  \label{covering}
 Let us consider the projection
\[\p_1 : Z_{\r, N_1} \times \C^* \to   Z_{\r, N_1}.\]
The restriction of $\p_{1 | S^{(1)}} : (S^{(1)}, o_1)  \to    (Z_{\r, N_1}, 0)$
defines an unramified finite covering of the torus $T_{N_1} \subset Z_{\r, N_1}$.
\end{proposition}

\begin{remark}\label{r56}
  Proposition \ref{covering} means that
the strict transform $(S^{(1)}, o_1)$ is a germ of \textit{toric
quasi-ordinary} (t.q.o.) hypersurface relative to the base $
Z_{\r, N_1}$ (see \cite{GP-Fourier}, Definition 5 and
\cite{GP-Can}). If the cone $\r$ is regular for the lattice
$N_1$ then   $(S^{(1)}, o_1)$ is a q.o.~hypersurface
in the classical sense. The notions and results in Section
\ref{qoh} generalize to the case of t.q.o.~singularities in the
relative hypersurface case. This applies in particular for the
definition of  characteristic exponents, semigroup, semi-roots,
etc. It is convenient to consider the characteristic exponents of
a (classical) q.o.~branch  $\z$ not only as $d$-tuples of
rational numbers, but as rational vectors inside a \textit{reference
cone} $\r^\vee$, which is regular for a \textit{reference lattice}
$M_0$. Then, if we drop the regularity condition, the statement of
Lemma \ref{expo} characterizes the t.q.o.~branches,  which
parametrize the t.q.o.~singularities in the relative hypersurface
case. See \cite{GP-Fourier}.
The arguments above apply similarly to the semiroots $S_j^{(1)}$.
\end{remark}

\begin{proposition}\cite[Proposition 19]{GP-Fourier}     \label{trans-semi}
 The following statements hold for $j=2,
\dots,g $.
\begin{enumerate}
 \item[(i)] The strict transform  $(S_j^{(1)}, o_1)$ of the semi-root $(S_j,0)$
is parametrized by  a t.q.o.~branch $\z_j^{(1)}$ which
has characteristic exponents
$\l_2 -\l_1, \dots, \l_j - \l_1$  and characteristic integers $n_2, \dots, n_g$,
with respect to the reference cone $\r^\vee$ and reference lattice $M_1$.
\item[(ii)]
The germ  $(S_j^{(1)}, o_1)$ is a $(j-1)$-th semi-root of $(S^{(1)}, o_1)$.
\end{enumerate}
\end{proposition}

If $g > 1$ it follows from Propositions \ref{covering} and \ref{trans-semi} that
the germ $(S_j^{(1)}, o_1)$ is defined by the vanishing of
\begin{equation} \label{fac-j}
 f_j^{(1)} := (x^{- n_1 \dots n_j \l_1} f_j   \circ \Pi_1)_{|Z_{\bar{\r}_1}}.
\end{equation}
The germ $(S_j^{(1)}, o_1)$ is
 parametrized by  setting
$y_1 = \z_j^{(1)}$, for $j=2, \dots, g$. We can expand $f_j^{(1)}$
as an element of $\C \{ (\r')^\vee \cap M_1' \}$ in such a way
that its Newton polyhedron $\Newton( f^{(1)} ) \subset
(M_1)_\R \times \R$ has a unique compact edge $\mathcal E_2$ with
vertices $(e_1 (\l_2 -\l_1), 0)$ and $(0, e_1)$. This expansion of
$f^{(1)}$ is of the form,
\begin{equation} \label{forma-f-1}
f^{(1)} = (y_1^{n_2} - c_2   x^{n_2 (\l_2 - \l_1}))^{e_2} + \cdots
\end{equation}
where $c_2 \in \C^*$ and the terms which are not written
lie above  the edge $ \mathcal E_2$.
The polyhedron    $\Newton( f^{(1)} )$
 defines a subdivision $\Sigma_2$ of $\r'$ with only two
cones of dimension $d+1$. These cones are  $\bar{\s}_2^+$ (resp.
$\bar{\s}_2^-$), which corresponds to the vertex
 $(e_1 (\l_2 -
\l_1), 0)$ (resp. $(0, e_1)$) of the polyhedron $\Newton (f^{(1)}
)$. With respect to these new coordinates  
we consider the toric modification
$
 \Pi_2: Z_3 \to Z_2$,
induced by the subdivision $\Sigma_2$.
We set $w_2 := y_1^{n_2} x^{- n_2 (\l_2 - \l_1)}$.
In the chart $Z_{\bar{\r}_2} = Z_{\bar{\r}_2, N_1'} \subset Z_3$  one can factorize
\begin{equation}
\label{factor1} (f^{(1)} \circ \Pi_2 )_{|Z_{\bar{\r}_2}}=  x^{e_1
(\l_2 - \l_1)} \left( (w_2 - c_2)^{e_3} + \cdots \right),
\end{equation}
where the terms which are not written vanish on $\orb
(\bar{\r}_2)$. The strict transform $S^{(2)}$,  defined then by
the vanishing of $ x^{-e_1(  \l_2 - \l_1) }   \cdot (f^{(1)} \circ
\Pi_2) $, is a germ at the point $o_2 \in   \orb (\bar{\r}_2)$
given by $w_2 - c_2 =0$.

Then, we iterate this procedure, obtaining for $2 < j \leq g$ a
sequence of local toric modifications $\Pi_{j-1}: Z_j \to
Z_{j-1}$, which are described by replacing the index $1$ by $j-1$
and the index $2$ by $j$ above. The toric modification $\Pi_{j-1}$
is defined in terms of the fan $\Sigma_{j-1}$ which subdivides the
cone $\r'$, with respect to the lattice $N'_{j-1} = N_{j-1} \times
\Z$, into cones $\bar{\s}^+_j$, $\bar{\s}^+_j$ and $\bar{\r}_j$
for $j=1, \dots, g$. We denote by $\Sigma_{g+1}$ the fan of faces
of the cone $\r'$, with lattice $N'_g$.

\begin{lemma}\label{lemmaExc-j}       We have the following equalities for  $1 \leq j \leq g$,
\begin{equation}\label{Exc-j}
(f_j \circ \Pi_{1} \circ \dots \circ \Pi_{j})_{|Z_{\bar{\r}_j}} =
x^{n_j \gamma_j} y_j \quad \mbox{ and } \quad (f \circ \Pi_{1}
\circ \dots \circ \Pi_{j})_{|Z_{\bar{\r}_j}} =      x^{e_{j-1}
\gamma_j}  f^{(j)}.
\end{equation}
\end{lemma}
\textit{Proof.}
This is consecuence of the decompositions
(\ref{factor}), (\ref{fac-j})  and (\ref{factor1}) at levels $2,\dots,j$, together
with Lemma \ref{P4eqexp}.
\hfill $\square$

\begin{remark}
The strict transform  $(S^{(g)}, o_g)$ of $(S,0)$ by the
composition $\Pi_1 \circ \cdots \circ \Pi_g$ is isomorphic to the
germ of toric variety $(Z_{\r, N}, 0)$, with $N= N_g$. This germ
is the normalization of $(S,0)$ (see \cite{GP-Fourier}).
\end{remark}

\begin{remark}     \label{iso-semi-i}
We define the lattice homomorphism $\phi_j \colon M_{j-1}' \to
M_j$, $(\l, a) \mapsto \l + a (\l_j - \l_{j-1})$. By Lemma
\ref{iso-semi} and induction we get that the restriction of the
dual homomorphism $\phi_j^* : N_j \to N'_{j-1}$ defines an
isomorphism of semigroups
\begin{equation} \label{ide}
\phi^*_{j| \r \cap N_j} \colon  \r \cap N_j \to \bar{\r}_j \cap
N'_{j-1},
\end{equation}
for $j=1, \dots, g$. This isomorphism induces a glueing of the
cones $\bar{\r}_j \in \Sigma_j$ and $\bar{\r}_{j-1} \in \Sigma_{j-1}$,
which preserves the lattice points. As a consequence of this
isomorphism we deduce that the germ $(Z_j, o_j)$ is analytically
isomorphic to the germ of normal toric variety $(Z_{\r', N_j'},
0)$. 
\end{remark}

\begin{remark}  \label{iso-complex}
 Glueing the fans
$\Sigma_j$, $j =1, \dots, g+1$ along the identifications
(\ref{ide}), defines a \textit{conic polyhedral complex
 with integral structure}  $\bar{\Theta}$
associated to the q.o.~hypersurface $(S,0)$ (see
\cite{GP-Fourier}).
\end{remark}

The theory of \textit{toroidal embeddings},
introduced by Kempf \textit{et al.} in \cite{Kempf} Chapter II,
is a generalization of the theory of toric varieties. A toroidal
embedding is a pair $(V, U)$,  where $V$ is a normal variety and
$U$ is an open smooth subvariety. In addition, the local algebra
of $V$ at any point $v \in V$ is formally analytically isomorphic
to the local algebra of an affine toric variety $Z$ with torus
$T$, in such a way that the ideal of $V \setminus U$ corresponds
to the ideal of $Z \setminus T$. The toroidal embedding is called
\textit{without self-intersection} if the irreducible components
of $V - U$ are normal. In that case the toroidal embedding can be
attached  to a \textit{conic polyhedral complex with integral
structure}, a collection $\Xi$ of abstract rational polyhedral
cones with certain glueing conditions. Their cells are in one to
one correspondence with the strata of the stratification of $V$
associated to $(V, U)$. A finite rational polyhedral subdivision
$\Xi'$  of $\Xi$ induces a modification $V' \rightarrow V$. This modification 
defines an isomorphism $U' \to U$ over $U$, for some open subset $U' \subset V'$,
in such a way that  the pair  $(V', U')$ is also a toroidal embedding. In particular, if  $\Xi'$ is a \textit{regular
subdivision} of $\Xi$
then the map $V' \to V$ is a resolution of singularities of $V$.

In the case studied here, the pair $(Z_g, U)$ consisting of the
normal variety $Z_g$ and the complement $U \subset Z_g$ of the
strict transform of the hypersurface $\{ x_1 \cdots x_d f_0 \cdots
f_g = 0 \} \subset \C^{d+1}$ by the modification $\Pi_{1} \circ
\cdots \circ \Pi_{g}$, defines a toroidal embedding without
self-intersection. Its conic polyhedral complex is $\bar{\Theta}$
and  this complex is isomorphic, by an isomorphism which preserves
the integral structure, to the fan $\Theta$ of Definition
\ref{complex} (see \cite{GP-Fourier}). Figure \ref{complex-figure}
represents the projectivization of the conic polyhedral complex
$\bar{\Theta}$ associated to a q.o.~surface with two
characteristic monomials. The letters $f_j$ (resp. $x_i$) label the vertices associated to the strict transforms of
$f_i =0$ (resp. $x_i=0$) by the modification $\Pi_{1} \circ
\Pi_{2}$.

\setlength{\unitlength}{0.5 mm}
\begin{figure}
\begin{center}
\begin{picture}(-60,40)(-6,-1)
\linethickness{0.4mm}

\drawline(-41,20)(-16,-20)
\drawline(-41,20)(-66,-20)
\drawline(-66,-20)(-16,-20)
\drawline(-55,-3)(-26,-3)
\drawline(-5,20)(-26,-3)
\drawline(-55,-3)(-5,20)

\drawline(-28,9.5)(-17,7)
\drawline(-28,9.5)(12,20)
\drawline(12,20)(-17,7)

\jput(-72,-23){{${x}_1$}}
\jput(-15,-23){{${x}_2$}}
\jput(-43,22){{${f}_0$}}
\jput(-5,22){{${f}_1$}}
\jput(12,22){{$f_2$}}

\end{picture}
\end{center}

\

\

\

    \caption[]{}\label{complex-figure}
\end{figure}

\begin{theorem}\label{P4thmres}{\rm (\cite{GP-Fourier}, Theorem 1)}
The proper morphism $\Pi : = \Pi_{1} \circ \cdots \circ \Pi_{g}$
is an embedded normalization of $(S,0) \subset (\C^{d+1}, 0)$. An
embedded resolution of $(S,0)$, 
\[
\Pi' \colon Z' \to (\C^{d+1}, 0) 
\]
is obtained by composition of $\Pi$
with the modification associated to a regular subdivision $\bar{\Theta}'$ of the
conic polyhedral complex $\bar{\Theta}$.
\end{theorem}

 \section{Contact of the arcs with the sequence of semi-roots}   \label{contact}

We consider the irreducible q.o.~polynomial $f \in \C\{x_1, \dots, x_d \} [y]$
together with its semi-roots $f_0, \dots, f_g =f$.

\begin{notation} \label{not
contact}
 We set
\[
 F := (x_1, \dots, x_d, f_0, \dots, f_g)
\]
We denote by $(V,0) \subset (\C^{d+1}, 0)$ the hypersurface germ  defined by
$x_1 \cdots x_d f_0 \cdots f_g =0$.
We denote by $\mathcal{H}$ the set of arcs
$
 \mathcal{H}:=   \Lcal ( \C^{d+1})_0 \setminus \Lcal (V)_0$.
If $\varphi \in \mathcal{H}$ then   the vector
\[
 {\ord} (F \circ \varphi) := ( \ord (x_1 \circ \varphi), \dots,  \ord (x_d \circ \varphi),
 \ord (f_0 \circ \varphi),  \dots, \ord (f_g \circ \varphi)),
\]
belongs to $\Z^{d+g+1}_{>0}$.
\end{notation}

In this section we deal with the problem of determining which are the vectors $k \in \Z^{d+g+1}$ such that
$k = {\ord} ( F \circ \varphi) $, for some $\varphi \in \mathcal{H}$.
For this purpose it is useful to pass by
a \textit{coordinate free} approach to deal
with arcs in toric varieties.

\subsection*{Description of arcs in toric varieties}

Let $N$ be a  lattice of rank $d$ and $\theta \subset N_\R$ be a $d$-dimensional rational strictly convex cone.
We denote by $O$ the origin of the toric variety $Z_\theta$.

An arc $\psi \in \Lcal (T_N)$ defines a group homomorphism
\begin{equation} \label{arc-torus}
 M \to \C[[t]]^*, \quad m \mapsto x^m \circ \psi,
\end{equation}
where $\C[[t]]^*$ denotes the group of units of $\C[[t]]$.
Conversely, any group homomorphism (\ref{arc-torus})
defines an arc in the torus $T_N$.

   The arc $\varphi \in \Lcal (Z_\theta) $ has
\textit{generic point} in the torus $T_N$ if and only if $x^m
\circ \varphi \ne 0$,  for any $m \in M =N^*$.
\begin{definition}
\label{ac-gen}
The functions
\[
 \orden (\varphi) : M \to \Z, \quad m \mapsto \ord ( x^m \circ \varphi) , \mbox{ and }
\ac (\varphi) : M \to \C^*, \quad m \mapsto \ac( x^m \circ \varphi)
\]
are  group homomorphisms.
\end{definition}

\begin{remark}   \label{remi}
 If $\varphi \in \Lcal (Z_\theta)$ has generic point in the torus
then we have:
\begin{enumerate}
\item[(i)] The vector $\orden (\varphi) $ belongs to $N$ and $\ac(\varphi)$ defines a closed point in the torus $T_N$.

 \item[(ii)]  The condition  $\varphi (0) =O$ holds
if and only if $\ord (x^m \circ \varphi) >0$, $\forall \, m \in \theta^\vee \cap M \setminus \{0 \}$, that is,
if and only if $\orden (\varphi) \in \stackrel{\circ}{\theta} \cap N$.

\item[(iii)] The map
\[
\mathrm{u} (\varphi) \colon M \to \C[[t]]^*, \quad m \mapsto t^{- \ord(\varphi) (m)} (x^m \circ \varphi)
\]
 is a group homomorphism, i.e., an arc in the torus $T_N$.

\item[(iv)] The homomorphisms
 $\ord(\varphi)$ and $\mathrm{u} (\varphi)$ determine the arc $\varphi$.
\end{enumerate}
\end{remark}

\begin{lemma}   \label{ord-ac1}
Let   $\Sigma$ be a fan subdividing the cone $\theta$. We denote
by     $\p_\Sigma : Z_\Sigma \to Z_\theta$ the toric modification
defined by $\Sigma$. If $\varphi  \in \Lcal (Z_\theta)_O $ is an
arc with generic point in the torus, there is a unique arc
$\varphi' \in \Lcal(Z_\Sigma)$ such that $\p_\Sigma \circ \varphi'
= \varphi$. The origin   $\varphi'(0)$  of the lifted arc
$\varphi'$ belongs to the orbit $\orb (\t)$ of the unique cone $\t
\in \Sigma$ such that $\orden (\varphi) \in \stackrel{\circ}{\t}$.
Moreover we have that
\begin{equation} \label{ord-ac}
  \orden (\varphi) = \orden (\varphi') \mbox{ and }  \ac (\varphi) = \ac (\varphi').
 \end{equation}
\end{lemma}
\textit{Proof.}
         By the valuative criterion of properness there is a unique arc $\varphi' \in \Lcal(Z_\Sigma)$
lifting $\varphi$. The map $\p_\Sigma$ is a proper modification which is the identity on the torus,
hence equality (\ref{ord-ac})   holds.
 The statement on the origin of the arc follows analogously as
Remark \ref{remi} (ii).
\hfill $\square$

\subsection*{The order of contact of arcs in $\mathcal{H}$ with the sequence $F$}

Now we apply the above notions  to the case of an arc $\varphi \in \mathcal{H}$.

As in Section \ref{resolution} we assume that $y = f_0 $.
If $\p$ denotes the projection (\ref{projection})
we get that $\p \circ \varphi$ is an arc in $\Lcal (\C^d)_0$ with generic point in the torus
$T_{N_0}$ of $\C^d = Z_{\r, N_0}$.
Setting  $\nu := \orden (\p \circ \varphi) \in N_0$ we deduce that
\begin{equation} \label{split}
\orden (\varphi)  = (\nu, \ord (y \circ  \varphi) ) \in N_0 \times
\Z = N_0' \mbox{ and } \ac (\varphi ) \colon M_0' \to \C^*,
\end{equation}
correspond to the group homomorphisms introduced in Definition
\ref{ac-gen}.  The coordinates of $\nu$ with respect to the basis
$\epsilon_1, \dots, \epsilon_d$ of $N_0$
 are
$(\ord (x_1 \circ \varphi), \dots, \ord (x_d \circ \varphi))$ (see
Notation \ref{basis}).

The toric map $\Pi_1$ of $\C^{d+1}$ induced by $\Sigma_1$ (see
Section \ref{resolution}) is a  proper modification, which is the
identity over the torus $T_{N_0'}$. Since $\varphi$ has generic
point in  $T_{N_0'}$, by the valuative criterion of properness
there is a unique lifting $\varphi^{(1)} \in \Lcal (Z_1)$ of
$\varphi$, i.e., there exists a unique $\varphi^{(1)} \in \Lcal
(Z_1)$ such that   $\varphi = \Pi_1 \circ \varphi^{(1)}$.
   The vector $\orden ( \varphi ) $ belongs to $\bar{\s}^+$ (resp. $\bar{\s}^-$) if and only if
\[
 \langle \nu, \l_1 \rangle \leq \mbox{ (resp. } \geq \mbox{  ) } \ord (y \circ \varphi).
\]

Proposition \ref{lift1}
characterizes the case when $\varphi^{(1)} (0) = {o_1}$.
\begin{proposition}    \label{lift1}
The following conditions are equivalent:
\begin{enumerate}
\item[(i)]  $\varphi^{(1)} (0) = o_1$,

\item[(ii)] $\ord (y_1 \circ \varphi^{(1)}) > 0$,

\item[(iii)] $\orden (\varphi) \in \mathrm{int}(\bar{\r}_1)$ and $(f_1)_{| \orden \varphi } ( \ac (\varphi) ) =0 $.
\end{enumerate}

\end{proposition}
\textit{Proof.}
We show first the equivalence (i) $\Leftrightarrow$ (ii). If $\varphi^{(1)} (0) = o_1$ then we get
$\ord (y_1 \circ \varphi^{(1)} )
> 0$,  since $y_1$ belongs to a
system of generators of the maximal ideal defining the point $o_1$
at $Z_1$. The other generators are monomials of the form $x^m$ for
$m \in \r^\vee \cap M_1$ and $m\ne 0$. For those vectors $m$ we
have that $\ord (x^m \circ \varphi ) = \ord ( x^m \circ
\varphi^{(1)})
>0 $ by (\ref{ord-ac}) and Remark \ref{remi}.

Now we deal with the equivalence (ii) $\Leftrightarrow$ (iii). The
condition $\orden (\varphi) \in \mathrm{int}(\bar{\r}_1)$ is
equivalent to $\varphi^{(1)} (0) \in \orb ( \bar{\r}_1)$. If $\orden
(\varphi) \in \mathrm{int}(\bar{\r}_1)$ then we get that
$(f_{1})_{| \orden(\varphi)} = y^{n_1} - c_{\l_1}^{n_1} x^{\l_1}$.
Since  $y_1 = (x^{-n_1 \l_1} f_1\circ \Pi_1 )_{| Z_{\bar{\r}_1}} $
we deduce  that the constant term of the series $y_1 \circ
\varphi^{(1)} $ is equal to
\[
\ac(x^{-n_1 \l_1} \circ \varphi ) \left( (f_{1})_{| \orden(\varphi)} (\ac
(\varphi) \right).
\]
This ends the proof of the equivalence between (ii) and (iii).
 \hfill $\square$

If  $\varphi \in \mathcal{H}$ then the series  $y_1 \circ
\varphi^{(1)}$ is nonzero. Hence $\varphi^{(1)}$ can be seen as
an arc with generic point in the torus $T_{N_1'}$, associated to
the toric structure of the germ $(Z_1, o_1)$ (see Lemma
\ref{iso-chart}). By the valuative criterion of properness there
is a unique lifting $\varphi^{(2)} \in \Lcal (Z_2)$ of
$\varphi^{(1)}$ such that $\Pi_2 \circ \varphi^{(2)} =
\varphi^{(1)}$. Notice also that $y_2 \circ \varphi^{(2)} \ne 0$
by (\ref{Exc-j}) and the definition of $\mathcal{H}$. By iterating
this argument, using that $\varphi \in \mathcal{H}$, there exists
a unique lifting of $\varphi^{(i)} \in \Lcal (Z_j)$ such that
$\Pi_1 \circ \cdots \circ \Pi_j \circ \varphi^{(i)} = \varphi$,
for $i=1, \dots,g$. 

We get different notions of
orders of the lifted arcs
 which depend on the various toric structures
appearing at toric resolution.
\begin{notation}         \label{split-i}
We denote by
$
\orden^{(i)} (\varphi^{(i)} ) \colon M_i' \to \Z $ and 
$\ac^{(i)} (\varphi^{(i)} ) \colon M_i' \to \C^*$,
the group homomorphisms introduced in Definition \ref{ac-gen},
with respect to the toric structure of the germ $(Z_i, o_i)$ (see
Remark \ref{iso-semi-i}). By convenience, we denote also $\varphi$
by $\varphi^{(0)}$ and $\orden(\varphi)$, $\ac (\varphi)$ in
(\ref{split}) by $\orden^{(0)} (\varphi^{(0)})$,
$\ac^{(0)} (\varphi^{(0)})$ respectively. Sometimes $y_0$ also denotes
$f_0$.
\end{notation}

\begin{proposition}    \label{liftj}
Suppose that $\varphi^{(i)}(0) = o_i$ for $i=1, \dots, j-1$.
The following conditions are equivalent:
\begin{enumerate}
\item[(i)]  $\varphi^{(j)} (0) = o_j$,

\item[(ii)] $\ord (y_j \circ \varphi^{(j)}) > 0$,

\item[(iii)] $\orden^{(j-1)} (\varphi^{(j-1)}) \in
\stackrel{\circ}{\bar{\r}_j}$ and $(f_j^{(j-1)})_{| \orden^{(j-1)}
(\varphi^{(j-1)}) } ( \ac (\varphi^{(j-1)}) ) =0 $.
\end{enumerate}
\end{proposition}
\textit{Proof.} The proof follows by induction using the same
arguments as in Proposition \ref{lift1}.
\hfill $\square$

Notice that if $ \ord ( y_j \circ \varphi^{(j)} ) > 0$ then the inequalities
                 $ \ord (y_i\circ \varphi^{(i)})  > 0$ hold for $0 \leq i < j$.

We introduce the  \textit{depth} of an arc $\varphi  \in
\mathcal{H}$, with respect to the modification $\Pi$ of Theorem
\ref{P4thmres}.
\begin{definition}    \label{depth}
If   $\ord (y_g \circ \varphi^{(g)}) > 0 $
we define $\mathrm{depth} (\varphi) := g+1$. Otherwise,  there exists a unique $1 \leq j \leq g$ such that
$\ord (y_j \circ \varphi^{(j)}) = 0$ and  $ \ord (y_{j-1} \circ \varphi^{(j-1)}) > 0$,
and then we set $\mathrm{depth}(\varphi) = j$.
\end{definition}

If     $\varphi \in \mathcal{H}$ is an arc of depth $j$ then
$\varphi^{(i)} \in \Lcal(Z_i)$ is an arc with origin at $o_i$, for
$i=0, \dots, j-1$.

\begin{proposition}       \label{orden-j}
 If $\varphi \in \mathcal{H}$ is an arc of depth $j>1$ and if  $\nu = \orden (\p \circ \varphi)$ then we have
\[
  \ord  (y_0 \circ  \varphi )  =  \langle \nu, \l_1 \rangle    \quad \mbox{ and }    \quad
  \ord (y_i \circ  \varphi^{(i)})   =   \langle \nu, \l_{i+1} - \l_i \rangle, \quad \mbox{ for } i = 2, \dots, j-1.
\]
\end{proposition}
\textit{Proof.}
This is consequence of Proposition \ref{liftj} and the definitions.
\hfill $\square$

\begin{lemma}       \label{split-i2}
If $\varphi \in \mathcal{H}$ is an arc of depth $j \geq 1$ with
$\nu = \orden (\p \circ \varphi)$, then, we have that
\[
 \orden^{(i)}  (\varphi^{(i)}) = (\nu, \ord  ( y_i \circ \varphi^{(i)} ) ) \mbox{ for }  i=0, \dots, j-1.
\]
\end{lemma}
\textit{Proof} The assertion follows by Proposition \ref{orden-j}
and induction. We get $\nu \in N_i$ using that
 $N_{g} \subset \dots \subset N_0$. \hfill $\square$

\begin{theorem}   \label{trop}
 We have the following equalities:
\begin{equation}              \label{prof}
\{ \ord (F \circ \varphi) \mid \varphi \in \mathcal{H} \} = \bigsqcup_{j=0}^{g+1} \bigsqcup_{\theta \in \Theta_j}
\stackrel{\circ}{\theta} \cap \Z^{d+g+1}
\end{equation}
and  if $1 \leq j \leq g+1$,
\begin{equation} \label{prof-j}
\{ \ord (F \circ \varphi) \mid \varphi \in \mathcal{H} \hbox{ {\rm and}
} \mathrm{depth} (\varphi) = j \} = \bigsqcup_{\theta \in
\Theta_j} \stackrel{\circ}{\theta} \cap \Z^{d+g+1}.
\end{equation}
\end{theorem}
\textit{Proof.} It is enough to prove the equality (\ref{prof-j}).
We assume first that $j < g+1$.
If $\varphi \in \mathcal{H}$ and $\mathrm{depth} (\varphi) = j  < g+1$
then there is a unique cone $\bar{\theta} \in \{ \bar{\r}_j, \bar{\s}_j^+, \bar{\s}_j^- \} \subset \Sigma_j$ such that
            \[
             \orden ^{(j-1)} (\varphi^{(j-1)})  \in \,  \mathrm{int}(\bar{\theta}) \cap N_{j-1}'.
            \]
We prove first that $\ord (F \circ \varphi)$ belongs to the
the cone $\theta \in \{
     {\r}_j, {\s}_j^+, {\s}_j^- \} = \Theta_j$.
Suppose that  $\bar{\theta} = \bar{\r}_j$. By Proposition \ref{orden-j}
we have that
 $\orden ^{(j-1)} ( \varphi^{(j-1)} ) = (\nu, r)$ with $r= \langle \nu, \l_j - \l_{j-1} \rangle$.
By Formula (\ref{Exc-j}) we get,
\[
  f_i \circ \varphi = (f_i \circ \Pi_1 \circ \cdots \circ \Pi_i)_{| Z_{\bar{\r}_i}}
   \circ \varphi^{(i)} = (x^{n_i  \g_i} \circ \varphi )
\cdot (y_i \circ \varphi^{(i)}),
\]
for $ i = 1, \dots, j-1$, hence, taking orders we obtain:
\begin{equation} \label{cle}
 \ord (f_i \circ \varphi)  = \langle \nu, n_i \g_i \rangle + \langle \nu, \l_{i+1} - \l_i \rangle =
\langle \nu, \g_{i+1}  \rangle.
\end{equation}
For $i=j$ we get
\begin{equation}\label{cle-2}
 \ord (f_j \circ \varphi)  = \langle \nu, n_j \g_j \rangle,
\end{equation}
since $\ord (y_j \circ \varphi^{(j)}) = 0 $ by hypothesis. For $g
\geq i > j$ we deduce the equality $( f_i \circ \Pi_1 \circ \cdots
\circ \Pi_j ) _{ |  Z_{\bar{\r}_i}}= x^{n_j \cdots n_i \g_j} \cdot
f_i^{(j)}$  from Proposition \ref{orden-j}. By construction, if
$\ord (y_j \circ \varphi^{(j)}) =0$ then it follows that $\ord (f_i
^{(j)} \circ \varphi^{(j)})  = 0$ and then
\[
 \ord (f_i  \circ \varphi) =    n_i \cdots n_{j+1} \langle \nu, n_j \g_j \rangle.
\]
We have shown that $\ord (F \circ \varphi) $ is an integral point
in the relative interior of the cone $\r_j \in \Theta_j$.

Conversely, we prove that any point  $k= (k_1, \dots, k_{d+g+1})
\in \stackrel{\circ}{\r}_j \cap \Z^{d+g+1}$ is of the form $k =
\ord (F \circ \varphi)$, for some arc $\varphi \in \mathcal{H}$.
By definition of $\r_j$ if $\nu := (k_1, \dots, k_d)  \in
\stackrel{\circ}{\r} \cap N_0$ we have the  equalities
\[
 k_{d+i} = \langle \nu, \g_i \rangle, \mbox{ for } i=1, \dots j,
\]
and
\[
 k_{d +j + l} = n_{j+1} \cdots n_{j + l-1}      \langle \nu, n_j \g_j \rangle, \mbox{ for } l=1, \dots g-j+1  .
\]
 Since the vectors $e_1, \dots, e_d, \g_1, \dots, \g_j$ span the
 lattice $M_j$,
 the formulas
above show that $\nu$ belongs to $N_j$.
We get also that  $ k_{d+j} - n_j k_{d+j -1} = \langle \nu, \l_j -\l_{j-1} \rangle >0$
and  $(\nu, k_{d+j} - n_j k_{d+j -1}  )   \in \bar{\r}_j \cap N_{j-1}'$.
Then, any arc $\psi \in \Lcal(Z_{j-1})_{o_{j-1}}$
with generic point in the torus $T_{N_{j-1}'}$ and such that
\[
           \orden^{(j-1)} (\psi) = (\nu, k_{d+j} - n_j k_{d+j -1}  )
\]
verifies that $\varphi:= \Pi_1 \circ \cdots \circ \Pi_{j-1} \circ
\psi $ belongs to $\mathcal{H}$, $\varphi^{(j-1)} = \psi$ and
$\ord ( F \circ \varphi) = k$.

Now we consider the case when $\orden ^{(j-1)} (\varphi^{(j-1)}) =
(\nu, r) $ belongs to the interior of ${\bar{\s}}_j^{\pm}$.
By the previous arguments we get 
\begin{equation}  \label{uso-1}
 \ord (f_i \circ \varphi ) = \langle \nu, \gamma_{i+1} \rangle, \mbox{ for } i = 0, \dots, j - 2, 
 \mbox{ and }
  \ord (f_{j-1} \circ \varphi ) = \langle \nu, n_{j-1} \gamma_{j-1} \rangle  + r, 
\end{equation}
where $r = \ord (y_{j-1} \circ \varphi^{(j-1)})$.
We deduce the formula 
\begin{equation}\label{uso-2}
\left\{
\begin{array}{ccccclcl}
  \langle \nu , n_{j-1} \g_{j-1} \rangle   &  < &  \ord (f_{j-1} \circ \varphi )  &  <  &      \langle \nu , \g_{j} \rangle  & \mbox{ if }
&   (\nu, r) \in \bar{\s}_j^-  ,
\\
  \langle \nu , \g_{j} \rangle  & < &   \ord (f_{j-1} \circ \varphi )  &  &          & \mbox{ if }
&   (\nu, r) \in \bar{\s}_j^+ .
\end{array}   \right.
\end{equation}
Notice also that 
\begin{equation} \label{uso-3}
 \ord (f_{j + l -1} \circ \varphi) =  n_j \cdots n_{j+l -1}  \langle \nu , n_{j-1} \g_{j-1} \rangle   + 
 \ord (f_{j + l -1} ^{(j-1)} \circ  \varphi^{(j-1)})
\end{equation}
where
\begin{equation}   \label{uso-4}
 \ord (f_{j + l -1} ^{(j-1)} \circ  \varphi^{(j-1)}) = 
 \left\{
\begin{array}{lcl}
n_j \cdots n_{j+l -1} \, r  & \mbox{ if }
&   (\nu, r) \in \bar{\s}_j^-  ,
\\
 n_j \cdots n_{j+l -1}   \, \langle \nu, \l_j - \l_{j-1} \rangle   & \mbox{ if }
&   (\nu, r) \in \bar{\s}_j^+ ,
\end{array}   \right.
\end{equation}
for  $1 \leq l \leq g-j+1$. 
We obtain from  the definitions and 
the formulas (\ref{uso-1}), (\ref{uso-2}), (\ref{uso-3}) and (\ref{uso-4})   that 
$\ord (F \circ \varphi)$ belongs to 
    $\mathrm{int} ({\s}_j^-)
\cap \Z^{d+g+1}$  (resp.  $\mathrm{int} ({\s}_j^+)
\cap \Z^{d+g+1}$)  if   $(\nu, r) \in \bar{\s}_j^-$  (resp.  if $(\nu, r) \in \bar{\s}_j^+$ ).  
We prove similarly
that all the points in the sets $ \mathrm{int} ( {\s}_j^{\pm} ) \cap
\Z^{d+g+1}$ are of this form.

The proof in the case of $j = g+1$ and $\bar{\theta}= {\s}_{g+1}$ is analogous.
\hfill $\square$

We end this section with some auxiliary results.

\begin{lemma}   \label{pos-trop}
The map 
${s}_j \colon \r \cap N_j \mapsto \Z^{d+ g +1}$ given by 
\[ 
\nu  \mapsto \nu + \sum_{i=1}^{j} \langle \nu, \g_i \rangle
\epsilon_{d+i}  + \sum_{i=j + 1}^{g+1} n_{j} \cdots n_{i - 1} \langle \nu, \g_{j} \rangle \epsilon_{d+i}, 
\]
defines an isomorphism of semigroups $\r \cap N_j \cong  \r_j \cap \Z^{d+g+1}$. 
\end{lemma}
\textit{Proof}. This is consequence of the proof of Theorem \ref{trop}. \hfill $\square$

\begin{lemma} \label{minus}
We fix an integer $1 \leq j \leq g$.
The linear map $\Upsilon_{j} : (N_{j-1}')_\R \to \R^{d+g+1}$
given by
\[
(\nu, r) \mapsto \nu + \sum_{i=1}^{j-1} \langle \nu, \g_i \rangle
\epsilon_{d+i} + ( r + n_{j-1} \langle \nu, \g_{j-1} \rangle)
\epsilon_{d+j}  + \sum_{i=j + 1}^{g+1} n_{j} \cdots n_{i - 1} ( r +
n_{j-1} \langle \nu, \g_{j-1} \rangle) \epsilon_{d+i},
\]
verifies that $\Upsilon_j (\bar{\s}_j^- \cap N_{j-1}') = \s_j^-
\cap \Z^{d+g+1}$ and $ \Upsilon_j (0,1) = \epsilon_{d+j} + \sum_{i
= j+1}^{g+1} n_{j} \cdots n_{i-1} \epsilon_{d + i }$. The
restriction of $\Upsilon_j$  to $N_{j-1}'$ is an isomorphism onto
the lattice spanned by $\bar{\s}_j^- \cap \Z^{d+g+1}$.
We set  $\hat{\theta} := \Upsilon_j ( \theta) $ for $\theta
\subset \r'$ a cone in  $N_{j-1}'$. The generating function of the
set $\mathrm{int}({\hat{\theta}}) \cap \Z^{d+ g+1}$ is of the form
$G_{\hat{\theta}} = P_{\hat{\theta}} \prod_v (1 - x^v)$ (see Lemma
\ref{gf}). We denote by $\tilde{P}_{\hat{\theta}}$ the image of
the polynomial $P_{\hat{\theta}}$ by the monomial map
(\ref{polymap}) and by $S_{\hat{\theta}}$ the rational function
defined by replacing $\theta$ by $\hat{\theta}$  in Formula
(\ref{stheta}). The following equality holds:
\begin{equation} \label{minus-rel}
S_{\s_j^-}  = S_{\hat{\r}'} - S_{\hat{\bar{\r}}_j} -
S_{\hat{\bar{\s}}_{j}^+}
\end{equation}
\end{lemma}
\textit{Proof.} The assertion about the map $\Upsilon_j$ is
consequence of the proof of Theorem \ref{trop}. The formula
(\ref{minus-rel}) holds because $S_{\s_j^-} =
S_{\hat{\bar{\s}}_{j}^-}$ and $G_{\hat{\r}'} = G_{\hat{\bar{\r}}_j}
+ G_{\hat{\bar{\s}}_{j}^+} + G_{\hat{\bar{\s}}_{j}^-}$ by
additivity of generating functions of integral points with respect
to the partition $ \mathrm{int}(\hat{\r}') = \mathrm{int} (
\hat{\bar{\r}}_j ) \sqcup \mathrm{int} (\hat{\bar{\s}}_{j}^+)
\sqcup \mathrm{int} ( \hat{\bar{\s}}_{j}^- )$.
\hfill $\square$

\section{Computation of the motivic measure of $H_{k, 1}$}\label{medida}

In this  section we compute the motivic measure of the set
\[
 H_{k,1} := \{ \varphi \in \mathcal{H}  \mid \ord (F \circ \varphi) = k, \ac (f \circ \varphi) = 1 \}.
\]
by using the change of variables formula for motivic integrals
(see Theorem \ref{P4changevariables}). We assume in this section
that  $k$ is fixed and that $H_k$ is non-empty, hence by Theorem
\ref{trop},  there exists a unique $j \in \{1, \dots, g+1 \}$ and
a unique cone $\theta \in \Theta_j$ such that  vector $k$ belongs
to $\stackrel{\circ}{\theta} \cap \Z^{d+g+1}$.

We fix a regular subdivision $\bar{\Theta}'$ of the conic
polyhedral complex $\bar{\Theta}$ associated to the toric
resolution of singularities, that is, a sequence $\Sigma_i'$ of
regular subdivisions of the fans $\Sigma_i$, which are compatible
with the glueing of faces mentioned in Remark \ref{iso-semi-i}. By
Proposition \ref{orden-j} the vector $\orden^{(i-1)}
(\varphi^{(i-1)}) $ belongs to $\bar{\r}_i$, for $i=0,\dots, j-1$
(see Definition \ref{split-i}). Let $\t_i \in \Sigma_i'$ be a
$d+1$-dimensional regular cone containing the vector $\orden^{(i-1)}
(\varphi^{(i-1)}) $. We can always chose a regular subdivision
$\bar{\Theta}'$ such that
\begin{equation} \label{cone-cond}
 \dim \t_i \cap \bar{\r}_i = d-1,       \mbox{ for } i = 1, \dots, j-1,    \quad \mbox{  
(see Figure \ref{cone-cond2})}.
\end{equation}

The cone $\t_i$ is generated by a basis  $a_1^{(i)}, \dots,
a_{d+1}^{(i)}$ of the lattice $N_{i-1}'$ and by the assumption
(\ref{cone-cond})  we can guarantee that $a_1^{(i)}, \dots,
a_{d}^{(i)} \in \bar{\r}_i$. By Lemma \ref{iso-semi} and Remark
\ref{iso-semi-i} the
 pair consisting of the cone
$\bar{\r}_i$ and the lattice generated by  $\bar{\r}_i \cap N_{i-1}'$
is isomorphic to  $(\r, N_i)$.
We identify these pairs under this isomorphism and we denote by
  $\mathcal{B}_i$
the basis of
$N_i'$
given by $\mathcal{B}_i := \{ (a_1^{(i)}, 0), \dots, (a_{d}^{(i)}, 0), (0,1) \}$ for $i=1, \dots, j-1$.
We set $\mathcal{B}_0 := \{ \epsilon_1, \dots, \epsilon_{d+1} \}$ the canonical basis of $\Z^{d+1}$.

We denote by $x_1^{(i)}, \dots, x_{d+1}^{(i)}$ the coordinates of $Z_{\t_i, N_{i-1}'} \cong \C^{d+1}$
associated to the dual basis of   $a_1^{(i)}, \dots, a_{d+1}^{(i)}$.
The toric morphism $\p_{\t_i}$ is written in this coordinates as:
\begin{equation}
\label{local-coor}
\begin{array}{lcllll}
x_1^{(i-1)} & = &
(x_1^{(i)}) ^ {a^{(i)}_{1,1}} &
(x_2^{(i)})^{a^{(i)}_{2,1}}    &
\cdots                         &
(x_{d+1} ^{(i)})^{a^{(i)}_{d+1,1}}
\\
\dots & \dots &    \dots & \dots &  \dots &    \dots
\\
x_{d}^{(i-1)} & = &   (x_1^{(i)})^{a^{(i)}_{1,d}}
&
(x_2^{(i)})^{a^{(i)}_{2,d}}
&
\cdots
&
(x_{d+1} ^{(i)})^{a^{(i)}_{d+1,d}}
\\
y_{i-1}        & = &
 (x_1^{(i)})^{a^{(i)}_{1,d+1}}
&
(x_2^{(i)})^{a^{(i)}_{2,d+1}}
&
\cdots
&
(x_{d+1} ^{(i)})^{a^{(i)}_{d+1,d+1}},
\\
\end{array}
\end{equation}
where
$(a^{(i)}_{r,1}, \dots, a^{(i)}_{r, d+1})$ denote the coordinates of the vector $a^{(i)}_{r}$, for $r=1, \dots, d+1$,
with respect to the basis     $\mathcal{B}_{i-1}$ of $N_{i-1}'$.

We denote by $\1^{(i)}$ the vector in $M_i'$ corresponding to the
monomial $x_1^{(i)} \cdots x_d^{(i)} y_{i}$, for $i=0,    \dots,
j-1$.

\begin{remark}\label{obstau}
We abuse of notation and denote by the same letter the lifting
$\varphi^{(i)}$ of $\varphi$ by the maps $\p_{\t_1} \circ \cdots
\circ  \p_{\t_{i}}$ and $\Pi_1 \circ \cdots \circ \Pi_{i}$. This
causes no confusion since our analysis depends only on the group
homomorphisms $\orden^{(i)} (\varphi^{(i)})$ and $\ac^{(i)}
(\varphi^{(i)})$ of Notation \ref{split-i}. These homomorphisms
coincide with
 those defined after taking the regular
subdivisions by Lemma \ref{ord-ac1}.
\end{remark}

\setlength{\unitlength}{0.4 mm}
\begin{figure}
\begin{center}
\begin{picture}(-60,40)(-6,-1)
\linethickness{0.4mm}

\put(-35,-10){\vector(1,0){30}}
\put(-35,-10){\vector(0,1){50}}
\put(-35,-10){\vector(-1,-1){20}}

\jput(-45,10){$\sigma_1$}

\drawline(-35,-10)(5,25)
\drawline(-60,10)(5,25)
\drawline(-60,10)(-35,-10)
\jput(5,25){\circle*{2}}
\jput(-60,10){\circle*{2}}

\jput(-13,20.5){\circle*{2}}

\dottedline{2}(-35,-10)(-13,20.5)
\drawline(-35,-10)(-13,30)
\drawline(-35,-10)(-3,23)
\drawline(-35,-10)(-25,18)
\dottedline(-13,30)(-3,23)
\dottedline(-13,30)(-25,18)

\jput(-5,27.5){$\tau$}
\jput(-12,23.5){$\underline{{\bf k}}$}

\end{picture}
\

\

\

\end{center}
\caption[]{}\label{cone-cond2}

\end{figure}

\begin{proposition}     \label{jacobian}
With the previous notations if $\varphi \in H_k$ is of depth $j$ then
\[
 \ord (\mathrm{Jac} (\p_{\t_1} \circ \dots \circ \p_{\t_{j-1}} ) \circ \varphi^{(j-1)} )  =  \xi_j (k)
\  - \langle  \orden^{(j-1)} (\varphi^{(j-1)} ) , \1^{(j-1)} \rangle.
\]
\end{proposition}
\textit{Proof.}
We have that the critical locus of the morphism $  \p_{\t_1} \circ \dots \circ \p_{\t_{j-1}} $ is defined by
\begin{equation} \label{jac-0}
 \mathrm{Jac}    (\p_{\t_1} \circ \dots \circ \p_{\t_{j-1}} ) =
(\mathrm{Jac}(\p_{\t_1}) \circ   \p_{\t_2} \circ \dots \circ \p_{\t_{j-1}} ) \cdots
(\mathrm{Jac} (\p_{\t_{j-1}})).
\end{equation}
From the expression (\ref{local-coor})
we get that
\begin{equation}                                                \label{jac-1}
\mathrm{Jac} (\p_{\t_i}) )=     ((x_1^{(i-1)}   \cdots x_{d}^{(i-1)} y_{i-1}) \circ \p_{\t_i} )
( x_1^{(i)} \cdots  x_{d+1}^{(i)})^{-1}.
\end{equation}
If $\varphi \in H_k$ our hipothesis guarantee that $\ord (
x_{d+1}^{(i)} \circ \varphi^{(i)} ) =0$, since $a_{d +1}^{(i)}
\notin \bar{\r}_i$. Composing with $\varphi^{(i)}$ at both sides
of (\ref{jac-1})
 and taking orders provides the equality:
\begin{equation}   \label{jac-2}
\ord      ( \mathrm{Jac} (\p_{\t_i})  \circ \varphi^{(i)} ) =
\langle \orden^{(i-1)} \varphi^{(i-1) }, \1^{(i-1)} \rangle
-
\langle \orden^{(i)} \varphi^{(i) }, \1^{(i)} \rangle
+
\ord (y_i \circ \varphi^{(i)}).
\end{equation}
By definition of the lifting by the toric morphisms we get that
$
 \p_{\t_{i+1}} \circ \dots \circ \p_{\t_{j-1}} \circ \varphi^{(j-1)} = \varphi^{(i)}$,
   for any $0 \leq i < j-1$.
We apply this observation
when we
take orders after  composing both sides of (\ref{jac-0}) with $\varphi^{(j-1)}$.
By using (\ref{jac-2})  we get that    $\ord (\mathrm{Jac} (\p_{\t_1} \circ \dots \circ \p_{\t_{j-1}} ) \circ \varphi^{(j-1)} )$
is equal to
\begin{equation}
\langle \orden^{(0)} ( \varphi ), \1^{(0)} \rangle - \langle  \orden^{(j-1)} \varphi^{(j-1)} , \1^{(j-1)} \rangle  +
\sum_{i=1}^{j-1}  \ord( y_i \circ \varphi^{(i)}).
\end{equation}

By Proposition \ref{orden-j}  and equality    (\ref{cle}) we
deduce that $ \ord (y_i \circ \varphi^{(i)}) = k_{d+i+1} - n_i
k_{d+i}. $ The assertion follows since $    \langle \orden^{(0)} (
\varphi ), \1^{(0)} \rangle  = k_1 + \cdots + k_{d+1}$. See
Definition \ref{xi}. \hfill $\square$

\begin{remark} \label{curvette}
If $\f$ is a regular function on a smooth algebraic variety $Z$ and $E$ is 
an irreducible divisor on $E$ then the order of vanishing of $\f$ along $E$ coincides with 
the order of $\f \circ \a$, where $\a$ is a \textit{curvette} at $E$, that is, $\a (0) \in E$ and 
$\a$ is a smooth arc transversal to $E$.  
\end{remark}

\begin{remark} \label{vanish}
With the notations  and hypothesis of Proposition \ref{jacobian} assume in addition
that $\orden^{(j-1)} ( \varphi^{(j-1) }) $ is one element of the basis 
$a_1^{(j-1)}, \dots, a_{d+1}^{(j-1)}$. Then we deduce that 
\begin{equation} \label{jac-3}
 \ord (\mathrm{Jac} (\p_{\t_1} \circ \dots \circ \p_{\t_{j-1}} )) \circ \varphi^{(j-1)} )  =  \xi_j (k)  - 1.
\end{equation}
Notice that such a vector, say $\orden^{(j-1)} (\varphi^{(j-1) }) = a_1^{(j-1)}$, 
corresponds to an irreducible divisor
$E$, which is defined on the chart $Z_{\t_{j-1}, N_{j-1}'}$ by $x_1^{(j-1)} = 0$.
By Lemma \ref{ord-ac1} there is a unique lifting $\varphi^{(j)}$ of the arc  $\varphi^{(j-1)} $
such that $ \p_{\t_{j-1}} \circ \varphi^{(j)} = \varphi^{(j-1)} $ and 
$\varphi^{(j)} $ is  a curvette at $E$. 
We abuse of notation and denote also by $E$ the strict transform of the divisor $E$ in $Z'$, 
where we recall that $\Pi'\colon Z' \to \C^{d+1}$ is the toric embedded resolution of $(S,0) \subset (\C^{d+1},0)$ (see Theorem \ref{P4thmres}).  
The 
order of vanishing of $\mathrm{Jac} (\Pi')$ (resp of $\Pi' \circ f$) along the divisor $E$ is equal to 
$\xi_j (k) -1$ (resp $\eta(k)$).
This follows from Formula (\ref{jac-3}) and Theorem \ref{trop}, by 
using that these orders of vanishing
along $E$ coincide with that of $\mathrm{Jac} (\p_{\t_1} \circ \dots \circ \p_{\t_{j-1}} )$
(resp. of $f \circ \p_{\t_1} \circ \dots \circ \p_{\t_{j-1}} $).
\end{remark} 

\begin{remark} \label{vanish-3}
We apply remark \ref{vanish} when $k = \nu_i^{(j)}$ for some $1 \leq i \leq d$ and $1 \leq j \leq g$,  and then $E_i^{(j)}$ denotes the divisor $E$.
See Remark \ref{vanish-2}.  
\end{remark}

Recall the notations introduced in Definitions \ref{integral} and
\ref{csigma}.
\begin{proposition}       \label{mu-k} If $k \in \stackrel{\circ}{\theta} \cap \Z^{d+g+1}$ for some
$\theta \in \Theta_j$, $j =1, \dots g+1$,
we have the following formula for the motivic measure of the set $H_{k,1}$:
\[
     \mu (H_{k,1}) = c_1 (\theta) (\L -1)^{\dim \theta -1} \L^{- \xi_j (k) }.
\]
\end{proposition}
\textit{Proof.}    We deal first with the case when the arcs in $H_k$ are of depth one.
The set $H_{k,1}$ is stable for $l \geq \max \{k_1, \dots, k_{d+1} \}$, thus the
motivic measure of $\mu (H_{k,1})$ is equal to
\begin{equation} \label{me0}
 \mu (H_{k,1}) = [\pi_{l} (H_{k,1} ) ]     \L^{ -l (d+1) }.
\end{equation}
Now we compute the class  $[\pi_{l} (H_{k,1} ) ]  \in \kmon$.
If $\varphi \in H_{k,1}$ we have the expansions
\[
 x_i \circ \varphi  = \ac (x_i \circ \varphi) t^{k_i} + c_{i, 1} t^{k_i + 1} + \cdots + c_{i, l - k_i} t^{l}  
\quad   {\mod t^{l+1}}.
\]
The coefficients $ c_{i, 1}, \dots, c_{i, l - k_i}$, for $i =1,  \dots,d+1$, are free and contribute to the class
$[\pi_{l} (H_{k,1} ) ]$ with the factor $ \L^{l(d+1) - \sum_{i=1}^{d+1} k_i}$.

The angular coefficients  $\ac(x_i \circ \varphi)$ determine the
point $\ac( \varphi) \in T_{N_0}$. Since $\varphi \in H_{k, 1}$ is
of depth one, the equality $\ac (f \circ \varphi) =  f_{| \orden
(\varphi)} (\ac(\varphi))$ holds. It follows  that the angular
coefficients contribute to the class
 $[\pi_{l} (H_{k,1} ) ]$
with the factor $[\{ \ac (\varphi) \in T_{N_{0}'} \mid  f_{| \orden
(\varphi)} (\ac(\varphi))  =1 \} ]$. The result in this case
follows by formula (\ref{forma-f}) and Definition \ref{csigma}.

Now we assume that the arcs in $H_k$ are of depth $1 < j \leq
g+1$. We set
\[
H_{k,1}^{(j-1)} := \{ \varphi^{(j-1)} \mid \varphi \in H_{k, 1} \}.
\]
We apply the change of variables formula for motivic integrals to
the transformation $\Psi : = \p_{\t_1} \circ \cdots \circ
\p_{\t_{j-1}}$ (see Theorem \ref{P4changevariables}).
\[
    \mu (H_{k,1}) = \int_{\varphi \in H_{k,1}} d\mu = \int_{\varphi^{(j-1)} \in  H_{k,1} ^{(j-1)} }
\L^{- \ord(\mathrm{Jac} (\Psi) \circ \varphi^{(j)})  } d\mu .
\]
By Proposition  \ref{jacobian} we get
that
\begin{equation} \label{me1}
 \mu (H_{k,1}) =   \L ^{ - \xi_j (k)
\  + \langle  \orden^{(j-1)} \varphi^{(j-1)} , \1^{(j-1)} \rangle}
\mu (H_{k,1}^{(j-1)}).
\end{equation}
Then the set
$H_{k,1}^{(j-1)}$ is stable
for
\[
 l \geq \max
\{ \ord (x_1^{(j-1)} \circ \varphi^{(j-1)} ) , \dots ,\ord (
x_d^{(j-1)} \circ \varphi^{(j-1)} ), \ord (y_{j-1} \circ
\varphi^{(j-1)} ) \},
\]
where the right-hand side of this formula is independent of the
choice of $\varphi \in H_{k,1}$. It follows that the motivic
measure of $\mu (H_{k,1} ^{(j-1)})$ is equal to
\begin{equation} \label{me2}
 \mu (H_{k,1}^{(j-1)}) = [\pi_{l} (H_{k,1}^{(j-1)} ) ] \,    \L^{ -l (d+1) }.
\end{equation}
The non-angular coefficients of the  the expansions of the arc
$\varphi^{(j-1)}$ with respect to the coordinates $x_1^{(j-1)} ,
\dots, x_d^{(j-1)} , y_{j-1}$, contribute to  the class  $[\pi_{l}
(H_{k,1}^{(j-1)} ) ]$ with the factor
\begin{equation} \label{me3}
 \L^{l(d+1)
- \langle \orden^{(j-1)} \varphi^{(j-1)} , \1^{(j-1)} \rangle }.
\end{equation}

 It remains to compute the contribution of the angular
coefficients.

 We assume
first that $j \leq g$. By (\ref{Exc-j}) we get that $(f \circ \Pi_1
\circ \cdots \circ \Pi_{j-1})_{| Z_{\bar{\r}_{j-1}} } = x^{e_{j-2} \g_{j-1}} f^{(j-1)}$
where
\[
f^{(j-1)} = (y_{j-1}^{n_j} - c x^{ n_j (\l_j - \l_{j-1} ) }
)^{e_j} + \cdots, \quad \mbox{ and} \quad f^{(j-1)}_j =
(y_{j-1}^{n_j} - c x^{ n_j (\l_j - \l_{j-1} ) } ) + \cdots,
\]
for some $ c \ne 0$ and where  the other terms have exponents
outside the compact edge of the Newton polyhedron. By Proposition
\ref{liftj} the vector $\orden^{(j-1)} (\varphi^{(j-1)} ) $ belongs
to $\bar{\r}_{j}$ hence  we deduce that $ (f^{(j-1)})_{|
\orden^{(j-1)} (\varphi^{(j-1)})} = ( (f^{(j-1)}_{j})_
{|\orden^{(j-1)} (\varphi^{(j-1)})}) ^{e_j} $ and
$(f^{(j-1)}_{j})_{|\orden^{(j-1)} (\varphi^{(j-1)})} = y_{j-1}^{n_j}
- c x^{ n_j (\l_j - \l_{j-1} ) }$.

Since $\varphi$ is of depth $j$ we get that $y_j \circ
\varphi^{(j)}$ is a series with non zero constant term. By
(\ref{Exc-j}) we have that $(f \circ \Pi_1 \circ \dots \circ
\Pi_{j-1} ) _{| Z_{\bar{\r}_{j-1}}} = x^{e_{j-2} \g_{j-1}}
f^{(j-1)}$. By Proposition
\ref{lift1} we obtain the equality
$
\ac (f \circ \varphi) = \left( x^{e_{j-2} \g_{j-1}} ( (f^{(j-1)}
 _{j})_{|\orden^{(j-1)} (\varphi^{(j-1)})} )
^{e_j} \right) (\ac^{(j-1)} (\varphi^{(j-1)}))$.
 It follows that angular coefficients contribute to the class    $[\pi_{l}
(H_{k,1}^{(j-1)} ) ]$ with the factor defined by the class of
\begin{equation} \label{me5}
\{  \ac^{(j-1)}  (\varphi^{(j-1)} ) \in T_{N_{j-1}'} \mid
x^{e_{j-2} \g_{j-1}} (f^{(j-1)})_{| \orden^{(j-1)} ( \varphi^{(j-1)} )
} (\ac^{(j-1)}  (\varphi^{(j-1)}) )  =  1. \}
\end{equation}
In Lemma \ref{igual} below we prove that this class is equal to
$c_1 (\theta) (\L-1)^{\dim \theta -1}$.

If $\varphi$ is of depth $g+1$ then we get that $ \ac (f \circ
\varphi) = \ac(  (x^{e_{g-1} \g_{g}}  y_g) \circ \varphi^{(g)})$.
In this case the contribution of the
angular coefficients to the class $[\pi_{l} (H_{k,1}^{(j-1)} ) ]$
is $(\L -1)^{d}$.

Summarizing this discussion,  by (\ref{me3}) we obtain that
\begin{equation} \label{me4}
[\pi_{l} (H_{k,1}^{(j-1)} ) ] = c_1 (\theta) (\L-1)^{\dim \theta
-1} \L^{ l (d+1) - \langle \orden^{(j-1)} (\varphi^{(j-1)}) ,
\1^{(j-1)} \rangle },
\end{equation}
and by (\ref{me2}) we get,
$
\mu ( H_{k,1}^{(j-1)} ) = c_1(\theta)  (\L-1)^{\dim \theta -1}
\L^{ - \langle \orden^{(j-1)} (\varphi^{(j-1)}), \1^{(j-1)} \rangle
}$.
Then, the result follows by (\ref{me1}) and (\ref{me0}).
\hfill $\square$

\begin{lemma}
\label{igual} With hypothesis and notations of Proposition  \ref{mu-k} and its proof,
 the  class of
(\ref{me5}) in $\kmon$ is equal to
\[
\left\{
\begin{array}{lll}
[ \{ (x, y) \in (\C^*)^2 \mid (y^{n_j} - x^{r_j} )^{e_j} =1 \} ]
(\L-1)^{d-1} & \mbox{ if }  &  k \in \r_j,
\\
{[ \mu_{n_j e_j} ]} (\L -1)^{d} & \mbox{ if }  &  k \in \s_j^{-} ,
\\
{[ \mu_{r_j e_j} ]} (\L -1)^{d} & \mbox{ if }  &  k \in \s_j^{+}.
\end{array}
\right.
\]
\end{lemma}
\textit{Proof.} If $k \in  \r_j$ then (\ref{me5}) defines the
hypersurface of $T_{N_{j-1}'} = T_{N_{j-1}} \times
\C^* $ with equation,
$
x^{e_{j-1} n_{j-1} \g_{j-1}} (y_{j-1}^{n_j} - c
x^{ n_j (\l_j - \l_{j-1}) } )^{e_j} = 1$,
for some $c \ne 0$.
We have that $ e_{j-2} \g_{j-1} = e_j n_{j-1} n_{j} \g_{j-1} $, hence by
(\ref{rel-semi}) we obtain
\begin{equation}       \label{semi2}
    n_{j-1} n_j \g_{j-1} +   n_j (\l_j - \l_{j-1})  = n_j \g_j.
\end{equation}
We make first the change of coordinates in the torus
defined by $y_{j-1} = x^{ -n_{j-1} \g_{j-1}} u_{d+1}$.   By  (\ref{semi2})
the set
(\ref{me5})
is defined by the equation
$(u_{d+1}^{n_j} - c x^{n_j \g_j})^{e_j} =1$.
By definition of $r_j$ there is a primitive vector $\tilde{\g}_1
\in N_{j-1}$ such that $n_j \g_j = r_j \tilde{\g}_1$ (see Definition \ref{integral}). We consider
a basis of the lattice $N_{j-1}$ of the form $\tilde{\g}_1,
\dots, \tilde{\g}_r$. Setting $u_i := x^{\tilde{\g}_i}$, for $i=1,
\dots, d$,  defines coordinates $(u_1, \dots, u_d)$ on the torus
$T_{N_{j-1}}$. Then we get coordinates   $(u_1, \dots, u_{d+1})$
on the torus
$T_{N_{j-1}'}$
such that the set
(\ref{me5}) becomes
\begin{equation}\label{semi4}
\{ u \in (\C^*)^{d+1}  \mid (u_{d+1}^{n_j} - c u_1^{r_j} )^{e_j} =
1. \}
\end{equation}
This ends the proof in this case.

If $k \in \s_j^{-}$ (resp. $k \in \s_j^{+}$) then (\ref{me5}) is the hypersurface of the torus $ T_{N_{j-1}'} $
defined by the equation,
$
x^{e_{j-2} \g_{j-1}} y_{j-1}^{n_j
e_j} = 1$,
(resp. $x^{e_jn_j \g_j } = 1$).
In the coordinates   $(u_1, \dots, u_d, u_{d+1})$ defined above
(\ref{me5}) becomes
$
\{ u \in (\C^*)^{d+1}  \mid u_{d+1}^{n_j e_j} = 1  \} $ 
( resp.  $ \{ u \in (\C^*)^{d+1}  \mid u_1^{r_j e_j} = 1 \}$)
and we are done.
  \hfill $\square$

\section{The proofs of the main results}\label{zf}

We apply first the results of the previous sections to give a
formula for the motivic zeta function of a q.o.~ polynomial in
terms of the characteristic data. The following Proposition is
crucial (see Definition \ref{xi}).

\begin{proposition}             \label{zf1}
We have that
$
 Z(f, T)_0 = \sum_{j=1}^{g+1} \sum_{\theta \in \Theta_j} \sum_{k \in \stackrel{\circ}{\theta} \cap \Z^{d+g+1}}
\mu( H_{k,1} ) T^{\eta(k)}$.
\end{proposition}
\textit{Proof.} The truncation map $\pi_n :
\mathcal{L}(\C^{d+1})_0 \rightarrow \mathcal{L}_n(\C^{d+1})_0$ is
a surjection because $\C^{d+1}$ is smooth. Let us denote by
$\mathcal{Z}_{n,1} \subset \mathcal{L}(\C^{d+1})_0$ the preimage of
$\mathcal{X}_{n,1}$ under $\pi_n$. The sets $\mathcal{Z}_{n,1}$
are stable at level $n$ and thus $\mu(\mathcal{Z}_{n,1})= [\pi_n
(\mathcal{Z}_{n,1})] \L^{-n(d+1)} = [\mathcal{X}_{n,1}]
\L^{-n(d+1)}$. We deduce that $Z(f,T)_0 = \sum_n
\mu(\mathcal{Z}_{n,1})T^n$.

For any integer $n \geq 0$ the sets $ \mathcal{Z}_{n,1}$ and
$\mathcal{Z}_{n,1} \cap \mathcal{H}$ have the same motivic measure
since their difference is a subset of arcs contained in a
hypersurface (see Lemma \ref{motmes}).  By Theorem \ref{trop}  we
have the partition
\[
\mathcal{Z}_{n,1} \cap \mathcal{H} = \bigsqcup_{j=1}^{g+1}
\bigsqcup_{\theta \in \Theta_j} \bigsqcup_{{ { k \in
\stackrel{\circ}{\theta} \cap \Z^{d+g+1}}} }^{ \eta(k)=n
}H_{k,1}. \] This partition may be non-finite. The result follows
by using the properties of the motivic measure (see Proposition
\ref{motmes}). \hfill $\square$

\begin{lemma}\label{P4lim} We have that 
$ \lim_{T \rightarrow \infty} (\L -1)^{\dim \theta -1} 
 S_\theta =   (-1)^{\ell_g + 1} (\L -1 )^{\ell_g}$
if $\theta = \s_{g+1}$. If   $1 \leq j \leq g$ we get that:
$$\lim_{T \rightarrow \infty} (\L -1)^{\dim \theta -1}  S_\theta = \left \{ \begin{array}{ll}
(-1)^{\ell_j }({\L}-1)^{{\ell_j -1}}& \hbox{if }  \theta = \s^+_j,
\\
(-1)^{\ell_j}({\L}-1)^{\ell_j -1}& \hbox{if }  \theta = \rho_j,
\\
(-1)^{\ell_{j-1} + 1 }({\L}-1)^{\ell_{j-1}}& \hbox{if } \theta = \s^-_j.
\end{array} \right .$$
\end{lemma}
\textit{Proof.} The functions $S_\theta$  from Theorem
\ref{main-zeta} are rational functions in
$\mathcal{M}^{\hat{\mu}}_\C[[T]]_{\rm sr}$ of degree $0$ with
respect to $T$.   If $1 \leq j \leq g$ the number of primitive
vectors $v \in \Z^{d+g+1}$ generating an edge of $\rho_j$ such
that $v_{d+g+1} =0 $ is equal to the number $d - \ell_j$ of coordinates
of $\lambda_j$ which are zero (see Lemma \ref{not00}). In the case
of the cone $\sigma^+_j$ this number is equal to $d+ 1 - \ell_j$ because
we have to consider the edge spanned by $\epsilon_{d+j}$. Then 
we apply Lemma \ref{simplicial}.

In the case of the cone $\s_j^-$ we use  Lemma \ref{minus}. It
follows by Lemma \ref{simplicial} that $\lim_{T \rightarrow
\infty}  (\L -1)^d S_ {\hat{\r}'} = (-1)^{\ell_{j-1} + 1} (\L-1)^{\ell_{j-1}}$
and $\lim_{T \rightarrow \infty} (\L -1)^d S_{\hat{\bar{\r}}_j}   = -
 \lim_{T \rightarrow \infty} (\L -1)^d S_{\hat{\bar{\s}}_{j}^+} $,
since the last coordinate of the primitive vector $\Upsilon_j
(0,1) \in \Z^{d+g+1}$ is nonzero. Then the assertion follows by equality
 (\ref{minus-rel}).
\hfill $\square$

We use the notations introduced in Section \ref{results}.

\bigskip

\textbf{Proof of Theorem \ref{main-zeta}.} We deal first with the
formula for the motivic zeta function.
It  is enough to compute the sum of the auxiliary series
 $ R_\theta := \sum_{k \in \stackrel{\circ}{\theta} \cap \Z^{d+g+1}}    \mu( H_{k,1} ) T^{\eta(k)} $,
for $\theta \in \Theta_j$ and $1 \leq j  \leq g+1$ (see
Proposition \ref{zf1}).
By Proposition  \ref{mu-k} we get that
$R_\theta =       c_1( \theta) (\L-1)^{\dim \theta -1}
\sum_{k \in \stackrel{\circ}{\theta} \cap \Z^{d+g+1}}
\L^{-\xi_j (k)} T^{\eta(k)}$. 
By  Lemma \ref{gf} and Notation \ref{gf2} we deduce that
$
 R_\theta =  c_1( \theta) (\L-1)^{\dim \theta -1}  S_\theta$, 
and then the conclusion follows.

The formula for the motivic Milnor fiber appears now as a
consequence of the formula for the motivic zeta function from Lemma
\ref{P4lim} and Definition \ref{csigma}.

We  prove now the formula for the naive motivic zeta function. We
consider the sets $H_{k} := \{ \varphi \in \mathcal{H} \mid \ord
(F \circ \varphi) = k\}$, for $k \in \Z^{d+g+1}_{\geq 0}$.  By
arguing as in Proposition \ref{zf1} we get the formula
\[
Z^{\mathrm{naive}} (f, T) = \sum_{j=1}^{g+1} \sum_{\theta \in
\Theta_j} \sum_{k \in \stackrel{\circ}{\theta} \cap \Z^{d+g+1}}
\mu( H_{k} ) T^{\eta(k)}.
\]
If $\theta \in \Theta_j$ the class of the set
$$\{ \ac^{(j-1)} (\varphi^{(j-1)} ) \in T_{N_{j-1}'} \mid
x^{e_{j-2} \g_{j-1}} (f^{(j-1)})_{| \orden^{(j-1)} ( \varphi^{(j-1)}
) } (\ac^{(j-1)} (\varphi^{(j-1)}) )  \not =  0 \}$$ in $\kvar$
equals
\[
\left\{
\begin{array}{lll}
(\L -1)^{d+1} & \mbox{ if }  &  k \in \s_j^{+}, \s_j^{-},\s_{g+1}
\\
(\L-1)^{d} (\L-2)  & \mbox{ if }  &  k \in \r_j.
\end{array}
\right.
\]
This allows  us to prove that the motivic measure of $H_k$ is $\mu
(H_{k}) = c(\theta) (\L -1)^{\dim \theta -1} \L^{- \xi_j (k)}$
following the argument of Proposition \ref{mu-k}.

To prove the formula for the topological zeta
function we need  Proposition \ref{ztop-rel}. 

If $1 \leq j \leq g+1$ and if $\theta \in \Theta_j$ is a $(d+1)$-dimensional simplicial
cone, that is, $\theta =\s_j^+$ if $j< g+1$ or $\theta = \s_{g+1}$ if $j = g+1$.
By definition, using notations  Lemma \ref{gf} and Definition \ref{xi},
if $v_1, \dots, v_{d+1} \in \Z^{d+g+1}$ are the primitive vectors generating the rays of
$\theta$
we get the formula:
\begin{equation}   \label{rel-1}
   c(\theta)({\L}-1)^{\dim \theta -1}  S_\theta  = (\sum_{v \in
   D_\theta} \L^{-\xi(v)-s \eta(v)})
\prod_{i=1}^{d+1} ({\L}-1) { (1 - \L^{-\xi(v_i)-s
\eta(v_{i})})^{-1}}
\end{equation}
By applying $\chi_{\mathrm{top}}$ to (\ref{rel-1}) provides the term
$J_{\theta}(f,s)$ defined in  \ref{ztop-not}, since
$ [\P^{a-1}_\C]^{-1} = (\L -1) ( \L^a - 1)  $
 and  $\chi_{\mathrm{top}}([\P^{a-1}_\C])=a$.
We deduce similarly that if $\theta = \r_j$ for $1 \leq j  \leq g$
then the result of applying    $\chi_{\mathrm{top}}$  to (\ref{rel-1}) is equal to
$-J_{\r_j} (f,s)$.  The negative sign appears when applying $\chi_{\mathrm{top}}$ to
the term $(\L -2)$ of
$c(\r_j)$.
Finally, we study the case of $\theta = \s_j ^-$ for $1\leq j \leq g$, which may be not simplicial in general.
We use Lemma \ref{minus}. By (\ref{minus-rel}) we have that
\begin{equation}
\label{rel-2}
  c(\theta)({\L}-1)^{\dim \theta -1}  S_\theta (\L^{-s})   = (\L - 1)^{d+1}  ( S_{\hat{\r}'} (\L^{-s})
    - S_{\hat{\bar{\r}}_j} (\L^{-s})   -
S_{\hat{\bar{\s}}_{j}^+} (\L^{-s})   ) .
\end{equation}
By Lemma \ref{minus} we get that
$\mathrm{mult} (\hat{\r}') = \mathrm{mult} (\r_{j-1})$, 
$\mathrm{mult} ( \hat{\bar{\s}}_{j}^+  ) = \mathrm{mult} ( {{\s}}_{j}^+ )  $. 
By Lemma \ref{minus}  we have the equalities $\xi_j (\Upsilon_j (0,1)) = 1$ and 
$\eta (\Upsilon_j (0,1))= n_j \dots n_g$. Taking into account
these observations we deduce that:
\[
\chi_{\mathrm{top}}      ( (\L - 1)^{d+1}    S_{\hat{\r}'} )
(\L^{-s})
 =  \frac{J_{\r_{j-1}} (f, s)}{1 + n_j \cdots n_g  s},
\quad \chi_{\mathrm{top}}      ( (\L - 1)^{d+1}
S_{\hat{\bar{\s}}_{j}^+} (\L^{-s})   ) = \frac{J_{\s_j^+} (f, s)
)} {1 + n_j \cdots n_g  s},
\]
where the factor $(1 + n_j \cdots n_g  s)^{-1}$ corresponds to the vector $\Upsilon_j (0,1)$.
Notice also that
$\chi_{\mathrm{top}}      ( (\L - 1)^{d+1}  S_{\hat{\bar{\r}}_j}) = 0$. See Notation \ref{ztop-not}.
\hfill $\square$

\begin{remark}            \label{minus-s}
We can expand  $\chi_{\mathrm{top}}      ( (\L - 1)^{d+1}     S_{{{\s}}_{j}^-} )$  in terms of
a simplicial subdivision of the cone $\s_{j}^-$. Without loss of generality we assume in addition
that the edges of the subdivision are exactly the edges
of the cone $\s_{j}^-$. Such a subdivision always exists (see \cite{Ewald} Chapter 5, Theorem 1.2).
If $\theta_1, \dots, \theta_l$ are the $(d+1)$-dimensional cones in the subdivision we argue as
before to obtain $
 \chi_{\mathrm{top}} ((\L -1)^{d+1} S_{\s_{j}^-} ) = \sum_{i =1}^l {\mathrm{mult} (\theta_i) }
{\prod_{v_i} (\xi_j (v_i) + s \eta (v_i))^{-1} }$,
where $v_i \in \Z^{d+g+1}$ runs through the primitive vectors generating the edges of $\theta_i$, for $i=1, \dots, l$.
Compare also with Lemma 5.1.1 \cite{DL-JAMS}.
\end{remark}

\section{Comparison with the method of Newton maps}    \label{comparison}

The method of Newton maps allows to describe the parametrization of a q.o.~hypersurface singularity.
It has been used by Artal
 \textit{et al.} \cite{ACNLM-AMS}
to study the poles of the 
motivic zeta function of a quasi-ordinary hypersurface singularity. We describe this method in geometric terms through the toric embedded resolution (see \cite{GP-Fourier}).
For simplicity in the exposition we restrict ourselves to the case of irreducible germs.
We keep notations of the  Sections  \ref{qoh}, \ref{sec-toric} and \ref{resolution}.
As in Section \ref{resolution} we assume that $y = f_0$ (see Definition \ref{semi}).

Recall that we use the following notations: 
the canonical basis of $\Z^d = N_0$ is $\epsilon_1, \dots, \epsilon_d$, 
the cone spanned by $\epsilon_1, \dots, \epsilon_d$ is $\r$, and 
$\varepsilon_1, \dots, \varepsilon_d$ denotes
the dual basis of the dual lattice $M_0$.

\begin{definition}
We define for $1 \leq j \leq g$ a sequence of vectors $\{ \epsilon_i^{(j)} \}_{i=1}^d$ 
defining basis of a lattice    $\tilde{N}_j \subset N_0$;  then we denote by $\tilde{M}_j \supset M_0$  
the dual lattice of 
$\tilde{N}_j$ equipped with the dual basis 
$\{ \varepsilon_i^{(j)} \}_{i=1}^d$ of  $\{ \epsilon_i^{(j)} \}_{i=1}^d$. 
The definition is given in term of the characteristic exponents by induction on $j$.   
We expand 
$\l_i = \sum_{i=1}^d   \frac{q_i^{(1)}}{p_i^{(1)}} \varepsilon_i $
 where $\frac{q_i^{(1)}}{p_i^{(1)}}$ is an irreducible fraction for $i=1, \dots, d$.
Then we  set 
\[ \epsilon_i^{(1)} =  p_i ^{(1)} \epsilon_i, \mbox{ for } i =1, \dots, d. \]
This defines the lattice $\tilde{N}_1$.   
With these notations we obtain
\begin{equation}  \label{change}
\varepsilon_i = p_i^{(1)} \varepsilon_i^{(1)}, \mbox{ for } i =1, \dots, d, \mbox{ and }
\l_1 = q_1^{(1)}
\varepsilon_1^{(1)} + \cdots +  q_d^{(1)}
\varepsilon_d^{(1)}, 
\end{equation}
hence we get that $M_1 \subset \tilde{M}_1$ and $\tilde{N}_1 \subset N_1$. 
We suppose that $1 < j \leq g$ and $\tilde{M}_{j-1}$ together with its basis
$\{ \varepsilon_i^{(j-1)} \}_{i=1}^d$ have being defined by induction. 
Then we expand the vector $\l_{j} -\l_{j-1} $ in the basis  $\{ \varepsilon_i^{(j-1)}\}_{i=1}^d$. 
We obtain
$\l_j - \l_{j-1} = \sum_{i=1}^d \frac{q_i^{(j)}}{p_i^{(j)}} \varepsilon_i^{(j-1)}$, 
where $\frac{q_i^{(j)}}{p_i^{(j)}}$
are irreducible fractions. Then, with these notations we define the basis $\{ \epsilon_i^{(j)} \}_{i=1}^d$ of 
the lattice $\tilde{N}_j$ by setting: 
 \begin{equation}    \label{vect-2}
   \epsilon_i^{(j)} =  p_i ^{(j)} \epsilon_i^{(j-1)}, \mbox{ for } i =1, \dots, d. 
 \end{equation}
We get that 
 \begin{equation}  \label{changej}
\varepsilon_i^{(j-1)}  = p_i^{(j)} \varepsilon_i^{(j)}, \mbox{ for } i =1, \dots, d, \quad 
\l_j - \l_{j-1}  = q_1^{(j)}
\varepsilon_1^{(j)} + \cdots +  q_d^{(j)}
\varepsilon_d^{(j)},
\end{equation}
$M_j \subset \tilde{M}_j$ and $\tilde{N}_j \subset N_j$. 
\end{definition}

By definition the cone $\r$ is regular for the lattice $\tilde{N}_j$. Hence
the $d$-dimensional toric variety $Z_{{\r}, \tilde{N}_j}$ is isomorphic to $\C^d$.
We have a $\Z$-bilinear pairing $N_j \times \tilde{M}_j \to  \Q$, $(\nu, \g)  \mapsto \langle \nu, \g \rangle$
which extends the duality pairings $\tilde{N}_j \times \tilde{M}_j \to \Z$ and $N_j \times M_j \to \Z$.
The inclusion of lattices $\tilde{N}_j \subset N_j$ induces a proper equivariant map
\begin{equation} \label{quotient}
\kappa_j:  Z_{\r, \tilde{N}_j} \longrightarrow Z_{\r, N_j},
\end{equation}
  which is the projection  for the quotient of  $Z_{{\r}, \tilde{N}_j}$
with respect to the natural action of the finite abelian group
$G_j := N_j/ \tilde{N_j}$.
This action is defined on monomials and extended by linearity.
If $\nu + \tilde{N}_j \in N_j /  \tilde{N}_j$ and $\g \in \tilde{M}_j$ then the action is
given by
$
(\nu + \tilde{N}_j) \cdot x^\g := e^{2 \pi i \langle \nu, \g \rangle}  x^\g$.
The subset of invariant monomials by this action is
$\{  x^\g \mid \g \in M_{j} \}$.
The restriction of the map (\ref{quotient}) to the torus $T_{\tilde{N}_j} \to T_{N_j}$ is an unramified
covering of degree $|G_j|$ and the kernel of this map is isomorphic to $G_j$.
See \cite[Corollary 1.16]{Oda}.

\begin{lemma}     \label{vect}
The primitive integral vectors of the cone $\bar{\r}_j$ with respect to the lattice $N_{j-1} \times \Z$
(and also with respect to the lattice $\tilde{N}_{j-1} \times \Z$) 
are      $( p_i ^{(j)} \epsilon_i^{(j-1)} ,  q_i^{(j)} )$ for $i =1, \dots, d$. 
\end{lemma}
\textit{Proof}. 
In Lemma \ref{iso-semi} and Remark \ref{iso-semi-i} we introduced for $1 \leq j \leq g$ the 
lattice homomorphisms $\phi_j : M_{j-1} ' \to M_j$ 
in such a way that we can identify $M_j$ with the quotient lattice $M_{j-1}' / \mathrm{ker} (\phi_j)$. 
The dual lattice homomorphism $\phi_j^* : N_j \to N_{j-1}'$ is injective and 
by (\ref{changej}) it 
verifies that 
\begin{equation} 
\phi_j^* ( \epsilon_i^{(j)} ) = ( p_i ^{(j)} \epsilon_i^{(j-1)} ,  q_i^{(j)} ), \quad  i =1, \dots, d,
\end{equation}
are the primitive integral vectors for the lattice $N_{j-1}'$  (and also for $\tilde{N}_{j-1} \times \Z$)
in the edges of the 
cone $\bar{\r}_j \subset (N_{j-1}')_\R$. \hfill $\square$

\begin{remark}
By the arguments in Lemma \ref{vect} we can identify $N_j$ with its image  
 $\mathrm{Im} (\phi_j^*) =  (\mathrm{ker} (\phi_j))^\perp$,
which is equal to the subset of  $N_{j-1}'$ in the linear hull of the cone $\bar{\r}_j$.
We deduce that  $\phi_j^* (\r) = \bar{\r}_j$   hence 
 the map $\phi_j^*$ defines an isomorphism of toric varieties 
$Z_{\r, N_j} \to Z_{\bar{\r}_j, \phi_j^* (N_j)}$. 
\end{remark}

The following  maps are defined in terms of the characteristic exponents for $j =1 , \dots, g$,
with respect to suitable choices of coordinates. 
\begin{definition}
 The Newton map 
$\mathcal{N}_j : \C^d \times \C^* \to \C^{d+1} $ is the monomial map 
\[ (x_1^{(j)},  \dots, x_d^{(j)}, z_j) \to (x_1^{(j-1)}, \dots, x_d^{(j-1)}, \tilde{y}_{j-1}) \]
given by:
   \begin{equation}  \label{newton-map}    
\left\{
\begin{array}{lcl}
x_i^{(j-1)}   & = & (x_i^{(j)})^{p_i^{(j)}},  \mbox{ for } i =1, \dots, d
\\
\tilde{y}_{j-1} & =  & (x_1^{(j)})^{q_1^{(j)}} \cdots (x_d^{(j)})^{q_d^{(j)}} z_j.
\end{array}
\right.
\end{equation}                                     
\end{definition}

We express the Newton map $\mathcal{N}_j$  in terms of the restriction of the map $\Pi_j$  to the
chart $Z_{\bar{\r}_{j} , N_{j-1}'}$. We  denote
this restriction also by $\Pi_j$.
By Lemma \ref{iso-semi} and Remark \ref{iso-semi-i} 
the coordinate ring of $Z_{\bar{\r}_j , N_{j-1}'} $
is isomorphic to $\C [\r^\vee \cap M_j ] [w_j ^{\pm 1} ]$ and
$Z_{\bar{\r}_j , N_{j-1}'} $ is equivariantly isomorphic to  $Z_ {{\r}, {N}_j} \times \C^*$.
Setting $w_j = z_j ^{n_j} $ defines a map
\[
 \psi_j : Z_{{\r}, N_j}  \times \C^* \to    Z_{{\r}, N_j}  \times \C^*
,\]  which is an
unramified covering of degree $n_j$. This map also appears in \cite{GP-Fourier} Section 3.2.1.

\begin{proposition}  \label{pro-newton}
The Newton map $\mathcal{N}_1$ is equal to the composite
\[
 Z_{{\r}, \tilde{N}_1}  \times \C^* \stackrel{\kappa_1 \times Id}{\longrightarrow}
Z_ {{\r}, {N}_1} \times \C^*
\stackrel{\psi_1}{\longrightarrow}
Z_ {{\r}, {N}_1} \times \C^*
\stackrel{\Pi_1}{\longrightarrow}       \C^{d+1}.
\]
\end{proposition}
\textit{Proof.} 
By using the formula (\ref{change}) the relation $y = (x_1^{(1)})^{q_1^{(1)}}
\cdots (x_d^{(1)})^{q_d^{(1)}} z_1$,
can be rewritten as $z_1 = y x^{-\l_1}$. Then apply that $w_1 = z_1^{n_1}$.
It follows that then composite map above is defined by the  homomorphism of semigroups 
$\Z^{d+1}_{\geq 0} \to ({\r}^\vee  \cap \tilde{M}_1) \times \Z$
which maps     $i$-th canonical basic vector of $\Z^{d+1}$
to  $( p_i^{(1)} \varepsilon_i^{(1)}, 0) $  if $ 1\leq i \leq d$
or to    $\sum_{i=1}^d q_i^{(1)}
(\varepsilon_i^{(1)}, 0) + (0,1)$ if $i=d+1$. 
Then  the result follows from the definitions.  
\hfill $\square$

\begin{remark}
In the  \cite[Definition 3.20]{ACNLM-AMS} the Newton map is defined by replacing 
$z_1$ in (\ref{newton-map}) by $z_1 - \a_1$ for some suitable $\a_1 \in \C^*$. 
This translation corresponds in our definition to the choice of  
a suitable point $o_1'$ in the fiber $\psi_1^{-1} (o_1)$. 
\end{remark}

We deduce some consequences of Proposition \ref{pro-newton} and the results stated in Section \ref{resolution}.

We know that the strict transform $S^{(1)}$ of $S$ by the modification $\Pi_1$ is a germ at the point
$o_1 \in  \{ 0 \} \times \C^* \subset Z_ {{\r}, {N}_1} \times \C^* $. It follows that the fiber
 of $o_1$ by $\psi_1$ consists of $n_1$ different points.
Let us fix one point ${o}_1' \in \psi_1^{-1} (o_1)$.     The map $\psi_1$ defines an isomorphism of germs
$(Z_ {{\r}, {N}_1} \times \C^* , o_1')  \to (Z_ {{\r}, {N}_1} \times \C^*, o_1)$.
The preimage by $ \kappa_1 \times Id $ of the point  $o_1'$ has only one point which we denote by
$\tilde{o}_1$.
We get that the germ of the strict transform $\tilde{S}^{(1)}$ of $(S,0)$ 
by the Newton map $\mathcal{N}_1$
at the point $\tilde{o}_1$
coincides with the  germ at $\tilde{o}_1$ of the preimage of   $(S^{(1)},o_1)$ by
the map $\psi_1 \circ (\kappa_1 \times Id)$. 
The germs defined in this way
by picking different points in the fiber $\psi^{-1} (o_1)$ are isomorphic.
Proposition \ref{pro-newton} also implies that the germ $(\tilde{S}^{(1)}, \tilde{o}_1)$ is
invariant by the action of the
group $G_1$.
The germ   $(S^{1},o_1)$ is a t.q.o.~ singularity (see Remark \ref{r56}). 
We deduce that  $(\tilde{S}^{(1)}, \tilde{o}_1)$
is a germ of q.o.~hypersurface singularity.
If $g=1$ then Proposition \ref{pro-newton}
implies that $(\tilde{S}^{(1)}, \tilde{o}_1)$ is smooth and
the partial toric resolution procedure and the Newton map procedure
end after one step. 

Assume that $g>1$. We take a local coordinate $\tilde{y}_1$ at $\tilde{o}_1$, which
defines the germ $(\tilde{S}^{(1)}_1, \tilde{o}_1)$, that is,  the strict transform of the
semi-root $(S_1,0)$. Notice that  $(\tilde{S}^{(1)}_1,
\tilde{o}_1)$ is invariant by the action of $G_1$.
Notice that the germ $(Z_{\r, \tilde{N_1}} \times \C^*, \tilde{o}_1)$ is isomorphic to 
$(\C^{d+1}, 0)$ with coordinates $x_1^{(1)}, \dots, x_d^{(1)}, \tilde{y}_1$. 
In Section \ref{resolution} we proved also that the germ 
 $(Z_{\r, {N_1}} \times \C^*, {o}_1)$
is isomorphic to 
$(Z_{\r, N_1} \times \C , 0)$
with the coordinates given by $y_1$ and $x^\g$ for $\g \in \r^\vee \cap M_1$.
From Proposition    \ref{pro-newton}, we deduce
that the Newton polyhedron of the series defining the 
strict transform 
$(\tilde{S}^{(1)}, \tilde{o}_1)$ 
with respect to the coordinates  $x_1^{(1)}, \dots, x_d^{(1)}, \tilde{y}_1$
coincides with 
the Newton polyhedron of  the series $f^{(1)} \in \C \{ \r^\vee \cap M_1 \} \{y_1 \}$ 
defining the strict transform $(S^{(1)}, o_1)$. 
It follows that both polyhedra define the same dual subdivision $\Sigma_2$.
By Lemma \ref{vect} the primitive integral vectors in edges of 
the cone $\bar{\r}_2 \in \Sigma_2$ 
are the same for the lattices $\tilde{N}_1 \times \Z$ and $N_1 \times \Z$. 
The Newton map $\mathcal{N}_2$ is defined in terms of 
these primitive vectors. We deduce, as in Proposition \ref{pro-newton}, that 
$
\Pi_2   \circ  \psi_2   \circ (\kappa_2 \times Id)  = (\kappa_1 \times Id)  \circ \mathcal{N}_2$.
If $g >2$ we iterate this procedure as in Section \ref{resolution}. We
can extend these notions to the level $1 < j < g$. 
     At the level $j-1$ we choose one point $\tilde{o}_{j-1}$ in the
          fiber $(\psi_{j-1} \circ (\kappa_{j-1} \times Id) )^{-1} (o_{j-1})$.
Then the Newton map 
$
\mathcal{N}_{j} : Z_{\r, \tilde{N}_{j}} \times \C^* \to Z_{\r, \tilde{N}_{j-1}} \times \C
$
can be also defined in terms of the Newton polyhedron of the strict transform of $f$
at the point $\tilde{o}_{j-1}$, with respect to 
the coordinates $x_1^{(j-1)}, \dots, x_d^{(j-1)}, \tilde{y}_{j-1}$ of $Z_{\r, \tilde{N}_{j-1}} \times \C$, 
where
$\tilde{y}_{j-1}$ denotes the strict transform of the semiroot $f_{j-1}$.

The interplay between the Newton maps and the toric resolution is summarized by:
\begin{proposition}   With the notations above for $1 \leq j \leq g$ we have that 
  \[
 \Pi_j   \circ  \psi_j \circ  (\kappa_j \times Id)  = (\kappa_{j-1} \times Id)  \circ \mathcal{N}_j.
\]
The strict transform $\tilde{S}^{(j)}$ at the point $\tilde{o}_j$ of the q.o. hypersurface $(S,0)$ by 
the iteration of the first $j$ Newton maps is invariant by the action of $G_j$.
The image of the germ $(\tilde{S}^{(j)}, \tilde{o}_j)$
by the map $(\kappa_j \times Id)  \circ  \psi_j$ is
equal to $(S^{(j)}, o_j)$ for $j=1, \dots, g$.
The germ $(\tilde{S}^{(g)}, \tilde{o}_g)$ is smooth and
the map $(\tilde{S}^{(g)}, \tilde{o}_g) \to (S^{(g)}, o_g)$ is
the \textit{canonical orbifold map} associated to the normalization $(S^{(g)}, o_g)$ of the germ $(S, 0)$.
\end{proposition}

\begin{remark}
    See \cite{PPP-Duke} for a topological description of  the canonical orbifold map.
\end{remark}

\begin{remark} \label{false2}
The germ  $(S^{(1)},o_1)$ is irreducible however if $d \geq 2$ and $|G_1| > 1$ the germ
$(\tilde{S}^{(1)}, \tilde{o}_1)$
 it is not analytically
irreducible in general (see Example \ref{false}). 
 It may happen that   $(\tilde{S}^{(1)}, \tilde{o}_1)$
have irreducible factors which are not invariant by the action of $G_1$. This
phenomenon of \textit{``false reducibility''} has been analyzed and explained in 
a combinatorial way in \cite{MGV}.
\end{remark}

\subsection*{Comparison of the sets of candidate poles}\label{comparepoles} 

Artal et al. define in \cite[Definition 3.28]{ACNLM-AMS} a set of pairs 
 $CP (f,\omega)$ where  $f$ is a quasi-ordinary polynomial, $\omega$ 
is certain differential form  and 
the  coordinates are choosen in a suitable way.
We assume here that $f \in \C \{ x_1,...,x_d \} [y]$ is
an irreducible q.o.~polynomial as in Section  \ref{resolution}, $y= f_0$ and 
 $\omega = dx_1 \wedge \cdots \wedge dx_d \wedge dy$.   
\begin{definition} 
We define  first 
$CP (f,\omega)$ as the list of pairs 
$\{ ( B_i^{(j)}, b_i^{(j)} ) \}_{i=1, \dots, d}^{j=1, \dots, g} \cup \{ (1,1) \}$,
where for $1 \leq i \leq d$ and $1 < j \leq g$ we set 
\[ 
B^{(1)}_i := e_0 {q}^{(1)}_i,  \quad B^{(j)}_i := p^{(j)}_i B^{(j-1)}_i + e_{j-1} {q}^{(j)}_i   \mbox{ and }
b^{(1)}_i:= {p}^{(1)}_i + {q}^{(1)}_i ,   \quad b^{(j)}_i:= {p}^{(j)}_i b_i^{(j-1)}  
+ {q}^{(j)}_i. 
\]
\end{definition}

The number $B_i^{(j)}$  (resp. $b_i^{(j)} -1$) 
is equal to the order of vanishing of the pull-back  of $f$ (resp. of $\omega$) under 
the sequence of first $j$ Newton maps 
along the divisor $x_i^{(j)} = 0$. These formulas follow easily
by induction using the definition of the Newton maps (see \cite{ACNLM-AMS}). 
The terms $1 - \L^{-a} T^A$ for $(A,a ) \in  CP (f,\omega)$
form a set of candidate poles of the  
zeta function $Z^{\mathrm{naive}} (f, T)$ \textit{as computed} in \cite{ACNLM-AMS}.

In order to compare with our results we consider the cone $\r$ embedded in 
$\R^{d+g+1} = \R^d \times \R^{g+1}$ as $\r_0 = \r \times \{ 0 \}$ (see Notation \ref{cones}). 
By Lemma \ref{pos-trop} if  $ 1 \leq i \leq d$ the vector  $\nu_i^{(j)} := 
{s}_j ( \epsilon_i^{(j)} )$ is primitive for the lattice 
$\Z^{d+g+1}$ and defines an edge of the cone $\r_j$.

\begin{lemma}        \label{comparison-2}
We have that $B_i^{(j)} = \eta (\nu_i^{(j)})$ and $b_i^{(j)} = \xi_j (\nu_i^{(j)})$ for $1 \leq i \leq d$ and 
$1 \leq j \leq g$. 
\end{lemma}
\textit{Proof.} We prove these equalities by induction on $j$. If $j=1$ they hold by the definitions. 
Assume that $j >1$ and that the result holds for $j-1$. 

We deal first with the equality for $B_i^{(j)}$. We argue as in the proof of Theorem \ref{trop}. 
The vector $\nu_i^{(j)} = \ord F \circ \varphi$ for some arc $\varphi$ of depth $j$. 
By Lemma \ref{vect} and (\ref{vect-2}) we get
 $\eta( \nu_i^{(j)} ) = e_{j-1} \langle \epsilon_i^{(j)}, \g_j \rangle $ and similarly 
 $\eta( \nu_i^{(j)} ) = e_{j-2} \langle \epsilon_i^{(j-1)}, \g_{j-1} \rangle$.
Since $\epsilon_i^{(j)} = p_i^{(j)} \epsilon_i^{(j-1)} $ 
we deduce using (\ref{rel-semi}) that
\[
 \eta(\nu_i^{(j)} ) =   p_i^{(j)} e_{j-2} \langle  \epsilon_i^{(j-1)}, \g_{j-1} \rangle + 
e_{j-1} \langle \epsilon_i^{(j)} , \l_j - \l_{j-1} \rangle  = p_i^{(j)} \eta( \nu_i^{(j-1)}) 
 + e_{j-1} q_i^{(j)}. 
\]

We deal with the equality for $b_i^{(j)}$ in a similar way. We have that 
$\xi_j (\nu_i^{(j)}) = \sum_{l=1}^d \langle \epsilon_i^{(j)}, \varepsilon_l^{(0)} \rangle + \sum_{l=1}^{j-1} 
(1 - n_l) \langle \epsilon_i^{(j)}, \g_l \rangle +  \langle \epsilon_i^{(j)}, \g_j \rangle$. 
By (\ref{vect-2})  and relations (\ref{rel-semi}) 
we get:
\[
 \xi_j   (\nu_i^{(j)}) =   p_i^{(j)}  \xi_{j-1} (\nu_i^{(j-1)})  + \langle \epsilon_i^{(j)} , \g_j - n_{j-1} \g_{j-1} \rangle = 
 p_i^{(j)}  \xi_{j-1} (\nu_i^{(j-1)})  + q_i^{(j)}.
\]

In both cases the result follows by the induction hypothesis. 
\hfill $\square$

\begin{proposition}        \label{comp-2}
We have that $
CP(f,w) = \cup_{j=1}^{g}  \cup_{i=1}^d \{ (\eta ( \nu_i^{(j)} ) , \xi_j (\nu_i^{(j)} ) \} \cup \{ (1, 1) \}$. 
\end{proposition}
\textit{Proof.} The equality of these two lists is consequence of Lemma \ref{comparison-2}. We use also that 
$\epsilon_{d+g+1} \in \s_{g+1}$ verifies that 
$\eta ( \epsilon_{d+g +1}) =1  $, $\xi_{g+1} (\epsilon_{d+g +1})  =1  $. 
\hfill $\square$

We discuss now some properties related with special vectors (see Definition \ref{special}).
\begin{remark}
By Definition \ref{special} the vector $\nu_i^{(j)} \in \r_j$ is special if
the conditions $1 \leq i ' \leq d$, $q^{(j')}_{i'} =0$ for $j'=1, \dots, j-1$
and $p^{(j')}_{i'} = n_j,  q^{(j')}_{i'} =1$, hold if and only if $i = i'$ and in addition we have that 
$q^{(j+1)}_i >0$ if $j  < g$. 
If $\nu^{(j)}_i$ is a special vector then   $\nu^{(l)}_i = \epsilon_i$ 
for $ 1 \leq l \leq j-1$. By Lemma \ref{vect} and Notation \ref{cones} we get that 
$\xi_j (\nu^{(j)}_i)= n_j +1$ and $\eta(\nu^{(j)}_i)= n_j \cdots n_g$.
In this case we have that  
$(\xi_j(\nu^{(j)}_i), \eta(\nu^{(j)}_i)) \ne (\xi_j(\nu^{(j)}_l), \eta(\nu^{(j)}_l))$  
for $1 \leq l \leq d$ and $l \ne i$.
\end{remark}

\begin{definition}
The local contribution $LC_j$ at level $j$ to the naive motivic zeta function is
$$LC_j:=(\L-1)^{d+1}  \big ( S_{\s^+_j}  + S_{\s^-_j} + S_{\s^-_{j+1}} \big ) + (\L -1)^d(\L - 2) S_{\r_j},$$
where if $j = g$ one must replace   $\s^-_{g+1}$ for $\s_{g+1}$.
\end{definition}

The following result is analogous to  \cite[Proposition 5.4]{ACNLM-AMS}:
\begin{proposition} \label{LC}    If $\nu_i^{(j)}$ is a special vector 
then there exists a polynomial $R_j \in \Z [ \L^{\pm 1} ] [T]$ such that 
\[
 LC_j = R_j \prod_{r= j-1}^{j+1}  
\prod_{l=1, l \ne i }^d (1 - L^{- \xi_r (\nu_l^{(r)} ) } T^{ \eta( \nu_l ^{(r)})})^{-1} 
\prod_{r= j-1, r \ne j }^{j+1}     (1 - L^{- \xi_r (\nu_i^{(r)} ) } T^{ \eta( \nu_i ^{(r)})})^{-1}. 
\]
\end{proposition}
\textit{Proof}. 
Assume first that $j < g$. By the proof of Theorem \ref{main-zeta} 
 there exists a polynomial  $ R_j '  \in \Z[\L^{\pm 1} ]$ 
such that 
$ LC_j = R_j '  \prod_{r= j-1}^{j+1}  
\prod_{l=1}^d (1 - L^{- \xi_r (\nu_l^{(r)} ) } T^{ \eta( \nu_l ^{(r)})})^{-1}$.

The vector $\bar{v}^{(j)}_i := \frac{(n_j -1)}{n_j}\nu^{(j)}_i + \frac{1}{n_j}e_{d+j}$  (resp. $\epsilon_i$)  belongs to $\s_j^+ \cap \Z^{d+g+1}$ (resp. $\s_j^- \cap \Z^{d+g+1}$).
We denote by $\Theta'$ the minimal subdivision of the fan $\Theta$ which contains the cones 
$\tau^+_j :=\r_j + \bar{v}^{(j)}_i \R_{\geq 0} \subset  \s_j^+$ and  $\tau^-_j :=\r_j + e_i \R_{\geq 0}  \subset \s_j^- $.

Notice that  $\xi_j(\bar{v}^{(j)}_i ) = n_j$ and $\eta(\bar{v}^{(j)}_i )= (n_j -1)n_{j+1} \cdots n_g$ and also $\xi_j (\epsilon_i) = 1$ and $\eta(\epsilon_i) =0$. 

It is easy to see that  $\bar{v}^{(j)}_i$ (resp. $e_i$) together with the elements of a basis of the lattice $\f_j ^* (N_j) $
are part of a basis of $\Z^{d+g+1}$ hence we obtain that 
\[
S_{\t_j^+} = S_{\r_j}  \L ^{-\xi_j(\bar{v}^{(j)}_i )  } T^{\eta(\bar{v}^{(j)}_i ) }  (1 - \L ^{-\xi_j(\bar{v}^{(j)}_i )  } T^{\eta(\bar{v}^{(j)}_i ) } )^{-1} \mbox{ and } 
S_{\t_j^+} = S_{\r_j}  \L ^{-1}  (1 - \L ^{-1  } )^{-1}.
\]
We can decompose then the local contribution as a sum 
$LC_j = (\L- 1) ^{d+1} ( S_{\s_{j+1} ^- } + \sum_{\t}  S_\t ) +  (\L -1)^d (\L - 2) S_{\r_j} $, where 
$\t \ne \r_j $ runs through the cones of $\Theta'$ which meet the interior of $\s_j^{-}$ or of $\s_j^{+}$. 

The linear map $\Upsilon_{j+1} : N_{j+1} ' \to \Z^{d+g+1}$ was defined in 
Lemma \ref{minus}. One has that
$\Upsilon_{j+1}(0,1)=\epsilon_{d+1} + \sum_{i=j+2}^{g+1} n_{j+1} \cdots n_{i-1}\epsilon_{d+i}$. 
Hence
$\xi_{j+1} ( \Upsilon_{j+1} (0, 1) ) = 1$ and $\eta(\Upsilon_{j+1} (0, 1) ) = n_{j+1} \cdots n_g$.
By Lemma \ref{minus} we deduce that 
$S_{\Upsilon_{j+1}(\r ') }= S_{\r_j} \frac{\L^{-1}T^{n_{j+1} \cdots n_g}}{1 - \L^{-1}T^{n_{j+1} \cdots n_g}}$.

By these observations and the equality 
(\ref{minus-rel}) 
 it is enough to check that 
\begin{center}
$ 
\begin{array}{lcl}
LC_j '  & := &   (\L- 1) ^{d+1} ( S_{ \Upsilon_{j+1}  (\r') } + S_{\t_j^+}  
+ S_{\t_j^-}  ) +  (\L -1)^d (\L - 2) S_{\r_j} 
\\
& = & R_j'' \frac{\L ^{-n_j   } T^{(n_j -1)n_{j+1} \cdots n_g }}{  (1 - \L ^{-n_j } T^{ (n_j -1)n_{j+1} \cdots n_g } )} 
\frac{\L ^{- 1} T^{  n_{j+1} \cdots n_g } }{ (1 - \L ^{-1 } T^{ n_{j+1} \cdots n_g } )}
\prod_{l=1, l \ne i} \frac{ L^{- \xi_r (\nu_l^{(r)} ) } 
T^{ \eta( \nu_l ^{(r)})} }{ (1 - L^{- \xi_r (\nu_l^{(r)} ) } T^{ \eta( \nu_l ^{(r)})})} 
\end{array}$, 
\end{center}
for some polynomial $R_j'' \in \Z[\L^{\pm 1}][T]$.
The term $ LC_j ' $ is equal to 
\begin{center}
$
\begin{array}{c} 
(\L-1)^{d}  \big ( S_{\r_j}  + (\L-1)S_{\r_j} 
\frac{\L^{-n_j}T^{(n_j-1)n_{j+1} \cdots n_g}}{1 - {\L^{-n_j}T^{(n_j-1)n_{j+1} \cdots n_g}}} + 
(\L-1)S_{\r_j}\frac{\L^{-1}T^{n_j \cdots n_g}}{1-\L^{-1}T^{n_{j+1} \cdots n_g}} + (\L - 2) S_{\r_j} \big )
\\
 = (\L-1)^{d+1}S_{\r_j}  \big ( 1 + \frac{\L^{-n_j}T^{(n_j-1)n_{j+1} 
\cdots n_g}}{1 - {\L^{-n_j}T^{(n_j-1)n_{j+1} \cdots n_g}}} + 
\frac{\L^{-1}T^{n_{j+1} \cdots n_g}}{1-\L^{-1}T^{n_{j+1} \cdots n_g}} \big )
\\
 =  
(\L-1)^{d+1}S_{\r_j}  \big ( \frac{1 - \L^{-(n_j+1)}T^{n_j 
\cdots n_g}}{(1 - \L^{-n_j}T^{(n_j-1)n_{j+1} \cdots n_g})(1-\L^{-1}T^{n_{j+1} \cdots n_g})} \big ). 
\end{array}$
\end{center}
Finally, notice that the greatest common divisor of the terms $1 - \L^{ -(n_j +1)} T^{n_j \cdots n_g}$ and 
  $1 - \L^{-n_j}T^{(n_j-1)n_{j+1} \cdots n_g}$ (resp. of $1 - \L^{ -(n_j +1)} T^{n_j \cdots n_g}$  and
 $1-\L^{-1}T^{n_{j+1} \cdots n_g} $)
is equal to one. 

The proof in the case $j = g$ follows by a similar argument. 
\hfill $\square$

\bigskip

\noindent
\textbf{Proof of Corollary \ref{poles}}.
We use the formula for the naive motivic zeta function
of Theorem \ref{main-zeta}. 
By definition if $1 \leq i \leq d$ there exists at most one integer $1 \leq j \leq g$ 
such that $\nu_i^{(j)}$ is special. 
If  $1 \leq i, l \leq d$, $i \ne l$ and the vectors $\nu_i^{(j)}$ and 
$\nu_l^{(r)}$ are special for $1 \leq j , r \leq g$ then $j \ne r$. 
In this case 
then the pairs $(\xi_j (\nu_i^{(j)} ), \eta(\nu_i^{(j)}) = (n_j +1, n_j \cdots n_g ) $ and 
$(\xi_{r} (\nu_l^{(r)}) , \eta (\nu_l^{(r)}) = (n_r +1, n_r \cdots n_g)$ are linearly independent, hence 
the greatest common divisor of $1 - \L^{ -n_j -1} T^{n_j \cdots n_g}$ and 
$1 - \L^{ -n_l -1} T^{n_l \cdots n_g}$ is equal to one. 
These conditions guarantee that we can apply Proposition \ref{LC}  any time 
we have special vectors. 

Then the proof for the topological zeta function follows by applying Proposition \ref{ztop-rel}.  
\hfill $\square$

\begin{remark}  \label{comp-4} By the results obtained in this section 
it is easy to see that   the set of \textit{strongly  candidate poles} 
for the motivic zeta function introduced in \cite{ACNLM-AMS} 
 is equal to 
\[ 
SCP(f,w) := \cup_{j=1}^{g}  \cup_{i=1}^d \{ (\eta ( \nu_i^{(j)} ) , 
\xi_j (\nu_i^{(j)} ) \mid \nu_i^{(j)}  \textrm{ non-special}  \}  \cup \{ (1, 1) \}.\]
The set $SCP(f,w)$ of candidate poles  coincides with the one we get from  Corollary \ref{poles}. 
This implies that the arguments in the  proof of the monodromy conjecture in \cite{ACNLM-AMS} Chapter 6, 
are not affected by 
the inaccuracy of the  formulas for the zeta functions in \cite{ACNLM-AMS} 
(see Remarks \ref{comp} and \ref{comp-1}). 
\end{remark}

\section{Example}

\begin{example}   \label{false}
The q.o. polynomial $f=(z^2 - xy^3)^4 - x^4y^{13}$  is and has 
characteristic exponents $\lambda_1=(1/2,3/2)$ and $\lambda_2=(1/2,7/4)$.
 The generators of the semigroup are $\gamma_1=(1/2,3/2)$ and $\gamma_2=(1, 13/4)$.
We have that $n_1 = 2$ and $n_2 = 4$,  $\ell_1=\ell_2=2$ and $r_1=r_2=1$.
The conic integral complex $\Theta$ associated to $f$ is
represented by the Figure \ref{figexample8.6}. We have that
 \[ c_1
(\s_1^+) = c_1 (\s_2 ^- ) = [\mu_4], \quad c_1 (\s_1^-) ) [\mu_8],
\quad c_1( \s_2^+ ) = c_1 (\s_3) =1\] and also 
$c_1( \r_1) = [\{(x,y) \in (\C^*)^2 | (y^2 - x)^4=1\}$, $c_1
(\r_2) = [\{(x,{y}) \in (\C^*)^2 | y^4 - x=1\}$.
 The multiplicities of the simplicial cones are 
$ \mathrm{mult} (\r_1 ) =  \mathrm{mult} (\s_2^-)  = 2$, $ \mathrm{mult} (\s_1^+ ) = 
\mathrm{mult} (\s_2^+)  = 4$,  
  and $\mathrm{mult} (\r_2 ) =  \mathrm{mult} (\s_3)  = 1$.

The motivic Milnor fibre is
$$S_{f,0} = c_1(\r_1) (1-\L) + [\mu_8]\L + c_1 (\r_2) (1-\L) - [\mu_4]\L(1-\L) + (1 - \L)^2.$$
The Hodge-Steenbrink spectrum is
$hsp(f,0)=\frac{1-t}{1-t^{1/8}}t - t = t^{\frac{9}{8}} + \cdots + t^{\frac{15}{8}} $.
According to Theorem \ref{main-zeta} the topological zeta function is 
\begin{center}
$Z_{\mathrm{top}}(f,s) =
\frac{24 s + 13}{(3+8s)(5+24s)} + \frac{22 + 96s}{(5+24s)(3+8s)(11+52s)}
- \frac{s}{(1+s)(3+8s)(11+52s)}.$ 
\end{center}

We compute this example with the 
method of \cite{ACNLM-AMS} with respect to the form $\omega=dx \wedge dy \wedge dz$.  
The Newton map followed by a suitable translation $\mathcal{N}_1$ is given by
 $x=x_1^2$, $y = y_1^2$, $z=
x_1y_1^3(z_1-1)$.
We get  $f\circ \mathcal{N}_1=x_1^8y_1^24(z_1^4-y_1^2 + \cdots) $  and 
$\omega \circ \mathcal{N}_1 =x_1^2y_1^4dx_1\wedge dy_1\wedge dz_1$.

Theorems 2.7 and 5.3 of \cite{ACNLM-AMS} imply that
$Z_{\mathrm{top}}(f,s) = \frac{24 s + 13}{(3+8s)(5+24s)} + Z_B (s) $, 
where   $\frac{24 s + 13}{(3+8s)(5+24s)}$ is called the Part A in \cite[Definition 2.3]{ACNLM-AMS}.  
The term 
 $ Z_B (s) $ is the sum of three contributions
corresponding to the three faces of the compact edge of the Newton polyhedron of 
$f \circ \mathcal{N}_1$.   The contribution of the two zero-dimensional faces is equal 
to $J_{\s^+_2}(f,s) + J_{\s^-_2}(f,s) = \frac{22 + 96s}{(5+24s)(3+8s)(11+52s)}$. 
The contribution of this edge equals
$\frac{-2s}{(1+s)(3+8s)(11+52s)}$ because it has integral lenght two. 

The disagreement with our
formula appears since the Newton map $\mathcal{N}_1$ is not birational and
the strict transform of $f$ by the Newton map is not irreducible
(see Remark \ref{false2}).

    \

\setlength{\unitlength}{0.5 mm}
\begin{figure}[h]
\begin{center}
\begin{picture}(-60,40)(-6,-1)
\linethickness{0.4mm}

\drawline(-85,-20)(-25,-20)

\drawline(-25,-20)(-37,5)
\drawline(-85,-20)(-77,10)
\drawline(-37,5)(-60,40)
\drawline(-77,10)(-60,40)
\drawline(-77,10)(20,40)

\drawline(-77,10)(-37,5)
\drawline(-77,10)(-35,23)
\drawline(-77,10)(-23,40)
\drawline(-77,10)(-18,23)

\drawline(-18,23)(-37,5)
\drawline(-18,23)(-23,40)
\drawline(-18,23)(20,40)

\jput(-85,-20){\circle*{2}}
\jput(-25,-20){\circle*{2}}
\jput(-60,40){\circle*{2}}
\jput(-77,10){\circle*{2}}
\jput(-37,5){\circle*{2}}
\jput(-18,23){\circle*{2}}
\jput(20,40){\circle*{2}}
\jput(-23,40){\circle*{2}}

\jput(-93,-20){$\epsilon_1$}
\jput(-20,-20){$\epsilon_2$}
\jput(-65,43){$\epsilon_3$}

\jput(-115,3){$(2,0,1,2,8)$}
\jput(-30,3){$(0,2,3,6,24)$}
\jput(-56,-8){${\sigma}^-_1$}
\jput(-61,24){${\sigma}^+_1$}
\jput(-39,13){${\sigma}^-_2$}

\jput(-29,29){${\sigma}^+_2$}
\jput(17,43){$\epsilon_5$}
\jput(-25,43){$\epsilon_4$}
\jput(-17,18){$(0,4,6,13,52)$}
\end{picture}
\end{center}

\

\

\caption[]{Projectivization of the fan $\Theta \subset
{\Z}^5_{\geq 0}$ of $f=(z^2-xy^3)^4 +
x^4y^{13}$.}\label{figexample8.6}
\end{figure}
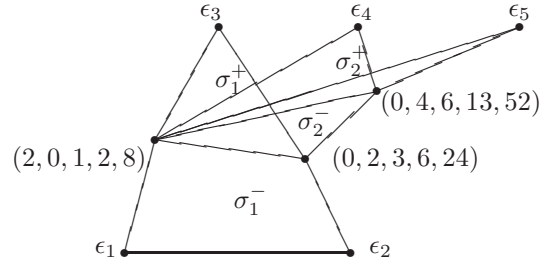
\end{example}

 \vspace{0.5cm}

\bibliographystyle{amsplain}

\begin{thebibliography}{10}

\bibitem{Abhyankar}
S.~S.~Abhyankar, \emph{On the ramification of algebraic functions}, Amer. J.
  Math. \textbf{77} (1955), 575--592. 


\bibitem{ACNLM-JLMS}
E.~Artal~Bartolo, Pi. Cassou-Nogu{\`e}s, I.~Luengo, and A.~Melle~Hern{\'a}ndez,
  \emph{The {D}enef-{L}oeser zeta function is not a topological invariant}, J.
  London Math. Soc. (2) \textbf{65} (2002), no.~1, 45--54.

\bibitem{ACNLM-AMS}
\bysame, \emph{Quasi-ordinary power series and their
  zeta functions}, Mem. Amer. Math. Soc. \textbf{178} (2005), no.~841. 
  



\bibitem{ACNLM-CM-08}
   \bysame,
  \emph{On the log-canonical threshold for germs of plane curves},
  Singularities {I}, Contemp. Math., vol. 474, Amer. Math. Soc., Providence,
  RI, 2008, pp.~1--14. 


\bibitem{DL-JAMS}
J.~Denef and F.~Loeser, \emph{Caract\'eristiques d'{E}uler-{P}oincar\'e,
  fonctions z\^eta locales et modifications analytiques}, J. Amer. Math. Soc.
  \textbf{5} (1992), no.~4, 705--720. 

\bibitem{Denef-LoeserIgusa}
\bysame, \emph{Motivic {I}gusa zeta functions}, J.
  Algebraic Geom. \textbf{7} (1998), no.~3, 505--537. 

\bibitem{Denef-LoeserInvent}
\bysame, \emph{Germs of arcs on singular algebraic varieties and motivic
  integration}, Invent. Math. \textbf{135} (1999), no.~1, 201--232.

\bibitem{Denef-LoeserBarca}
\bysame, \emph{Geometry on arc spaces of algebraic varieties}, European
  {C}ongress of {M}athematics, {V}ol. {I} ({B}arcelona, 2000), Progr. Math.,
  vol. 201, Birkh\"auser, Basel, 2001, pp.~327--348. 

\bibitem{Denef-LoeserLef}
\bysame, \emph{Lefschetz numbers of iterates of the monodromy and truncated
  arcs}, Topology \textbf{41} (2002), no.~5, 1031--1040. 

\bibitem{Denef-LoeserMcKay}
\bysame, \emph{Motivic integration, quotient singularities and the {M}c{K}ay
  correspondence}, Compositio Math. \textbf{131} (2002), no.~3, 267--290.


\bibitem{DuBois-Michel}
P.~Du~Bois and F.~Michel, \emph{The integral {S}eifert form
  does not determine the topology of plane curve germs}, J. Algebraic Geom.
  \textbf{3} (1994), no.~1, 1--38. 

\bibitem{Durfee}
A.~H.~Durfee, \emph{Fibered knots and algebraic singularities}, Topology
  \textbf{13} (1974), 47--59. 

\bibitem{Ewald}
G.~Ewald, \emph{Combinatorial convexity and algebraic geometry},
  Graduate Texts in Mathematics, vol. 168, Springer-Verlag, New York, 1996.

\bibitem{Fulton}
W.~Fulton, \emph{Introduction to toric varieties}, Annals of Mathematics
  Studies, vol. 131, Princeton University Press, Princeton, NJ, 1993, The
  William H. Roever Lectures in Geometry. 

\bibitem{Gau}
Y.-N.~Gau, \emph{Embedded topological classification of quasi-ordinary
  singularities}, Mem. Amer. Math. Soc. \textbf{74} (1988), no.~388, 109--129,
  With an appendix by Joseph Lipman. 

\bibitem{GP-Can}
P.~D. ~Gonz{\'a}lez~P{\'e}rez, \emph{Singularit\'es quasi-ordinaires toriques et
  poly\`edre de {N}ewton du discriminant}, Canad. J. Math. \textbf{52} (2000),
  no.~2, 348--368. 

\bibitem{GP-Fourier}
\bysame, \emph{Toric embedded resolutions of
  quasi-ordinary hypersurface singularities}, Ann. Inst. Fourier (Grenoble)
  \textbf{53} (2003), no.~6, 1819--1881. 

\bibitem{GPSemi}
\bysame, \emph{The semigroup of a quasi-ordinary
  hypersurface}, J. Inst. Math. Jussieu \textbf{2} (2003), no.~3, 383--399.

\bibitem{MGV}
M.~Gonz\'alez~Villa, \emph{Funci\'on zeta mot\'\i vica de singularidades
  quasiordinarias irreducibles}, Tesis Doctoral, Universidad Complutense de
  Madrid, 2010. 

\bibitem{Guibert}
G.~Guibert, \emph{Espaces d'arcs et invariants d'{A}lexander}, Comment. Math.
  Helv. \textbf{77} (2002), no.~4, 783--820. 



\bibitem{GLM-Duke}
G.~Guibert, F.~Loeser, and M.~Merle, \emph{Iterated vanishing
  cycles, convolution, and a motivic analogue of a conjecture of {S}teenbrink},
  Duke Math. J. \textbf{132} (2006), no.~3, 409--457. 

\bibitem{Jung}
H.W.E.~Jung, \emph{Darstellung der funktionen eines algebraischen k\"orpers
  zweier unabh\"angiger ver\"anderlichen $x$, $y$ in der umgebung einer stelle
  $x = a$, $y = b$}, J. Reine Angew. Math. \textbf{133} (1908), 289--314.

\bibitem{Kempf}
G.~Kempf, F.~F. ~Knudsen, D.~Mumford, and B.~Saint-Donat, \emph{Toroidal
  embeddings. {I}}, Lecture Notes in Mathematics, Vol. 339, Springer-Verlag,
  Berlin, 1973. 


\bibitem{Kulikov}
V.~S. Kulikov, \emph{Mixed {H}odge structures and singularities},
  Cambridge Tracts in Mathematics, vol. 132, Cambridge University Press,
  Cambridge, 1998. 


\bibitem{MR0344248}
A.~Landman, \emph{On the {P}icard-{L}efschetz transformation for algebraic
  manifolds acquiring general singularities}, Trans. Amer. Math. Soc.
  \textbf{181} (1973), 89--126. 




\bibitem{Lipman}
J.~ Lipman, \emph{Topological invariants of quasi-ordinary singularities},
  Mem. Amer. Math. Soc. \textbf{74} (1988), no.~388, 1--107. 

\bibitem{Lipman-Eq}
\bysame, \emph{Equisingularity and simultaneous resolution of singularities},
  Resolution of singularities ({O}bergurgl, 1997), Progr. Math., vol. 181,
  Birkh\"auser, Basel, 2000, pp.~485--505. 

\bibitem{McEwan-Nemethi}
L.~J.~McEwan and A.~N{\'e}methi, \emph{The zeta function of a
  quasi-ordinary singularity}, Compos. Math. \textbf{140} (2004), no.~3,
  667--682. 

\bibitem{Milnor}
J.~Milnor, \emph{Singular points of complex hypersurfaces}, Annals of
  Mathematics Studies, No. 61, Princeton University Press, Princeton, N.J.,
  1968. 



\bibitem{Navarro}
V.~Navarro~Aznar, \emph{Sur la th\'eorie de {H}odge-{D}eligne}, Invent. Math.
  \textbf{90} (1987), no.~1, 11--76. 

\bibitem{Oda}
T.~Oda, \emph{Convex bodies and algebraic geometry}, Ergebnisse der
  Mathematik und ihrer Grenzgebiete (3) [Results in Mathematics and Related
  Areas (3)], vol.~15, Springer-Verlag, Berlin, 1988

\bibitem{Chris-Steenbrink}
C.~A.~M.~Peters and J.~H.~M.~Steenbrink, \emph{Mixed {H}odge
  structures}, Ergebnisse der Mathematik und ihrer Grenzgebiete. 3. Folge. Vol.~52,
  Springer-Verlag, Berlin, 2008. 

\bibitem{PPP-Duke}
P.~Popescu-Pampu, \emph{On the analytical invariance of the semigroups of
  a quasi-ordinary hypersurface singularity}, Duke Math. J. \textbf{124}
  (2004), no.~1, 67--104. 


\bibitem{Saito-MHM}
M.~Saito,  \emph{Mixed {H}odge modules}, Publ. Res. Inst. Math. Sci.
  \textbf{26} (1990), no.~2, 221--333. 


\bibitem{Saito-MA91}
\bysame, \emph{On {S}teenbrink's conjecture}, Math. Ann. \textbf{289} (1991),
  no.~4, 703--716. 


\bibitem{Saito00}
\bysame, \emph{Exponents of an irreducible plane curve singularity},
  arXiv:math/0009133v2 (2000). 


\bibitem{Sakamoto}
K.~Sakamoto, \emph{The {S}eifert matrices of {M}ilnor fiberings defined by
  holomorphic functions}, J. Math. Soc. Japan \textbf{26} (1974), 714--721.


\bibitem{SSS}
R. Schrauwen,  J.~H.~M.~Steenbrink, and J.~Stevens, \emph{Spectral pairs and the
  topology of curve singularities}, Complex geometry and {L}ie theory
  ({S}undance, {UT}, 1989), Proc. Sympos. Pure Math., vol.~53, Amer. Math.
  Soc., Providence, RI, 1991, pp.~305--328. 

\bibitem{stanley}
R.~P.~Stanley, \emph{Enumerative combinatorics. {V}ol. 1}, Cambridge
  Studies in Advanced Mathematics, vol.~49, Cambridge University Press,
  Cambridge, 1997. 

\bibitem{Steenbrink76}
J.~H.~M. Steenbrink, \emph{Mixed {H}odge structure on the vanishing
  cohomology}, Real and complex singularities,
 Sijthoff and
  Noordhoff, Alphen aan den Rijn, 1977, pp.~525--563. 

\bibitem{Steenbrink-Asterisque}
\bysame, \emph{The spectrum of hypersurface singularities}, Ast\'erisque
  (1989), no.~179-180, 11, 163--184.



\bibitem{Varchenko}
A.~N. Varchenko, \emph{Asymptotic {H}odge structure on vanishing
  cohomology}, Izv. Akad. Nauk SSSR Ser. Mat. \textbf{45} (1981), no.~3,
  540--591, 688. 
\end{thebibliography}

\providecommand{\bysame}{\leavevmode\hbox to3em{\hrulefill}\thinspace}
\providecommand{\MR}{\relax\ifhmode\unskip\space\fi MR }
\providecommand{\MRhref}[2]{%
  \href{http://www.ams.org/mathscinet-getitem?mr=#1}{#2}
}
\providecommand{\href}[2]{#2}

\end{document}